\documentclass[12pt,twoside,leqno]{article}
\usepackage{amsmath}
\usepackage{amssymb}
\usepackage{amsxtra}
\usepackage{amscd}
\usepackage{amsthm}
\usepackage[mathscr]{eucal}

\theoremstyle{plain}
\newtheorem{thm}[subsection]{Theorem}
\newtheorem{prop}[subsection]{Proposition}
\newtheorem{cor}[subsection]{Corollary}
\newtheorem{lem}[subsection]{Lemma}

\theoremstyle{definition}

\newtheorem{rem}[subsection]{Remark}
\newtheorem{para}[subsection]{}

\newenvironment{pf}{\proof[\proofname]}{\endproof}

\begin{document}

\title{Attempts on SGA  for non-commutative rings. By Kazuya Kato. \newline
With Appendix by Takako Fukaya and Kazuya Kato.}


\maketitle

\newcommand\Cal{\mathcal}
\newcommand\define{\newcommand}

\define\gp{\mathrm{gp}}%
\define\fs{\mathrm{fs}}%
\define\an{\mathrm{an}}%
\define\mult{\mathrm{mult}}%
\define\Ker{\mathrm{Ker}\,}%
\define\Coker{\mathrm{Coker}\,}%
\define\Hom{\mathrm{Hom}\,}%
\define\Ext{\mathrm{Ext}\,}%
\define\rank{\mathrm{rank}\,}%
\define\gr{\mathrm{gr}}%
\define\cHom{\Cal{Hom}}
\define\cExt{\Cal Ext\,}%

\define\cB{\Cal B}
\define\cC{\Cal C}
\define\cD{\Cal D}
\define\cO{\Cal O}
\define\cS{\Cal S}
\define\cM{\Cal M}
\define\cG{\Cal G}
\define\cH{\Cal H}
\define\cI{\Cal I}
\define\cE{\Cal E}
\define\cU{\Cal U}
\define\cF{\Cal F}
\define\cN{\Cal N}
\define\fF{\frak F}
\define\Dc{\check{D}}
\define\Ec{\check{E}}

\newcommand{\cP}{{\Cal {P}}}
\newcommand{\cV}{{\Cal {V}}}
\newcommand{\cZ}{{\Cal {Z}}}
\newcommand{\N}{{\mathbb{N}}}
\newcommand{\Q}{{\mathbb{Q}}}
\newcommand{\Z}{{\mathbb{Z}}}
\newcommand{\R}{{\mathbb{R}}}
\newcommand{\C}{{\mathbb{C}}}
\newcommand{\bN}{{\mathbb{N}}}
\newcommand{\bQ}{{\mathbb{Q}}}
\newcommand{\bF}{{\mathbb{F}}}
\newcommand{\bZ}{{\mathbb{Z}}}
\newcommand{\bP}{{\mathbb{P}}}
\newcommand{\bR}{{\mathbb{R}}}
\newcommand{\bC}{{\mathbb{C}}}
\newcommand{\bbQ}{{\bar \mathbb{Q}}}
\newcommand{\ol}[1]{\overline{#1}}
\newcommand{\too}{\longrightarrow}
\newcommand{\respect}{\rightsquigarrow}
\newcommand{\compatible}{\leftrightsquigarrow}
\newcommand{\upc}[1]{\overset {\lower 0.3ex \hbox{${\;}_{\circ}$}}{#1}}
\newcommand{\Gmlog}{\bG_{m, \log}}
\newcommand{\Gm}{\bG_m}
\newcommand{\ep}{\varepsilon}
\newcommand{\Spec}{\operatorname{Spec}}
\newcommand{\val}{{\mathrm{val}}} 
\newcommand{\n}{\operatorname{naive}}
\newcommand{\bs}{\operatorname{\backslash}}
\newcommand{\Gal}{\operatorname{{Gal}}}
\newcommand{\gal}{{\rm {Gal}}({\bar \Q}/{\Q})}
\newcommand{\galp}{{\rm {Gal}}({\bar \Q}_p/{\Q}_p)}
\newcommand{\gall}{{\rm{Gal}}({\bar \Q}_\ell/\Q_\ell)}
\newcommand{\wep}{W({\bar \Q}_p/\Q_p)}
\newcommand{\wel}{W({\bar \Q}_\ell/\Q_\ell)}
\newcommand{\Ad}{{\rm{Ad}}}
\newcommand{\BS}{{\rm {BS}}}
\newcommand{\even}{\operatorname{even}}
\newcommand{\End}{{\rm {End}}}
\newcommand{\odd}{\operatorname{odd}}
\newcommand{\GL}{\operatorname{GL}}
\newcommand{\np}{\text{non-$p$}}
\newcommand{\g}{{\gamma}}
\newcommand{\G}{{\Gamma}}
\newcommand{\Lam}{{\Lambda}}
\newcommand{\La}{{\Lambda}}
\newcommand{\lam}{{\lambda}}
\newcommand{\la}{{\lambda}}
\newcommand{\uL}{{{\hat {L}}^{\rm {ur}}}}
\newcommand{\uQp}{{{\hat \Q}_p}^{\text{ur}}}
\newcommand{\sel}{\operatorname{Sel}}
\newcommand{\dt}{{\rm{Det}}}
\newcommand{\Sig}{\Sigma}
\newcommand{\fil}{{\rm{fil}}}
\newcommand{\SL}{{\rm{SL}}}
\newcommand{\spl}{{\rm{spl}}}
\newcommand{\st}{{\rm{st}}}
\newcommand{\Isom}{{\rm {Isom}}}
\newcommand{\Mor}{{\rm {Mor}}}
\newcommand{\bg}{\bar{g}}
\newcommand{\id}{{\rm {id}}}
\newcommand{\cone}{{\rm {cone}}}
\newcommand{\al}{a}
\newcommand{\ChL}{{\cal{C}}(\La)}
\newcommand{\Image}{{\rm {Image}}}
\newcommand{\toric}{{\operatorname{toric}}}
\newcommand{\torus}{{\operatorname{torus}}}
\newcommand{\Aut}{{\rm {Aut}}}
\newcommand{\Qp}{{\mathbb{Q}}_p}
\newcommand{\barQp}{{\mathbb{Q}}_p}
\newcommand{\Qpur}{{\mathbb{Q}}_p^{\rm {ur}}}
\newcommand{\Zp}{{\mathbb{Z}}_p}
\newcommand{\Zl}{{\mathbb{Z}}_l}
\newcommand{\Ql}{{\mathbb{Q}}_l}
\newcommand{\Qlur}{{\mathbb{Q}}_l^{\rm {ur}}}
\newcommand{\F}{{\mathbb{F}}}
\newcommand{\eps}{{\epsilon}}
\newcommand{\epsLa}{{\epsilon}_{\La}}
\newcommand{\epsLaVxi}{{\epsilon}_{\La}(V, \xi)}
\newcommand{\epsOLaVxi}{{\epsilon}_{0,\La}(V, \xi)}
\newcommand{\Qplin}{{\mathbb{Q}}_p(\mu_{l^{\infty}})}
\newcommand{\otimesQplin}{\otimes_{\Qp}{\mathbb{Q}}_p(\mu_{l^{\infty}})}
\newcommand{\galFl}{{\rm{Gal}}({\bar {\Bbb F}}_\ell/{\Bbb F}_\ell)}
\newcommand{\gallur}{{\rm{Gal}}({\bar \Q}_\ell/\Q_\ell^{\rm {ur}})}
\newcommand{\galFF}{{\rm {Gal}}(F_{\infty}/F)}
\newcommand{\galFv}{{\rm {Gal}}(\bar{F}_v/F_v)}
\newcommand{\galF}{{\rm {Gal}}(\bar{F}/F)}
\newcommand{\epsVxi}{{\epsilon}(V, \xi)}
\newcommand{\epsOVxi}{{\epsilon}_0(V, \xi)}
\newcommand{\plim}{\lim_
{\scriptstyle 
\longleftarrow \atop \scriptstyle n}}
\newcommand{\sig}{{\sigma}}
\newcommand{\ga}{{\gamma}}
\newcommand{\del}{{\delta}}
\newcommand{\Vss}{V^{\rm {ss}}}
\newcommand{\Bst}{B_{\rm {st}}}
\newcommand{\Dpst}{D_{\rm {pst}}}
\newcommand{\Dcrys}{D_{\rm {crys}}}
\newcommand{\DdR}{D_{\rm {dR}}}
\newcommand{\Fin}{F_{\infty}}
\newcommand{\Kla}{K_{\lambda}}
\newcommand{\Ola}{O_{\lambda}}
\newcommand{\Mla}{M_{\lambda}}
\newcommand{\Det}{{\rm{Det}}}
\newcommand{\Sym}{{\rm{Sym}}}
\newcommand{\LaSa}{{\La_{S^*}}}
\newcommand{\cX}{{\cal {X}}}
\newcommand{\MHG}{{\frak {M}}_H(G)}
\newcommand{\tauMla}{\tau(M_{\lambda})}
\newcommand{\Fvur}{{F_v^{\rm {ur}}}}
\newcommand{\Lie}{{\rm {Lie}}}
\newcommand{\cL}{{\cal {L}}}
\newcommand{\cW}{{\cal {W}}}
\newcommand{\fq}{{\frak {q}}}
\newcommand{\cont}{{\rm {cont}}}
\newcommand{\SC}{{SC}}
\newcommand{\Om}{{\Omega}}
\newcommand{\dR}{{\rm {dR}}}
\newcommand{\crys}{{\rm {crys}}}
\newcommand{\hatSig}{{\hat{\Sigma}}}
\newcommand{\rdet}{{{\rm {det}}}}
\newcommand{\ord}{{{\rm {ord}}}}
\newcommand{\BdR}{{B_{\rm {dR}}}}
\newcommand{\BdRO}{{B^0_{\rm {dR}}}}
\newcommand{\Bcrys}{{B_{\rm {crys}}}}
\newcommand{\Qw}{{\mathbb{Q}}_w}
\newcommand{\barkappa}{{\bar{\kappa}}}
\newcommand{\oppLa}{{\Lambda^{\circ}}}
\newcommand{\bG}{{\mathbb{G}}}
\newcommand{\Rep}{{{\rm Rep}}}
\newcommand{\red}{{{\rm red}}}
\newcommand{\nspl}{{{\rm nspl}}}
\newcommand{\Proj}{{{\rm Proj}}}
\newcommand{\Alg}{{{\rm Alg}}}
\newcommand{\Pic}{{{\rm Pic}}}
\newcommand{\Sets}{{{\rm Sets}}}
\newcommand{\cl}{{{\rm cl}}}
\newcommand{\ct}{{{\rm cont}}}
\newcommand{\m}{{{\frak m}}}
\newcommand{\nn}{{{\frak n}}}
\newcommand{\fin}{{{\rm fin}}}

\begin{abstract} 
We define non-commutative schemes by using prime ideals  of non-commutative rings, and discuss the \'etale cohomology, the  Betti cohomology, and the fundamental groups of  non-commutative schemes. For non-commutative schemes which are finite over centers, we prove the finiteness theorem for the higher direct images in \'etale cohomology theory, and the comparison theorem between \'etale cohomology and Betti cohomology.  In Appendix,   for non-commutative schemes over finite fields which are finite over centers and satisfy a certain condition, $L$-functions are expressed by using \'etale cohomology with compact supports. 

\end{abstract}

\section*{Contents}

\noindent \S\ref{sec0}. Introduction

\noindent \S\ref{sec1}. Algebras over non-commutative rings

\noindent \S\ref{sec2}. Non-commutative schemes

\noindent \S\ref{sec3}. Basics of non-commutative schemes, 1

\noindent \S\ref{sec4}. Basics of non-commutative schemes, 2 

\noindent \S\ref{sec5}. \'Etale cohomology

\noindent \S\ref{sec6}.  Betti cohomology

\noindent \S\ref{sec7}. Fundamental groups

\noindent \S\ref{sec8}. Chow groups, relation to  class field theory

\noindent \S\ref{App}. Appendix: Zeta and $L$-functions of non-commutative schemes By Takako Fukaya

\setcounter{section}{-1}
\section{Introduction}\label{sec0}

\begin{para}  For a ring $A$ which need not be commutative, let  $\Spec(A)$ be the set of all prime ideals of $A$. 
Then for ring homomorphism $h: A\to B$ and $\frak p \in \Spec(B)$, $h^{-1}(\frak p)$ need not be a prime ideal of $A$, unlike the case of commutative rings.
This is probably the reason of the fact that 
in non-commutative algebraic geometry, people do not use the space $\Spec(A)$ so much. 
However, Procesi \cite{Pr} shows that if $h$ satisfies the condition

\medskip

(*)  As a ring,  $B$ is generated by $h(A)$ and the centralizer $C_B(A)= \{b\in B \;|\; h(a)b=bh(a)\; \text{for all $a\in A$}\}$ of $A$ in $B$,

\medskip
\noindent
then we have a map $$\Spec(B)\to \Spec(A)\;;\; \frak p \mapsto h^{-1}(\frak p).$$

We say $B$ is an $A$-algebra if this condition (*) is satisfied. We will show that if $B$ and $C$ are $A$-algebras, the tensor product $B\otimes_A C$ has a unique ring structure such that we have $(b_1\otimes  c_1)(b_2\otimes c_2)= (b_1b_2)\otimes (c_1c_2)$  if $b_1, b_2\in B$ and $c_1, c_2\in C$  and if  either $b_1, b_2\in C_B(A)$ or $c_1,c_2\in C_C(A)$ (Theorem. \ref{prod}). Many things about  algebras over commutative rings are generalized to $A$-algebras in this sense.

Consider the category ($\Alg$) of rings in which a morphism is a ring homomorphism $A\to B$ such that $B$ is an $A$-algebra. Using prime ideals of rings, we will define a category of non-commutative schemes which are endowed with sheaves of (not necessarily commutative) rings, and we will have a contra-variant functor $$A \mapsto \Spec(A)$$   from the full subcategory of $(\Alg)$ consisting of $A$ which is finitely generated as an algebra over the center $Z(A)$, to the category of non-commutative schemes. See Section \ref{sec2}.

\end{para}

\begin{para} 
For a non-commutative scheme $X$, 
 we define in Section \ref{sec5} the \'etale cohomology groups of $X$ and in Section \ref{sec7} the fundamental group of $X$, imitating SGA 4 and SGA 1, respectively, generalizing these groups for schemes. 
 The main results of this paper are the following (1)--(5). 
 
 We consider a condition ({\bf F}) on a non-commutative scheme and a condition ({\bf F}$_\C$) on a non-commutative scheme over $\C$  (cf. \ref{bfF}). 
 
 (i) If $R$ is an excellent commutative ring (resp. finitely generated commutative $\C$-algebra) and $A$ is an $R$-algebra such that $A$ is finitely generated as an $R$-module, every open subspace $U$ of $\Spec(A)$ satisfies ({\bf F}) (resp. ({\bf F}$_\C$)). 
 
 (ii) For a  non-commutative scheme $X$ (resp. $X$ over $\C$), the condition  ({\bf F}) (resp. ({\bf F}$_\C$)) is that $X$ has an open covering whose every member is isomorphic to some $U$ in (i). 
 
 The following (1) and (2) are the  non-commutative versions of the finiteness theorem of Gabber ((\cite{ILO} Exp. XIII, Exp. XXI) and the comparison theorem of Artin ((\cite{SGA4} vol. 3, Exp. XVI, 4)), and in fact proved by the reductions to their theorems. 
 
 \medskip
 
 (1) The finiteness theorem.
 
{\bf Theorem} (Theorem \ref{fin2}). Let $X$and $Y$ be quasi-compact non-commutative schemes satisfying ({\bf F}) and let $f:X\to Y$ be a morphism which is locally of finite presentation (\ref{lfp1}). Let $\cF$ be a constructible sheaf (\ref{const1}) of abelian groups (resp. sets, resp. groups) on the \'etale site of $X$. Then  $Rf_*\cF$ (resp. $f_*\cF$, resp. $R^1f_*\cF$)  is constructible. 

\medskip

(2) Comparison theorem. 

For a non-commutative scheme $X$ over $\C$ satisfying ({\bf F}$_\C$), we define a topological space $X_{cl}$ (\ref{cl1}). In the case $X$ is a scheme of finite type over $\C$, $X_{cl}$ is the set $X(\C)$ of $\C$-valued points of $X$ endowed with the classical topology. 

{\bf Theorem} (Theorem \ref{e=B}). Let $X$ and $Y$ be quasi-compact non-commutative schemes over $\C$ satisfying ({\bf F}$_\C$) and let $f:X\to Y$ be a morphism which is locally of finite presentation. Let $\cF$ be a constructible sheaf of abelian groups (resp. sets, resp. groups) on  the \'etale site of $X$. Then we have an isomorphism 
$$(R^mf_*\cF)_{cl} \overset{\cong}\to R^mf_{cl,*}\cF_{cl}$$
for every $m$ (resp. for $m=0$, resp. for $m=1$), where $f_{cl}$ denotes the continuous map $X_{cl}\to Y_{cl}$ induced by $f$ and $(\;)_{cl}$ for  sheaves denote the pullbacks to these spaces. 

\medskip
As a corollary of this theorem, for non-commutative schemes over $\C$ satisfying ({\bf F}$_\C$), 
we have a comparison theorem (Theorem \ref{e=B3}) for fundamental groups.

 \medskip
 
 (3) Proper base change theorem.
 
 We also prove a kind of proper base change theorem (Theorem \ref{pbc}, Theorem \ref{pbc2}) for non-commutative schemes satisfying ({\bf F}) and for morphisms satisfying a certain condition (\ref{indep}.1). 

\medskip

(4) Chow groups.

 In Section \ref{sec8}, we consider Chow groups of non-commutative schemes satisfying ({\bf F}), and consider a relation of $CH_0$ with class field theory (Theorem \ref{CFT}). 

\medskip

(5) Zeta and $L$-functions.

In Appendix with Takako Fukaya, basing on the study of Hasse zeta functions of non-commutative rings \cite{Fu}, the theory of zeta functions and $L$-functions for non-commutative schemes is given.  It is shown that for a non-commutative scheme $X$ over a finite field satisfying ({\bf F}) and a certain condition, the zeta function and $L$-functions of $X$ are expressed by using the  compact support \'etale cohomology (Theorem \ref{Hczeta}). This should be a fragment of SGA 5 for non-commutative schemes. 

Discussion of the authors on the zeta and $L$-functions for non-commutative rings was the motivation of the study of this paper.

\medskip

Our hope is 
 to develop theories as in SGA (S\'eminaire de G\'eom\'etrie Alg\'ebrique du Bois Marie) for non-commutative schemes.
 However we have to remark a negative present aspect of our work. Strong results  are obtained in this paper only for non-commutative schemes satisfying ({\bf F}). It seems that more ideas are needed to treat non-commutative rings with small centers nicely.

\end{para}
\begin{para} The authors thank Benson Farb and Nobushige Kurokawa for encouragements. 

The authors are partially supported by NSF Award 1601861.

\end{para}

We assume a ring has the identity  element  written by $1$ and ring homomorphisms respect $1$.

\section{Algebras over non-commutative rings}\label{sec1}

Let $A$ be a (not necessarily commutative) ring.
\begin{para}\label{Aalg1}

By an {\it $A$-algebra}, we mean a ring $B$ endowed with a ring homomorphism $A\to B$ satisfying the following equivalent conditions (i)--(iv). Let $C_B(A)$ be the centralizer $ \{b\in B\;|\; ab=ba\;\text{in $B$ for all} \; a\in A\}$ of $A$ in $B$.

(i) As a ring, $B$ is generated by the image of $A$ and $C_B(A)$. 

(ii) As a left $A$-module, $B$ is generated by $C_B(A)$.

(iii) As a right $A$-module, $B$ is generated by $C_B(A)$.

(iv) As a ring endowed with a ring homomorphism from $A$, $B$ is isomorphic to  $A\langle T_{\lambda} \; (\lambda\in \Lambda)\rangle/I$  for some index set $\Lambda$ and 
 some two-sided ideal $I$ of the non-commutative polynomial ring $A\langle T_{\lambda} \; (\lambda\in \Lambda)\rangle$. 
 Here $T_{\lambda}$ are indeterminates which do not commute with each other but commute with elements of $A$. 
  
\medskip

 If $A$ is commutative, this terminology $A$-algebra coincides with the usual terminology $A$-algebra.

\end{para}

\begin{para}\label{Aalg2}  By a {\it relatively commutative $A$-algebra} (r.c. $A$-algebra for short), we mean a ring $B$ endowed with a ring homomorphism $A\to B$ satisfying the following equivalent conditions (i)--(iv).

(i) As a ring, $B$ is generated by the image of $A$ and the center $Z(B)$ of $B$. 

(ii) As a left $A$-module, $B$ is generated by $Z(B)$.

(iii) As a right $A$-module, $B$ is generated by $Z(B)$. 

(iv) As a ring endowed with a ring homomorphism from $A$, $B$ is isomorphic to $A[T_{\la}\;;\;\la\in \La]/I$ for some index set $\La$ and for some two-sided ideal $I$ of $A[T_{\la}\;;\;\la\in \La]$. Here $T_{\la}$ are indeterminates which commute with each other 
and with elements of $A$. 

 \medskip
 
  If $A$ is commutative, an r.c. $A$-algebra is nothing but  a commutative $A$-algebra. 
 
 \end{para}

\begin{prop} (1) If $B$ is an $A$-algebra and $C$ is a $B$-algebra, $C$ is an   $A$-algebra.

(2) If $B$ is an r.c.  $A$-algebra and $C$ is an r.c. $B$-algebra, $C$ is an r.c.  $A$-algebra.
\end{prop}
\begin{pf} These are seen by using the condition (iv) in \ref{Aalg1} or \ref{Aalg2}.
\end{pf}

\begin{para} We denote by (Rings) the usual category of rings. We denote by (Alg) (resp. (Alg$^{r.c}$) the category of rings for which a morphism is a ring homomorphism $A\to B$ such that $B$ is an $A$-algebra (resp. an r.c. $A$-algebra). 

For a ring $A$, 
we denote by (Rings$/A$) the category of rings endowed with a ring homomorphism from $A$, in which a morphism is a ring homomorphism which commutes with the homomorphisms from $A$. Similarly, (Alg$/A$) (resp. (Alg$^{r.c}/A$)) denotes the category of objects over $A$ in the category  (Alg) (resp.  (Alg$^{r,c}$)).

\end{para}
 
 The following tensor product is related, through \ref{rcthm}, to the base change (\ref{lpr}, \ref{fibp}) in the category of non-commutative schemes. 
 \begin{prop}\label{prod} Let $B$ and $C$ be $A$-algebras. 
 
(1) There is a unique ring structure on $B\otimes_ A C$ such that $(b_1 \otimes c_1)(b_2\otimes c_2)=b_1b_2\otimes c_1c_2$ if $b_1, b_2\in B$ and $c_1, c_2\in C$ and if either $b_1,b_2\in C_B(A)$ or $c_1,c_2\in C_C(A)$.

(2) For an object $D$ of (Rings$/A$),  we have a bijection
 $$\Hom_{(\text{Rings}/A)}(B \otimes_A C, D) \to $$ $$= \{(f,g)\in \Hom_{(\text{Rings}/A)}(B, D) \times \Hom_{(\text{Rings}/A)}(C, D)\;|\; $$ $$ f(b)g(c)=g(c)f(b) \;\text{for all}\; b\in C_B(A), c\in C_C(A)\},$$
 where an element $h$ of the first set corresponds to an element  $(f,g)$ of the second set by
 $$f(b)=h(b \otimes 1), \; g(c)=h(1\otimes c), \quad h(b \otimes c)= f(b)h(c) \quad (b\in B, c\in C).$$
Note that in the definition of the second set, we can replace ``$b\in C_B(A), c\in C_C(A)$'' by ``$b\in B, c\in C_C(A)$'', and also by ``$b\in C_B(A), c\in C$'', not changing the second set. 

(3) 
 With the ring homomorphism  $B\to B\otimes_A C\;;\; b \mapsto b \otimes 1$, $B\otimes_A C$ is a $B$-algebra. 
  With the ring homomorphism $C\to B \otimes_A C\;;\; c \mapsto 1 \otimes c$, $B\otimes_A C$ is a $C$-algebra. 
  With the ring homomorphism $A \to B \otimes_A C\;;\; a \mapsto a\otimes 1=1\otimes a$, $B\otimes_A C$ is an $A$-algebra.

 \end{prop}
 
 \begin{pf} Take an isomorphism $B \cong A\langle T_{\lambda} \; (\lambda\in \Lambda)\rangle/I$  over $A$ as in (iv) of \ref{Aalg1}. This isomorphism induces  $B\otimes_A C\cong 
 C\langle T_{\lambda} \; (\lambda\in \Lambda)\rangle/IC$ where $IC$ denotes the right ideal of $C\langle T_{\lambda} \; (\lambda\in \Lambda)\rangle$ generated by the image of $I$. In the ring $C\langle T_{\lambda} \; (\lambda\in \Lambda)\rangle$, elements of $C_C(A)$ commute with all elements in the image of $B$ because they commute with $T_{\lam}$ and also with elements of $A$. Since $C$ is generated by $C_C(A)$ and $A$ as a ring, this shows that $IC$ is a two-sided ideal of $C\langle T_{\lambda} \; (\lambda\in \Lambda)\rangle$. Hence $B\otimes_A C$ has a ring structure. By the definition of this ring structure, we have $(b_1 \otimes c_1)(b_2 \otimes c_2)= b_1b_2\otimes c_1c_2$ for $b_1, b_2\in B$ and $c_1, c_2\in C_C(A)$. Since $C$ is generated by $C_C(A)$ and $A$ as a ring, this shows that $(b_1 \otimes c_1)(b_2 \otimes c_2)=b_1b_2\otimes c_1c_2$ for $b_1, b_2\in C_B(A)$ and $c_1,c_2\in C$. This proves (1). This also proves that we have the map $h\mapsto (f, g)$ from the first set to the second set in (2). If we have an element $(f, g)$ of the second set 
 in (2), the map $B\otimes_\Z C\to D\;;\; b\otimes c \mapsto f(b)g(c)$ induces a map $h:B \otimes_A C\to D$. The induced map $C\langle T_{\lambda} \; (\lambda\in \Lambda)\rangle\to D$ is a ring homomorphism because $f(b)g(c)=g(c)f(b)$ for all $b\in C_B(A)$ and $c\in C$. Hence we have (2). 
 
 (3) is easily deduced from (1), (2). 
  \end{pf}
   
  \begin{rem} We do not have  $(b_1\otimes c_1)(b_2 \otimes c_2)= (b_1b_2 \otimes c_1c_2)$ for $b_1,b_2\in B$ and $c_1, c_2\in C$ in general. If this holds in the case $A=B=C$, then for all $a_1, a_2\in A$, we have in $B\otimes_A C$
  $$a_1a_2\otimes 1=a_1\otimes a_2= (a_1\otimes 1)(1\otimes a_2) =(1\otimes a_1)(a_2\otimes 1)=a_2 \otimes a_1= a_2a_1 \otimes 1,$$ that is, $a_1a_2=a_2a_1$ in $A=A\otimes_A A$ and this would imply that $A$ is commutative. 
  
  \end{rem}

\begin{prop}\label{BC=CB} Let the notation be as in \ref{prod}.  Then there is a unique isomorphism $B \otimes_A C \cong C \otimes_A B$ of $A$-algebras which sends  $b\otimes c$ to $ c\otimes b$ for all $b\in C_B(A)$ and $c\in C_C(A)$. For $b\in B$ and $c\in C$, this sends $b\otimes c$ to $c\otimes b$ if either $b\in C_B(A)$ or $c\in C_C(A)$. 
\end{prop}

\begin{pf} By \ref{prod} (2), the functor (Rings$/A) \to$ (Sets)   represented by $B \otimes_A C$ and that represented by $C\otimes_A B$ are identified, and the isomorphism 
$B\otimes_A C \cong C\otimes_A B$ which gives the identification is described as in \ref{BC=CB}. 
\end{pf}

\begin{rem} The isomorphism in \ref{BC=CB} does not send  $b \otimes c$ to $ c\otimes b$ for $b\in B$ and $c\in C$ in general. If this holds in the case $A=B=C$, then for $a_1, a_2\in A$, the element $a_1a_2\otimes 1= a_1\otimes a_2$ of $B \otimes_A C$ would be sent to the element $a_2\otimes a_1=a_2a_1\otimes 1$ of $C\otimes_A B$, but
 since $B\otimes_A C=A=C\otimes_A B$, this would imply that $a_1a_2=a_2a_1$ in $A$ and hence that $A$ is commutative.

\end{rem}

\begin{prop}\label{ZZ} 

If $B$ is an $A$-algebra, the image of the center $Z(A)$ of $A$ in $B$ is contained in $Z(B)$.  

\end{prop}

\begin{pf} The image of an element of $Z(A)$ in $B$ commutes with elements of $C_B(A)$ and also with the images of all element of $A$. Hence it commutes with all elements of $B$. 
\end{pf}

\begin{cor} Let $B$ be an r.c. algebra over $A$ and let $C$ be an $A$-algebra. Let $h: B\to C$ be a morphism in (Rings$/A$). Then $h$ is a morphism in (Alg$/A$) if and only if $h(Z(B))\subset Z(C)$. 
\end{cor}

 \begin{cor}\label{A[T]}  Let $B$ be an $A$-algebra. Then for a two-sided ideal $I$ of $A[T_{\lam}\;;\;\lam\in \La]$, we have a bijection  $$\Hom_{(\text{Alg}/A)}(A[T_{\lam}\;;\;\lam\in \La]/I, B)$$ $$\to \{(x_{\lam})_{\lam}\in Z(B)^{\La}\;;\; f((x_{\lam})_{\lam})=0\;\text{for all}\; f\in I\}\;;\; h\mapsto (h(T_{\lam}))_{\lam}.$$
 
 \end{cor}

 \begin{prop}\label{rcthm} Let $B$ be an r.c. $A$-algebra and let  $C$ be an  $A$-algebra.
 
 (1) $C_B(A)=Z(B)$. 
 
 (2) The $C$-algebra $B \otimes_A C$ is  relatively commutative.

 (3) The $A$-algebra $B\otimes_A C$ is the push-out (i.e. co-fiber product) of $B \leftarrow A \rightarrow C$ in the category $(\Alg)$.

 (4) The category ($\Alg^{r.c.}$) has push-outs. If $B$ and $C$ are  r.c. $A$-algebras, $B\otimes_A C$ is the push-out of $B \leftarrow A \rightarrow C$ in $(\Alg^{r.c})$. 
 \end{prop}

\begin{pf} (1) Since $B$ is generated by $Z(B)$ and the image of $A$ as a ring, an element of $C_B(A)$ commutes all elements of $B$.

(2)  The ring $B\otimes_A C$ is generated by $Z(B)\otimes 1$ and $1\otimes C$ as a ring. Because $B\to B \otimes_A C$ is a morphism in (Alg) (\ref{prod} (3)), $Z(B) \otimes 1$ is contained in the center of $B\otimes_A C$ by \ref{ZZ}. 

(3) For a ring $D$ and for a morphism $B\to D$ in (Alg), $C_B(A)$ is cent to $Z(D)$ by (1) and by \ref{ZZ}. By this, (3) follows from \ref{prod} (2). 

(4) follows from (3). 
\end{pf}

 \begin{prop}\label{Morita} Assume that  a ring $A$ and a ring $A'$ are Morita equivalent. Then:
 
 (1)  The category of $A$-algebras  and the category of $A'$-algebras  are equivalent.
 
 (2) The category of r.c. $A$-algebras  and the category of r.c. $A'$-algebras  are equivalent.

 \end{prop}
\begin{pf}

 There is a $(A, A')$-bi-module $P$
  which is finitely generated faithful and projective over these rings, and we have $A'=\text{End}_A(P)^{\circ}$ and $A= \text{End}_{(A')^{\circ}}(P)$.   Here $ (\;)^{\circ}$ denotes the opposite ring. We have the functor from the category of $A$-algebras to the category of $A'$-algebras given by
 $B \mapsto B':= \text{End}_B(B \otimes_A P)^{\circ}$. The converse functor is obtained by  $B'\mapsto B =\text{End}_{(B')^{\circ}}(P \otimes_{A'} B')$. 
 
 This functor induces also the equivalence in (2). 
 \end{pf}

\begin{para}\label{mat} 

In the case $A'$ is the matrix algebra $M_n(A)$ ($n\geq 1$), the equivalence in \ref{Morita} is given by sending an $A$-algebra $B$ to the $M_n(A)$-algebra $B'=M_n(B)= A'\otimes_A B$. The converse functor is given by $B'\mapsto C_{A'}(B')$.

\end{para} 

\begin{prop}\label{Az} Let $A$ be a commutative ring and let $A'$ be an Azumaya algebra over $A$. Then the category of $A$-algebras and the category of $A'$-algebras are equivalent by the functor $B\mapsto B'=A'\otimes_A B$ from the former category  to the latter category. The converse functor is given by $B'\mapsto B=C_{A'}(B')$. This equivalence induces an equivalence between the categories of r.c. algebras. 
 \end{prop}

\begin{pf} \'Etale locally on  $\Spec(A)$, $A'$ is isomorphic to $M_n(A)$. Hence we are reduced to \ref{mat} by \'etale descent for the \'etale topology of $\Spec(A)$.  
\end{pf}

\section{Non-commutative schemes}\label{sec2}

\begin{para} Recall that a prime ideal of a ring $A$ (which may not be commutative) is a two-sided ideal $\frak p\neq A $ of $A$ such that if $I$ and $J$ are two-sided ideals of $A$ and if $IJ\subset\frak p$, then we have either $I\subset \frak p$ or $J\subset \frak p$. 

Let $\Spec(A)$ be the set of all prime ideals of $A$. It has the Zariski topology (Jacobson topology). A closed subset is a set $V(I)=\{\frak p\in \Spec(A)\;|\; I \subset \frak p\}$ given for  a two-sided ideal $I$ of $A$. We have $\cap_{I\in S} V(I) =V(\sum_{I\in S} I)$ for a                                                                                                                                                  set $S$ 
of two-sided ideals of $A$ and $V(I) \cup V(J)=V(IJ)$ for two-sided ideals $I$ and $J$ of $A$.

For $f\in A$, let $D(f)=\{\frak p\in \Spec(A)\;|\; f\notin \frak p\}$. Then the topology of $\Spec(A)$ is the weakest topology for which $D(f)$ is open for every $f\in A$. For $\frak p\in \Spec(A)$, $D(f)$ ($f\in A$) such that $\frak p\in D(f)$ form a base of neighborhoods of $\frak p$.

\end{para}

\begin{para}\label{Spec}
For a  ring homomorphism $A\to B$, the inverse image of a prime ideal of $B$ need not be a prime ideal of $A$. A standard example is the ring homomorphism $k \times k \to M_2(k)$ for a commutative field $k$,  the embedding as diagonal matrices. The ideal $(0)$ of $M_2(k)$ is a prime ideal, but its inverse image is the ideal $(0)$ of $k \times k$ which is not a prime ideal. 

However, by Procesi \cite{Pr}, for a morphism $h: A \to B$ in (Alg), we have a map
$\Spec(B) \to \Spec(A)\;;\; \frak p\mapsto h^{-1}(\frak p)$. This map is continuous.

\end{para}

\begin{para}\label{spi} By a {\it non-commutative space of prime ideals}, we mean a triple $(X, \cO_X, (\frak p(x))_{x\in X})$, where $X$ is a topological space, $\cO_X$ is a sheaf of rings (which need not be commutative) on $X$, and $\frak p(x)$ for each $x\in X$ is a prime ideal of the stalk $\cO_{X,x}$  satisfying the following conditions (i) and (ii).

(i) For each $x\in X$, the center $Z(\cO_{X,x})$ of $\cO_{X,x}$ is a local ring, and $Z(\cO_{X,x})\cap \frak p(x)$ is the maximal ideal of $Z(\cO_{X,x})$. 

(ii) If  $U$ is an open set of $X$ and if $f\in \cO_X(U)$, the set $\{x\in U\;|\; f_x\notin \frak p(x)\}$ is open in $U$.

\medskip

For non-commutative spaces of prime ideals $X=(X, \cO_X, (\frak p(x))_{x\in X})$ and $Y=(Y, \cO_Y, (\frak p(y))_{y\in Y})$, a morphism $X\to Y$ means a pair of a continuous map $f: X \to Y$ and a homomorphism of sheaves of rings $f^{-1}(\cO_Y) \to \cO_X$ on $X$ satisfying the following condition (iii).

(iii) For each $x\in X$ with image $y$ in $Y$, the ring homomorphism $\cO_{Y, y}\to \cO_{X,x}$ is a morphism of $(\Alg)$, and $\frak p(y)$ is the inverse image of $\frak p(x)$.

\end{para}

\begin{para}\label{F1} We consider the following

\medskip

Condition ({\bf f}) on a ring $A$: 

As a $Z(A)$-algebra, $A$ is finitely generated.

\end{para}

\begin{prop}\label{flatc} Let $A$ be a ring satisfying ({\bf f}) and let $R$ be a flat commutative ring over $Z(A)$. Then $R$ is the center of $R \otimes_{Z(A)} A$.

\end{prop}

\begin{pf} Take finitely many elements  $a_i\in A$ ($1\leq i\leq n$) 
which generate                                                                                                                  
 $A$ over $Z(A)$ as a ring.  We have an exact sequence $0 \to Z(A) \to A \overset{h}\to A^n$, where $h(x)=(a_ix-xa_i)_i$. Since $R$ is flat over $Z(A)$, the sequence $0\to R \to R \otimes_{Z(A)} A \overset{h}\to (R\otimes_{Z(A)} A)^n$ is exact.\end{pf}

\begin{para}\label{specspi1} For a ring $A$ which satisfies ({\bf f}) in \ref{F1},  we have
a non-commutative space of prime ideals
$(\Spec(A), \cO_A, (\frak p(x))_{x\in \Spec(A)})$,  where $\cO_A$ and $\frak p(x)$ are as follows.

Define the sheaf $\cO_A$ of rings on  $\Spec(A)$ as follows. Let $Z(A)$ be the center of $A$, and let $\cO_{Z(A)}$ be the structure sheaf of $\Spec(Z(A))$. Let $\cO_A$ be the inverse image of the sheaf $\cO_{Z(A)}\otimes_{Z(A)} A$ on $\Spec(Z(A))$ under the canonical map $\Spec(A) \to \Spec(Z(A))$ (\ref{Spec}).

Let $x\in X$ and let $y$ be the image of $x$ under $\Spec(A) \to \Spec(Z(A))$ and let $Z(A)_y$ be the local ring of $Z(A)$ at $y$. Then $Z(A)_y$ is the center of $\cO_{A,x}$ by \ref{flatc}. Let $\frak q$ be the prime ideal of $A$ corresponding to $x$, and let $\frak p(x):= Z(A)_y \otimes_{Z(A)} \frak q\subset \cO_{A,x}$. 

We denote this non-commutative  space of prime ideals simply  by $\Spec(A)$.

\end{para}

\begin{rem} We have the evident  canonical ring homomorphism $A \to \Gamma(\Spec(A), \cO_A)$. 
This map need not be bijective. See \ref{phil0} (1).

\end{rem}

\begin{thm}\label{specspi2} Let $A$ be a ring which is finitely generated as a $Z(A)$-algebra and let $X$ be a non-commutative space of prime ideals. Consider the map
$$\text{Mor}(X, \Spec(A)) \to \Hom_{(\text{Rings})}(A, \cO(X))$$
(Mor denotes the set of morphisms of non-commutative spaces of prime ideals, and  $\cO(X):=\Gamma(X, \cO_X)$) which sends a morphism $\psi$ on the left hand side to the composition $A\to \cO(\Spec(A)) \to \cO(X)$ where the second arrow is induced by $\psi$. Then this map is injective and the image consists of all homomorphisms satisfying  the following condition (i).

(i) For every $x\in X$, the composition $A \to \cO(X) \to \cO_{X,x}$ is a morphism in $(\Alg)$. 
\end{thm}

\begin{pf} Let $H$ be the subset of $\Hom(A, \cO(X))$ consisting of all elements which satisfy (i).   Let $Y=\Spec(A)$, $Z=\Spec(Z(A))$.

Let $\psi: X\to Y$ be a morphism and let $h: A\to \cO(X)$ be the associated ring homomorphism. We prove $h\in H$. 
Let $x\in X$, and let $y\in Y$ and $z\in Z$ be the images of $x$. 
Then  the map $A\to \cO_{X,x}$ is the composition  $A\to  \cO_{Z,z}\otimes_{Z(A)} A = \cO_{Y,y}\to \cO_{X,x}$. The first arrow of this composition is a morphism in $(\Alg)$. Since the second arrow is a morphism of $(\Alg)$, the composition is also a morphism in $(\Alg)$.  

Let $h\in H$. We define a morphism $\psi: X \to Y$ of non-commutative spaces of prime ideals associated to $h$ as follows. For $x\in X$, let $\psi(x)\in Y$ be the inverse image of $\frak p(x)$ under $A\to \cO_{X,x}$. Then $\psi$ is continuous by the condition (ii) in \ref{spi}. For $x\in X$, since $Z(A) \to \cO_{X,x}$ is a morphism in $(\Alg)$, the image of $Z(A)$ in $\cO_{X,x}$ is contained in the center $Z(\cO_{X,x})$ of $\cO_{X, x}$ by \ref{ZZ}. 
Let $f\in Z(A)$ and let $D_X(f)=\{x\in X\;|\; f_x\notin \frak p(x) \;\text{in}\; \cO_{X,x}\}$. 
Then for $x\in D_X(f)$, $f_x\in Z(\cO_{X,x})$ is not contained in the intersection of $\frak p(x)$ and $Z(\cO_{X,x})$ which is the maximal ideal of  the local ring $Z(\cO_{X,x})$. Hence $f_x$ is invertible in $Z(\cO_{X,x})$. By this,  we have a homomorphism from the pullback of $\cO_Z$ on $X$ to $Z(\cO_X)$. Hence we have a homomorphism from the inverse image of $\cO_Z \otimes_{Z(A)} A$ on $X$ to $\cO_X$. Thus we have a morphism of ringed spaces $X\to Y$, which we still denote by $\psi$,  and this is a morphism of non-commutative spaces of prime ideals. 

It is easy to see that $H\to \text{Mor}(X, Y)\;;\; h\mapsto \psi$ and $\text{Mor}(X,Y)\to H\;;\; \psi\mapsto h$ are the inverse maps of each other. 
\end{pf}

\begin{para} By a {\it non-commutative scheme}, we mean a non-commutative  space of prime ideals $X$ satisfying the following condition: For each $x\in X$, there is an open neighborhood $U$ of $x$ such that as a non-commutative space of prime ideals, $U$ is isomorphic to an open  subspace of $\Spec(A)$ (\ref{specspi1}) for some ring $A$ which satisfies ({\bf f}).

A morphism of non-commutative schemes is a morphism as non-commutative spaces of prime ideals. 

\end{para}

\begin{para} 
 We have the contra-variant functor $A \mapsto \Spec(A)$  from the full subcategory of (Alg) consisting of rings satisfying ({\bf f})   to the category of non-commutative schemes. 
 
This functor is not fully faithful. See \ref{S3} (2).

\end{para}

\begin{para}\label{S3} {\bf Example.} 
Let $k$ be a commutative field of characteristic $3$ and  let $A$ be the group ring $k[S_3]$ of the symmetric group $S_3$ of degree $3$. We show: 

\medskip

(1) The canonical map $A\to \cO(X)$ is not bijective. 

(2) The canonical map $\Hom_{(\Alg)}(A, A) \to \text{Mor}(\Spec(A),\Spec(A))$ is not bijective. 

\medskip

 The group $S_3$ is defined by generators $\alpha, \beta$ with relations $\alpha^2=1$, $\beta^3=1$, $\alpha \beta \alpha^{-1}=\beta^{-1}$. 
Let $I$ be the two-sided ideal of $A$ generated by $\beta-1$. Then $I^3=0$ and $A/I\cong k\times k$ where $\alpha\in A$ goes to $(1,-1)\in k \times k$.  The set  $\Spec(A)$ 
consists of  two elements $\frak p_1$ and $\frak p_2$ where $\frak p_i$ is the composition of $A\to A/I\cong k \times k$ and the $i$-th projection $k\times k \to k$. The center $Z(A)$ is an Artinian local ring and we have $Z(A) \cong k[u, v]/(u^2, v^2, uv)$, where $u$ corresponds to $1+\beta+\beta^2$ and $v$ corresponds to $\alpha(1+\beta+\beta^2)$. 

As a non-commutative scheme,  $\Spec(A)$ is the disjoint union of its two open subspaces $U_1=\{\frak p_1\}$ and $U_2= \{\frak p_2\}$ though $A$ is not a direct product of two non-zero rings. The maps $A\to \cO(U_i)$ ($i=1,2$) are isomorphisms.  

Hence $\cO(\Spec(A))= A \times A$, and $A\to \cO(\Spec(A))$ is the diagonal map $A\to A \times A$ which is not an isomorphism. This shows (1). 

By \ref{specspi2}, we have  $ \text{Mor}(\Spec(A), \Spec(A))= \Hom_{(\Alg)}(A, \Gamma(\Spec(A), \cO_A))= \Hom_{(\Alg)}(A, A \times A)= \Hom_{(\Alg)}(A, A)\times \Hom_{(\Alg)}(A,A)$. The diagonal map $\Hom_{(\Alg)}(A, A) \to \Hom_{(\Alg)}(A, A) \times \Hom_{(\Alg)}(A, A)$  is not bijective, for the set $\Hom_{(\Alg)}(A, A)$ has an element which is not the identity map (for example, the ring homomorphism $x \mapsto \alpha x\alpha^{-1}$).  This shows (2).

\end{para}

\begin{rem}\label{phil0} 

 The author does not think that (1), (2) in \ref{S3} are  weak points of our formulation of non-commutative geometry.

The fact $\Spec(A)$ is the disjoint union of $U_1$ and $U_2$ in \ref{S3}  is important for the relation with zeta functions of non-commutative rings. 
 \'Etale cohomology theory was created by Grothendieck in his efforts to have a cohomology theory which  explains properties of  zeta functions of varieties over a finite field predicted in Weil conjectures. Our dream is to have an \'etale cohomology theory of non-commutative rings which can explain 
 zeta functions of non-commutative rings (\cite{Fu}) over finite fields, and we hope  that this paper is a starting point.   In Appendix of this paper, T. Fukaya gives a partial result (\ref{Hczeta}).  
 
In \ref{S3}, let $k$ be the algebraic closure of $\F_3$ and let $A_0= \F_3[S_3]$.  Then the zeta function $\zeta_{A_0}(s)$ of $A_0$ is $ (1-3^{-s})^{-2}$ where each $(1-3^{-s})^{-1}$ corresponds to each point of $\Spec(A_0)$.  Because $\Spec(A)$ is the disjoint union of $U_1$ and $U_2$, the \'etale cohomology group 
 $H^0(A _{et}, \Q_{\ell})$ is two-dimensional, and we have the good relation $\zeta_{A_0}(s) = \text{det}(1- \varphi^{-1}_3 u\;|\; H^0(A_{et}, \Q_{\ell}))^{-1}$ of the zeta function of $A_0$ and the \'etale cohomology of $\Spec(A)$,  where $u=3^{-s}$ and $\varphi_3\in \Gal(\bar k/k)$ is  the Frobenius operator which acts here as the identity map on the cohomology group.

If $X$ is a ringed space  such that $\cO(X)=A$,  $X$ must be connected because $\cO(X)$ is not a product of two non-zero rings, and hence $H^0(X)$ should be 
one-dimensional for any reasonable cohomology theory, and hence we would not have a presentation of $\zeta_{A_0}(s)$ by using the cohomology of $X$. Such $X$ is not nice for us in the relation to the zeta function.

\end{rem}

\section{Basics of non-commutative schemes, 1}\label{sec3}

We consider the non-commutative analogue \ref{SpecB}  of the fact that a quasi-coherent sheaf $\cB$ of commutative rings  on a scheme $S$ defines a scheme $\Spec(\cB)$  over $S$. 

We also consider the ``tensor products'' and the  fiber products in the category of non-commutative schemes.

In this section \ref{sec3},  $X$ denotes a non-commutative scheme.

\begin{para} For a ring $A$, by an $A$-algebra of finite presentation, we  mean an $A$-algebra which is isomorphic over $A$ to $A\langle T_1, \dots, T_n\rangle/(f_1,\dots, f_m)$ for some $n, m$ and $f_1, \dots, f_m\in A\langle T_1,\dots, T_n\rangle$. 

By an  $\cO_X$-algebra locally of finite presentation (an l.f.p. $\cO_X$-algebra for short), we mean a sheaf of rings over $\cO_X$ on $X$ which is locally isomorphic to $\cO_X\langle T_1,\dots, T_n\rangle/(f_1, \dots, f_m)$ for some $n, m$ and sections $f_1, \dots, f_m$ of $\cO_X\langle T_1,\dots, T_n\rangle$.

\end{para}

\begin{lem}\label{lfp0} Let $\cB$ be an l.f.p. $\cO_X$-algebra.

Assume that we are given a ring $A$ satisfying ({\bf f}) and an isomorphism between $X$ and an open subspace of $\Spec(A)$. Then locally on $X$, there are an $A$-algebra $B$ of finite presentation and an isomorphism $\cB\cong \cO_X\otimes_A B$ over $\cO_X$. 

\end{lem} 

\begin{pf} Let $x\in X$, and let $z\in \Spec(Z(A))$ be the images of $x$. Take an open neighborhood $U$ of $x$ and an isomorphism $\cB|_U \cong \cO_U\langle T_1, \dots, T_n\rangle/(f_1, \dots, f_m)$. Since $\cO_{X,x}\cong \cO_{Z,z}\otimes_{Z(A)} A$, there is $h\in Z(A)$ which is invertible at $z$ such that the stalk $f_{i,x}$ comes from an element $g_i$  of $A[1/h]\langle T_1, \dots, T_n\rangle$ for  $1\leq i\leq m$.  Replacing $X$ by an open neighborhood of $x$, we may assume that $f_i$ coincides with the pullback of $g_i$. Let $B:= A[1/h]\langle T_1, \dots, T_n\rangle/(g_1,\dots, g_m)$ and $U'=U\cap D(h)$. Then $\cB|_{U'} \cong \cO_{U'}\otimes_A B$, and $B$ is an $A$-algebra of finite presentation because $A[1/h]=A[T]/(hT-1)$. 
\end{pf}

\begin{prop}\label{SpecB} Let $\cB$ be an l.f.p.  $\cO_X$-algebra. Then there is a non-commutative scheme $\Spec(\cB)$ over $X$ which represents the following contra-variant functor $F$ from the category of non-commutative spaces of prime ideals  over $X$ to ($\Sets$). 

Let $f: Y \to X$ be a non-commutative scheme over $X$. Then $F(Y)$ is the set of all homomorphisms  $f^{-1}(\cB) \to \cO_Y$ of sheaves of rings over $f^{-1}(\cO_X)$ on $Y$ such that for each $y\in Y$ with image $x$ in $X$, the induced map $\cB_x\to \cO_{Y,y}$ is  a morphism in ($\Alg$). 

\end{prop}

\begin{pf} Working locally on $X$, we may assume that we have $A$, $B$, and an isomorphism $\cB\cong \cO_X \otimes_A B$ as in \ref{lfp0}.
Then the inverse image  of $X$ in $\Spec(B)$, which is an open subspace of $\Spec(B)$,  satisfies the condition of $\Spec(\cB)$  by \ref{specspi2} (which we apply by taking the present $B$ and an open set $U$ of $Y$ as $A$ and $X$ there, respectively, to have the relation of $\Mor(U, \Spec(B))$ between $\Hom_{(\text{Rings})}(B, \cO(U))$). 
\end{pf}

\begin{para}\label{lfp1}

Let $f:Y\to X$ be a morphism of non-commutative schemes. We say $f$ is locally of finite presentation (l.f.p. for short)  if the following  condition (i)  is satisfied. 

(i) Locally on $X$ and of $Y$, there is an l.f.p. $\cO_X$-algebra $\cB$  such that $Y$ is isomorphic over $X$ to an open subspace of $\Spec(\cB)$.

\end{para}

\begin{lem}\label{lfp2} (1) If $f: Y\to X$ and $g: Z\to Y$ are l.f.p., the composition $f\circ g : Z\to X$ is also l.f.p. 

(2) If $Y\to X$ and $Z\to X$ are l.f.p.  and if $f: Z\to Y$ is a morphism over $X$, then $f$ is l.f.p. 
\end{lem} 

\begin{pf} (1) By the proof of \ref{SpecB}, we may assume that there are a ring $A$ satisfying ({\bf f}), an $A$-algebra $B$ of finite presentation, an isomorphism between $X$ and an open subspace of $\Spec(A)$, and an isomorphism  over $X$ between $Y$ and an open subspace of the inverse image of $X$ in $\Spec(B)$. We may assume also that there are an l.f.p.  $\cO_Y$-algebra $\cC$ and an isomorphism over $Y$ between $Z$ and an open subspace of $\Spec(\cC)$. By \ref{lfp0} which we apply by taking $(Y, B, \cC)$ as $(X, A, \cB)$ of \ref{lfp0}, we have that locally on $Y$ and on $Z$, there is a $B$-algebra $C$ of finite presentation such that $\cC\cong \cO_Y\otimes_B C$ over $\cO_Y$. Then $Z$ is isomorphic over $X$ to an open subspace of $\Spec(C)$. 

(2) We may assume that $X=\Spec(A)$ for some ring $A$ satisfying ({\bf f}) and $Y=\Spec(B)$ for some $A$-algebra $B$ of finite presentation. By \ref{lfp0}, we may assume that there is a $A$-algebra $C$ of finite presentation and that $Z$ is an open subspace of $\Spec(C)$. Let $z\in Z$. Then there are an element $h$ of the center of $C$ such that $z\in D(h)$, an open neighborhood $U$ of $z$ in $D(h)$, and an $A$-homomorphism $B\to C[1/h]$ which induces the morphism $U\to Y$. Since $B$ and $C$ are of finite presentation over $A$, $C[1/h]$ is of finite presentation over $B$.  
\end{pf}

\begin{thm}\label{lpr}  Assume we are given morphisms $f: Y\to X$ and $g:Z\to X$ of  non-commutative schemes. Assume at least one of $f$ and $g$ is l.f.p. Then there is a non-commutative scheme over $X$, which we denote by $Y\otimes_X Z$, which represents the following contra-variant functor $F$ from the category of non-commutative spaces of prime ideals over $X$ to (Sets). 

For a non-commutative space of prime ideals  $W$ over $X$, 
$F(W)$ is the set of pairs $(u,v)$ of morphisms $u: W\to Y$ and $v: W\to Z$ of non-commutative spaces of prime ideals over $X$ satisfying the following condition:
For each $w\in W$ with images $y\in Y$, $z\in Z$ and $x\in X$,  and for  each $a\in \cO_{Y,y}$ and $b\in \cO_{Z, z}$ which commute with the image of every element of $\cO_{X,x}$, the images of $a$ and $b$ in $\cO_{W,w}$ commute. 

If $g$ is l.f.p., the canonical morphism $Y \otimes_X Z \to Y$ is l.f.p.
\end{thm}

\begin{pf} We may assume that $Z=\Spec(\cB)$ for an l.f.p.  $\cO_X$-algebra $\cB$. Then  $Y\otimes_X Z$ is given as $\Spec(f^*\cB)$, 
where $f^*\cB= \cO_Y\otimes_{f^{-1}(\cO_X)} f^{-1}(\cB)$. 
\end{pf}

\begin{para}\label{otimes2} Let $A\to B$ and $A\to C$ be morphisms in (Alg) and assume $A$, $B$ and $C$  satisfy ({\bf f}). Assume at least one of $A\to B$ and $A\to C$ is of finite presentation. Then 
$\Spec(B)\otimes_{\Spec(A)} \Spec(C)=\Spec(B \otimes_A C)$.

\end{para}

\begin{para}\label{rcmor} Let $f:Y\to X$ be a morphism of non-commutative schemes. We say $f$ is r.c. (relatively commutative) if for every $y\in Y$ and for $x:=f(y)$, $\cO_{Y,y}$ is an r.c. $\cO_{X,x}$-algebra. 

Clearly, if $f$ and $g:Z\to Y$ are r.c., then $f\circ g$ is r.c.
\end{para}

\begin{prop}\label{AzX} Let  $S$ be a scheme, let $\cB$ be a sheaf of Azumaya algebras over $\cO_S$ on $S$, and assume that $X=\Spec(\cB)$. Then the functor $X\otimes_S$ gives an equivalence from  the category of non-commutative schemes over $S$ to the category of non-commutative schemes over $X$. 
It gives an equivalence between the category of schemes over $S$ and the category of r.c. non-commutative schemes over $X$. 
\end{prop}

\begin{pf}
This follows from \ref{Az} and the fact that for an Azumaya algebra $A$ over $R$, the sets of two-sided ideals are in bijection as is well known (it is deduced  from the case $A=M_n(R)$ by the \'etale descent for the \'etale topology of $\Spec(R)$)  and hence $\Spec(A) \to \Spec(R)$ is a homeomorphism. 
\end{pf}

We consider fiber products in the category of non-commutative schemes. 

\begin{thm}\label{fibp} Let $f:Y\to X$ and $g:Z\to X$ be morphisms of non-commutative schemes. Assume that at least one of $f$ and $g$ is l.f.p. and assume that at least one of $f$ and $g$ is  r.c.  Then $Y\otimes_X Z$ is  the fiber product of $Y\to X \leftarrow Z$ in the category of non-commutative spaces of prime ideals. Write $Y\otimes_X Z$ by $Y\times_X Z$ in this case. If $g$ is r.c., the canonical morphism $Y\times_X Z\to Y$ is r.c.

\end{thm}

\begin{pf} This follows from  \ref{rcthm} (3) and from  the following \ref{cS}. 
\end{pf}

\begin{lem}\label{cS}
 Let $Y\to X$ and $Z\to X$ be morphisms of non-commutative schemes, and assume at least one of these morphisms is l.f.p. Let $W=Y\otimes_X Z$, let $w\in W$, and let $y$, $z$, and $x$ be the images of $w$ in  $Y$, $Z$, and $X$, respectively. Let $\cS$ be the set of all elements of the center of $\cO_{Y,y}\otimes_{\cO_{X,x}} \cO_{Z,z}$ whose images in $\cO_{W,w}$ are invertible. Then we have an isomorphism $\cS^{-1}(\cO_{Y,y} \otimes_{\cO_{X,x}} \cO_{Z,z})\overset{\cong}\to \cO_{W,w}$. 

\end{lem}

\begin{pf} We may assume that $Z=\Spec(\cB)$ for an l.f.p. $\cO_X$-algebra $\cB$. 

In the case $Y=X$, this  follows from the construction of $\Spec(\cB)$ in the proof of \ref{SpecB}. The general case follows from the case $Y=X$  and the construction of $Y\otimes_X Z$ in the proof of \ref{lpr}. 
\end{pf}

\begin{lem}\label{lfp3}  

Let $Y\to X$ be an l.f.p.  morphism and assume that it is r.c. Assume that we are given a ring $A$ satisfying ({\bf f}) and an isomorphism between $X$ and an open subspace of $\Spec(A)$. Then locally 
 on $X$ and on $Y$, there exist an $A$-algebra $B$ of the form $A[T_1, \dots, T_n]/(f_1, \dots, f_m)$ and an isomorphism over $X$ between $Y$ and an open subspace of the inverse image of $X$ in $\Spec(B)$.

\end{lem}

\begin{pf}  By \ref{lfp0}, we may assume that there are an $A$-algebra $B$ of finite presentation and an isomorphism over $X$ between $Y$ and an open subspace of the inverse image of $X$ in $\Spec(B)$. Take a finite presentation of $B$ as a quotient of $A\langle T_1, \dots, T_r\rangle$ and let $b_i$ be the image of $T_i$ in $B$. Let $y\in Y$ and let $s$ be the image of $Y$ in $S:= \Spec(Z(B))$. Since $\cO_{Y, y}=\cO_{S,s} \otimes_{Z(B)} B$ and this is relatively commutative over $A$, there are finitely many elements $t_1,\dots, t_n$ of $Z(B)$ and an element $h$ of $Z(B)$ which is invertible at $s$  such that all $b_i$ are contained in the subring of $B[1/h]$ generated by the image of $A$ and $t_1, \dots, t_n$ and $1/h$. We have the surjective homomorphism $A[T_1, \dots, T_n, T_{n+1}]\to B[1/h]$ over $A$ which sends  $T_i$ to $b_i$ for $1\leq i\leq n$ and $T_{n+1}$ to $1/h$. 
Since $B[1/h]$ is of finite presentation as an $A$-algebra, the kernel of this surjection is a finitely generated two-sided ideal. Furthermore, for $U =D(h)\subset Y$, $U$ is isomorphic over $X$ with an open subspace of the inverse image of $X$ in $\Spec(B[1/h])$. 
\end{pf}

\begin{para}\label{diag} For a morphism $f:Y\to X$ of non-commutative schemes which is l.f.p. and  r.c.,  we have the diagonal morphism $\Delta_f: Y\to Y\times_X Y$ which corresponds the pair $(Y\to Y, Y\to Y)$ of the identity morphisms.

\end{para}
We consider immersions (open immersions, l.f.p. closed immersions, and l.f.p. immersions) of non-commutative schemes,

\begin{lem}\label{clim} For a morphism $f:Y\to X$ of non-commutative schemes, the following conditions (i) and (ii) are equivalent.

(i) Locally on $X$,  there are a ring $A$ satisfying ({\bf f}), a finitely generated two-sided ideal $I$ of $A$, and an isomorphism between $X$ and an open subspace of $\Spec(A)$, such that $Y$ is isomorphic over $X$  to the inverse image of $X$ in $\Spec(A/I)$ regarded as an open subspace of  $\Spec(A/I)$.

(ii) $Y$ is isomorphic over $X$  to $\Spec(\cB)$ over $X$ for some l.f.p. $\cO_X$-algebra $\cB$  such that the map $\cO_X\to \cB$ is surjective.

\end{lem}

\begin{pf} Assume (i). Then locally on $X$, with $A$ and $I$ as in (i), we have $Y=\Spec(\cB)$, where $\cB= \cO_X/I\cO_X$. When we regard $\cO_Y$ naturally as a sheaf on $X$ via the homeomorphism from $Y$ to the image of $Y$ in $X$, which is a closed subset of $X$, $\cB$ is identified with the image of $\cO_X\to \cO_Y$ and hence is independent of the choice of $(A,I)$. Hence $\cB$ is glued on the whole $X$ and we have $Y=\Spec(\cB)$ globally on $X$.

Assume (ii). Locally, take a ring $A$ satisfying ({\bf f}) and an isomorphism between $X$ and an open subspace of $\Spec(A)$. Then locally on $X$, for some element $h$ of $Z(A)$ which is invertible on $X$ and for some  $f_1, \dots, f_m\in A[1/h]$, we have $\cB= \cO_X/(f_1, \dots, f_n)= \cO_X \otimes_{A[1/h]} A[1/h]/(f_1, \dots, f_m)$. Hence we have (i).  
\end{pf}

\begin{para}\label{im} 

(1) A morphism of non-commutative schemes $Y\to X$ is called an open immersion if $Y$ is isomorphic over $X$ to an open subspace of $X$.

(2) A morphism of non-commutative schemes $Y\to X$ is called a closed immersion locally of finite presentation (l.f.p. closed immersion for short) if it satisfies the equivalent  conditions in \ref{clim}.

If $f:Y\to X$ is an l.f.p. closed immersion, $f$ gives a homeomorphism from $Y$ to a closed subset of $X$.

(3) A morphism of non-commutative schemes $Y\to X$ is called an l.f.p. immersion if it is a composite $Y\overset{a}\to U\overset{b}\to X$ such that $a$ is an l.f.p. closed immersion and $b$ is an open immersion.

\end{para}

\begin{lem} (1)  If $f: Y\to X$ and $g: Z\to Y$ are open (resp. l.f.p. closed,  resp. l.f.p.) immersions, the composition $f\circ g: Z\to X$ is an open (resp.   l.f.p. closed, resp. l.f.p.) immersion.

(2) If $f: Z\to X$ is an open (resp. l.f.p. closed, resp. l.f.p.) immersion and if $Y\to X$ is a morphism of non-commutative schemes, the induced morphism $Y\times_X Z \to Y$ is an open (resp. l.f.p. closed, resp. l.f.p.) immersion. (Note that an l.f.p. immersion is r.c. and hence we have $Y\times_X Z$.)
\end{lem}
\begin{pf} (1)  for open immersion is clear. 

(1) for l.f.p. closed immersion.  We may assume that $X$ is an open subspace of $\Spec(A)$ for a ring $A$ satisfying ({\bf f}) and $Y=X\times_{\Spec(A)} \Spec(A/I)$ for a finitely generated two-sided ideal $I$ of $A$. We have $Z=\Spec(\cC)$ for some $\cO_Y$-algebra $\cC$ locally of finite presentation such that $\cO_Y\to \cC$ is surjective. Working locally on $X$, we may assume that there are an element $h$ of $Z(A/I)$ which is invertible on $Y$  and elements $f_1, \dots, f_m\in A$ such that $\cC=\cO_Y/(h^{-1}f_1, \dots, h^{-1}f_m)$. Let $J$ be the ideal of $A$ generated by $I$ and $f_1, \dots, f_m$. Then we have $Z=X\times_{\Spec(A)} \Spec(A/J)$.  

(1) for l.f.p. immersion. Since the composition of open (resp. l.f.p. closed) immersions is an open (resp. l.f.p. closed) immersion as we have just proved,  it is sufficient to prove that  if $f: Y\to X$ is an l.f.p. closed immersion and if $g: Z\to Y$ is an open immersion, 
then $f\circ g: Z\to X$ is an l.f.p. immersion. There is an open set $U$ of $X$ such that $Z=Y\cap U$. Then $f\circ g$ is the composition $Z\to U \to X$, where the first arrow is  an l.f.p. closed immersion (because it is the restriction of $Y\to X$ to $U$) and the second arrow is an open immersion.  

(2) is clear. 
\end{pf}

\begin{para}\label{sep}
Let $f:Y\to X$ be an r.c. and l.f.p. morphism of non-commutative schemes. Consider the diagonal morphism $\Delta_f: X\to Y\times_X Y$ (\ref{diag}) which is l.f.p. by \ref{lfp2} (2). 

We say $f$ is separated (resp. unramified) if  $\Delta_f$ is an l.f.p. closed immersion (resp. open immersion).

\end{para}

\begin{lem}\label{sep2}

(1) If $Y\to X$ and $Z\to Y$  are separated (resp. unramified) (then these morphisms are assumed to be r.c. and l.f.p.), then the composition $Z\to X$ is separated (resp. unramified).

(2) If $f: Z\to X$ is separated (resp. unramified), then for every morphism $Y\to X$ of non-commutative schemes, the induced morphism $Y\times_X Z\to Y$ is separated (resp. unramified).

(3) Let $f: Y\to X$ and $g:Z \to Y$ be morphisms of non-commutative schemes. Assume that $f$ is  separated and  $f\circ g$ is an l.f.p. immersion. Then $g$ is an l.f.p.  immersion.
\end{lem}

\begin{pf} 
(1) The diagonal morphism  $Z\to Z\times_X Z$ is the composition of $Z\to Z\times_Y Z$ and  $Z\times_Y Z \to Z\times_X Z$, in which the second arrow  is the base change of $Y\to Y\times_X Y$.  
 
 (2) is evident.

(3) The morphism $g$ is the composition $Z\to Y\times_X Z\to Y$. The first arrow is a base change of  the diagonal $Y\to Y\times_X Y$ and hence is an l.f.p. closed immersion. The second morphism is a base change of $Z\to X$ and hence an l.f.p. immersion. 
\end{pf}
\begin{lem}\label{unrlem}  Let $f:Y\to X$ be an unramified morphism. Let $y\in Y$, $x=f(y)$, $A=\cO_{X,x}$, $B=\cO_{Y, y}$. Let $\m$ be a maximal two-sided ideal of $A$ (which may not be $\frak p(x)$) and assume that  $A/\m$ is an Azumaya algebra over a commutative field $k$. Then for some $n\geq 1$, there are commutative fields $k_1, \dots, k_n$ which are finite separable extensions of $k$  such that $B/B\m\cong A/\m \otimes_k (k_1\times \dots \times k_n)$ over $A/\m$. (Note that the left ideal $B\m$ of $B$ is a two-sided ideal of $B$ because $B$ is an $A$-algebra.)

\end{lem}
 \begin{pf} Let $\cS$ be the set of all elements of the center of $B\otimes_A B$ whose images in $B$ under $B\otimes_A B\to B\; ; \; b\otimes c \mapsto bc$ are invertible. Then since the diagonal morphism $Y\to Y\times_X Y$ is an open immersion, $\cS^{-1}(B\otimes_A B) \overset{\cong}\to B$ by \ref{cS} applied to  the case $Y$ and $Z$ are $Y$, $y$ and $z$ are $y$, and $w=y\in Y\subset Y\times_X Y$.  
 By \ref{Az}, $B/B\m= A\otimes_k E$ for some  commutative ring $E$ over $k$. Hence we have   
  $\cS^{-1}(B/B\m \otimes_{A/\m} B/B\m) = B/B\m$, that is, 
  $\cS^{-1}(E \otimes_k E)\overset{\cong}\to E$, where $\cS$ is the set of all elements of $E\otimes_k E$ whose images are invertible in $E$. Hence $E\cong k_1\times\dots\times k_n$ for some finite separable extensions $k_1, \dots, k_n$ of $k$. 
\end{pf}

\begin{para}\label{irre} 
As in \cite{Bo} Chap. 2, Section 4.1,  we say a topological space $T$ is irreducible if every finite intersection of non-empty open subsets of $X$ is non-empty.

If $X$ is a non-commutative scheme, there is a bijection from $X$ to the set of all irreducible closed subsets of $X$ which sends $x\in X$ to the closure of $\{x\}$ in $X$. In fact this is reduced to the case $X=\Spec(A)$ for a ring $A$, and this case is well-known.  
\end{para}

\begin{para}\label{spr0}
 Recall that a ring $A$ is called a prime ring if the ideal $\{0\}$ of $A$ is a prime ideal. 
 
 We say that a non-commutative scheme  $X$ is prime if it is irreducible and locally on $X$, $X$ is isomorphic to an open subspace of $\Spec(A)$ for some prime ring $A$.
 \end{para}

\begin{lem}\label{spr1} Let $X$ be a prime non-commutative scheme.

(1) Let $\eta$ be the point of $X$ whose closure is $X$. Then the prime ideal $\frak p(\eta)$ of $\cO_{X,\eta}$ is $(0)$. For every non-empty open set $U$ of $X$, the map $\cO_X(U) \to \cO_{X,\eta}=\cO_{X,\eta}/\frak p(\eta)$ is injective.

(2) Let $A$ be a ring satisfying ({\bf f}), and let $I$ be a two-sided ideal of $A$. Let $f: X\to \Spec(A)$ be a morphism whose image  is contained in $\Spec(A/I)\subset \Spec(A)$. Then $f$ factors through a morphism $X\to \Spec(A/I)$  in a unique way. 
\end{lem}

\begin{pf} (1) We may assume $X=\Spec(A)$ for a prime ring $A$ satisfying ({\bf f}). Note that $Z(A)$ is an integral domain. The map $\cO_X(U) \to \prod_{x\in U} \cO_{X, x}$ is injective. For $x\in U$, the map $\cO_X(U)\to \cO_{X, \eta}$ is the composition $\cO_X(U)\to \cO_{X,x}\to \cO_{X, \eta}$ and the last arrow is $Z(A)_{\frak p} \otimes_{Z(A)} A\to Z(A)_{(0)} \otimes_{Z(A)} A$ where $\frak p$  is the image of $x$ in $\Spec(Z(A))$. Since all non-zero elements of $Z(A)$ are non-zero-divisors of $A$, this last arrow is injective. 

(2) By \ref{specspi2}, it is sufficient to prove that $I$ is killed by the induced homomorphism $A\to \cO(X)$. By (1),  it is sufficient to prove that $I$ is killed by the map $A\to \cO_{X,\eta}=\cO_{X,\eta}/\frak p(\eta)$. Let $y\in \Spec(A/I)$ be the image of $\eta$. Then the last map factors as $A\to \cO_{\Spec(A), y}/\frak p(y) \to \cO_{X,\eta}/\frak p(\eta)$ and the first arrow here  kills $I$. 
\end{pf}

\begin{prop}\label{spr3} Let $X$ be a non-commutative scheme and let 
$C$ be an irreducible  closed subset of $X$. Then:

(1)  $C$ has a unique  structure of a prime non-commutative scheme with which $C$ represents the functor on the category of prime non-commutative schemes
$$Z\mapsto \{\text{morphisms}\;f: Z\to X\; \text{such that}\; f(Z)\subset C\}.$$
The natural morphism $C\to X$ is underlain by the inclusion map from $C$ to $X$.

(2) Assume that $X$ is covered by open subspaces $\Spec(A_{\lam})$ such that for every $\lam$, every two-sided ideal 
of $A_{\lam}$ is finitely generated. Then with the structure of $C$ in (1), the canonical morphism $C\to X$ is an l.f.p. closed immersion.

\end{prop}

\begin{pf} (1) We can work locally and hence we may assume that $X$ is an open subspace of $\Spec(A)$ for a ring $A$ satisfying ({\bf f}) and that $C= X\cap \Spec(A/\frak p)$ for a prime ideal $\frak p$ of $A$. Then $C$ is endowed with the structure of a non-commutative scheme as an open subspace of $\Spec(A/\frak p)$. Then $C$ becomes a prime non-commutative scheme. The rest follows from  \ref{spr1} (2). 

(2) is clear. 
 \end{pf}
 
  \begin{rem} If $C$ is a closed subset of a scheme $S$, $C$ has a canonical structure of a scheme, the reduced closed subscheme structure. The author does not know whether a close subset $C$ of a non-commutative scheme $X$ has such a canonical structure, except the case of an irreducible closed subset discussed above. For example, consider the following $X$ and $C$. 
 We have   $\Spec(\Z_3[S_3])= \{\m,\m', \frak p, \frak p', \frak q\}$, where $\m= (\alpha-1, \beta-1, 3)$, $\m'=(\alpha+1,\beta-1, 3)$, $\frak p= (\alpha-1, \beta-1)$, $\frak p'=(\alpha+1, \beta-1)$, $\frak q= (\beta^2+\beta+1)$ with  $\alpha$, $\beta$ as in \ref{S3},  $\m$ and $\m'$ are closed points, the closure of $\frak p$ is $\{\m, \frak p\}$, the closure of $\frak p'$ is $\{\m', \frak p'\}$, and the closure of $\frak q$ is $\{\m,\m', \frak q\}$. Let $X$ be the open subspace of $\Spec(A)$ 
consisting of $\m, \frak p, \frak q$, and let $C$ be $X$ itself. Then $C$ can be regarded as the non-commutative scheme $X$ but $C$ has also the structure as an open subspace of $\Spec(\Z_3[S_3]/(\frak p\cap \frak q))= \{\m, \m', \frak p, \frak q\}$. These two structures are different because the stalk $\cO_{C,\m}$ for the former is $\Z_3[S_3]$ but $\cO_{C,\m}$ for the latter is $\Z_3[S_3]/(\frak p\cap \frak q)$. Since $A$ and $A/(\frak p\cap \frak q)$ are semi-prime rings (that is, the intersection of all prime ideals is $(0)$) and semi-prime rings are  non-commutative analogues of commutative reduced rings, these two structures on $C$  are both like reduced scheme structure, and which one should be regarded as the best  canonical structure is  not clear to the author.
 
 \end{rem}

\section{Basics of non-commutative schemes, 2}\label{sec4}

We discuss flat morphisms of non-commutative schemes. We then consider specially  non-commutative schemes which are ``finite over the center''. We discuss the openness of a flat morphism (\ref{a3}).

\begin{para}\label{flat1} Let $h:A\to B$ be a homomorphism of rings. 

We say $h$ is right (resp. leftt) flat  if through $h$,  $B$ is flat as a right (resp. left) $A$-module. (That is, the functor $B \otimes_A$ is exact for left (resp. right) $A$-modules.)
We say $f$ is flat if it is right flat and left flat. 
\end{para}

\begin{para}\label{flat2}  Let $f: X\to Y$ be a morphism of non-commutative schemes. 

We say $f$ is flat (resp. right flat, resp. left flat) if $\cO_{Y, f(x)}\to \cO_{X,x}$ is flat (resp. right flat, resp. left flat)  for every $x\in X$. 
\end{para}

In this Section \ref{sec4}, we consider the right flatness. But the left flatness is treated in the same way by considering opposite rings. 

\begin{lem}\label{lemflat}  (1) If $X\to Y$ and $Y\to Z$ are right flat morphisms of non-commutative schemes, the composition $X\to Z$ is right flat. 

(2) Let  $X\to Y$ be a right flat morphism. Then for a morphism $Y'\to Y$ of non-commutative schemes, if at least one of $X\to Y$ and $Y'\to Y$ is l.f.p., the induced morphism $X\otimes_Y Y'\to Y'$ is right flat.

\end{lem}

\begin{pf}

(1) is clear. 

(2)  Let $x'$ be a point of $X':=X\otimes_Y Y'$, let $x$ be the image of $x'$ in $X$, let $y'$ be the image of $x'$ in $Y'$, and let $y$ be the image of $y'$ in $Y$. Let $P=\cO_{X,x}\otimes_{\cO_{Y,y}} \cO_{Y', y'}$ and $Q=\cO_{X',x'}$. Then by \ref{cS}, there is a multiplicative subset $\cS$ in the center of $P$ such that the canonical map $P\to Q$ induces an isomorphism $\cS^{-1}P \overset{\cong}\to Q$. Since $P$ is right flat over $\cO_{Y', y'}$, $Q$ is right flat over $\cO_{Y',y'}$. \end{pf}

The following \ref{exist} should be well known. 
\begin{lem}\label{exist} Let $A$ be a prime ring  and assume that $A$ is a subring of a ring $B$. Then there is a prime ideal $\frak p$ of $B$ such that $A\cap \frak p= (0)$.

\end{lem}

\begin{pf} Consider the set $\frak S$ of all two-sided ideals $I$ of $B$ such that $A\cap I=(0)$. By Zorn's lemma, $\frak S$ has a maximal element $\frak p$. We prove that $\frak p$ is a prime ideal. Assume $I$ and $J$ are two-sided ideals of $B$ such that $\frak p\subset I$ and $\frak p \subset J$ and $IJ\subset \frak p$. Since $(A\cap I)(A \cap J)=(0)$ and $(0)$ of $A$ is a prime ideal, we have either $A\cap I=(0)$ or $A\cap J=(0)$. Assume $A\cap I=(0)$. 
By the property of $\frak p$ as a maximal element, we have $I =\frak p$. 
\end{pf}

\begin{prop}\label{ffus} Let $A\to B$ be a morphism of ($\Alg$) and assume it is faithfully right flat (that is, it is right flat and $B\otimes_A M \neq 0$ for every non-zero left $A$-module $M$).  
Then the morphism $f: \Spec(B)\to \Spec(A)$ is surjective. If $A$ and $B$ satisfy ({\bf f}) and $g: Y\to \Spec(A)$ is a morphism of non-commutative schemes and if either $f$ or $g$ is l.f.p., the map $\Spec(B) \otimes_{\Spec(A)} Y\to Y$ is surjective. 
\end{prop}

\begin{pf}

 By \ref{otimes2} (applied to $C=\cO_{Y, y}/\frak p(y)$ for $y\in Y$), it is sufficient to prove that $\Spec(B) \to \Spec(A)$ is surjective. 
Let $\frak p$ be a prime ideal of $A$. Note that since $B$ is an $A$-algebra, the left ideal $B\frak p$ is a two-sided ideal $B\frak p B$ of $B$. 
We show that  $A/\frak p\to B/B \frak p $ is injective. Let $M$ be the kernel of this map. By applying $B\otimes_A$ to the exact sequence $0\to M \to A/\frak p \to B/B\frak p$, we have an exact sequence $0\to B \otimes_A M \to B/B\frak p \to B \otimes_A B/\frak p$.  The  map $B/B \frak p  \to B \otimes_A B/B \frak p$  is injective because it has a section $x\otimes y\mapsto xy$. Hence $B\otimes_A M=0$, and hence $M=0$. By \ref{exist} applied to the injection $A/\frak p\to B/B \frak p $, there is a prime ideal $\frak q$ of $B$ whose inverse in $A$ is $\frak p$. 
\end{pf}

In the rest of this Section \ref{sec4}, we consider an algebra $A$ over a commutative ring $R$ which is finitely generated as an $R$-module, and non-commutative schemes which are locally an open subspace of $\Spec(A)$ for such $A$. 

For a ring $A$, let $\max(A)$ be the set of all maximal two-sided ideals of $A$. 

\begin{lem}\label{fin} Let $R$ be a commutative ring and let $A$ be an $R$-algebra such that $A$ is a finitely generated $R$-module.

(1) Let $\frak q\in \Spec(A)$ and let $\frak p\in \Spec(R)$ be the image of $\frak q$. Then $\frak q\in \max(A)$ if and only if $\frak p\in \max(R)$.  

(2) Let $I$ be a two-sided ideal of $A$ and let $I'\subset R$ be the inverse image of $I$ under $R\to A$. Then the image of $V(I)\subset \Spec(A)$ under $\Spec(A)\to \Spec(R)$ is $V(I')$.

(3) The map $\Spec(A)\to \Spec(R)$ is a closed map.

(4) If $R\to A$ is injective,  the map $\Spec(A) \to \Spec(R)$ is surjective.

\end{lem}

\begin{pf} (1) First assume $\frak q\in \max(A)$. We prove that $R/\frak p$ is a field. Let $a$ be a non-zero element of $R/\frak p$. Since the center of $A/\frak q$ is a field, $a$ is invertible in $A/\frak q$. Since the inverse $b$ of $a$ in $A/\frak q$ is integral over $R/\frak p$ and satisfies $\sum_{i=0}^n c_ib^i=0$ for some $c_i\in R/\frak p$ with $c_n=1$, we have $b= -\sum_{i=0}^{n-1} c_ia^{n-1-i}\in R/\frak p$. 

 Next assume $\frak p\in \max(R)$. Then $A/\frak p A$ is a finite-dimensional algebra over a commutative field $R/\frak p$. Hence every prime ideal of $A/\frak p A$ is a maximal two-sided ideal. Hence $\frak q\in \max(A)$. 

(2) 
It is clear that the map $\Spec(A)\to \Spec(R)$ sends $V(I)$ to $ V(I')$. Let $\frak p\in V(I')$. Since $R/I'\subset A/I$, $A/I\otimes_R R_{\frak p}\supset R/I'\otimes_R R_{\frak p}=(R/I')_{\frak p} \neq 0$. Take a maximal two-sided ideal $\m$ of this non-zero ring $A/I \otimes_R R_{\frak p}$. By (1), the inverse image of $\m$  in $(R/I')_{\frak p}$ is the maximal ideal of $(R/I')_{\frak p}$. Let $\frak q$ be the prime ideal of $A$ defined to be the inverse image of $\frak m$.  Then $\frak q\in V(I)$ and $\frak q\mapsto \frak p$ under $\Spec(A)\to \Spec(R)$. 

(3) follows from (2). 

(4) follows from the case $I=0$ of (2). 
\end{pf}

\begin{para}\label{idem} Let $A$ be a ring such that there is a commutative local ring $R$ and $A$ is an $R$-algebra and $A$ is finitely generated as an $R$-module.  

By \ref{fin} (1),  $\max(A)$ is a finite set. 

Let $\m\in \max(A)$.  By an idempotent for $(A, \m)$, we mean an element $e$ of $A$ such that $e^2=e$, $e\notin \m$, and $e\in \m'$ for all $\m' \in \max(A)\smallsetminus \{\m\}$.

\end{para}

\begin{prop}\label{idem2} Let $A$, $R$, $\m$ be as in \ref{idem}.

(1) Assume $R$ is excellent and henselian.  Then there is an idempotent for $(A, \m)$.

(2) Let $e$ be an idempotent for $(A, \m)$.  Let $I$ be a two-sided  ideal of $A$. Then $I\subset \m$ if and only if $e\notin I$. 

(3) Let $e$ be an idempotent for $(A, \m)$. Then $D(e)= \{\frak p\in \Spec(A)\;|\; \frak p\subset \m\}$. In particular, $D(e)$ is the smallest neighborhood of $\m$ in $\Spec(A)$. 
\end{prop}

\begin{pf} Let $\m_R$ be the maximal ideal of $R$.

(1) We first consider the case $R$ is a complete Noetherian local ring. In this case, $A=\varprojlim_n A/\m_R^nA$. 
Let $J=\Ker(A \to \prod_{\frak m} A/\frak m)$, where $\frak m$ ranges over all maximal two-sided ideals of $A$. Then $J^n\subset m_RA$ for some $n\geq 1$ and hence $A\overset{\cong}\to \varprojlim_n A/J^n$. We have an element $e_1\in A/J$ whose image in $A/\m$ is $1$  and  whose image in $A/\m'$ is $0$ for every $\m\in \max(A)\smallsetminus \{\m\}$. We can lift an idempotent $e_n\in A/J^n$ to an idempotent $e_{n+1}\in A/J^{n+1}$ by the standard method.
(If $\tilde e_n$ denotes a lifting of $e_n$ to $A$ and if $\tilde e_n^2=\tilde e_n+ x$ with $x\in J^n$, then for $\tilde e_{n+1}:=\tilde e_n+ (1-2\tilde e_n)x$, we have $\tilde e_{n+1}^2  \equiv  \tilde e_{n+1} \bmod J^{n+1}$.)

We next reduce the proof of (1) to the case $R$ is complete. Let $\hat R=\varprojlim_n R/\m_R^n$ be the completion of $R$. Then an idempotent $\hat e$  for $(A \otimes_R \hat R \m )$ exists. There is a finitely generated subring $R'$ of $\hat R$ over $R$ such that $\hat e$ comes from an idempotent  $e'$ of $A':= A\otimes_R R'$ satisfying $e'-1 \in \m A'$ and 
$e'\in \m' A'$  for all $\m'\in \max(A)\smallsetminus \{\m\}$. By Artin's approximation theorem generalized by Popescu \cite{Po}, there is an $R$-homomorphism $R'\to R$. Let $e$ be the image of $e'\in A'$ under $A'=A\otimes_R R'\to A$. Then $e$ is an idempotent for $(A, \m)$.

(2) If $I\subset \m$, then $e\notin I$ clearly. Assume $I$ is not contained in $\m$. We prove $e\in I$. Let $B=A/I$.  Then the image $\bar e$ of $e$ in $B$ is contained in all maximal two-sided ideals of $B$. Since the kernel of $B/\m_RB \to \prod_{\m'}  B/\m'$, where $\m'$ ranges over all maximal two-sided ideals of $B$, is a nil-ideal, we have $\bar e^n \in \m_RB$ for some $n\geq 1$. Hence $\bar e=\bar e^n\in \m_RB$. Let $P=B\bar e\subset B$. Then we have a homomorphism 
$B\to P\;;\; b\mapsto b\bar e$ and the  composition $P\to B \to P$ is the identity map. The induced homomorphisms $P/\m_RP \to B/\m_RB\to P/\m_RP$ have the properties  that  the first map is the zero map and the composition is the identity map. Hence $P/\m_RP=0$. By Nakayama's lemma, we have $P=0$. Hence $\bar e=0$.

(3) Straightforwards from (2). 
\end{pf}

\begin{para}\label{cU(x)} We will often consider the following situation.

Let $R$ be an excellent  strict local ring, let $S=\Spec(R)$, let $A$ be an $R$-algebra which is finitely generated as an $R$-module, and let $\m \in \max(A)$. Denote $\m$ also by $x$. We denote the open subspace $\{\frak p\in \Spec(A)\;|\; \m \supset \frak p\}$ of $\Spec(A)$ (\ref{idem2} (3)) by $\cU(x)$. This is  the smallest neighborhood of $x$ in $\Spec(A)$. 

\end{para}

\begin{para}\label{bfF}  The condition ({\bf F}) on a non-commutative scheme $X$.

 ({\bf F}): Locally on $X$, there are an excellent commutative ring $R$ and an $R$-algebra $A$ which is finitely generated as an $R$-module, and an open immersion $X\to \Spec(A)$.

 The condition  ({\bf F}$_T$) on a non-commutative scheme over $T$, where $T$ is an excellent scheme. 

({\bf F}$_T$)  Locally on $X$, there are a commutative ring $R$ over $T$ which is of finite type over $T$, an $R$-algebra $A$ which is finitely generated as an $R$-module, and an open immersion $X\to \Spec(A)$ over $T$. 

\end{para}

The following \ref{Fres} is clear. 

\begin{lem}\label{Fres} Let $X$ be a non-commutative scheme satisfying ({\bf F}) and let $x\in X$. Then $\cO_{X,x}$ is a finitely generated module over the local excellent ring $Z(\cO_{X,x})$, and $\kappa(x):=\cO_{X,x}/\frak p(x)$ is a finite-dimensional simple algebra over the residue field of $Z(\cO_{X,x})$. 

\end{lem}

\begin{prop}\label{Fprop} Let $f:X\to Y$ be a morphism of non-commutative schemes. Assume $X$ and $Y$ satisfy ({\bf F}). 
 Assume $f$ is surjective. Let $g:Y'\to Y$ be a morphism of non-commutative schemes and assume that either $f$ or $g$ is l.f.p. Then the morphism $X\otimes_Y Y' \to Y'$ is surjective. 

\end{prop}

\begin{pf} 

Let $y'\in Y'$. We prove that there is an element of $X\otimes_Y Y'$ whose image in  $Y'$ is $y'$. 
Let $y$ the image of $y'$ in $Y$ and let $x$ be an element $X$ whose image in $Y$ is $y$. Let $A=\cO_{Y, y}/\frak p(y)$, $B= \cO_{X, x}/\frak p(x)$, $C=\cO_{Y',y'}/\frak p(y')$.  Replacing $Y$ with $\Spec(A)$, $X$ with $\Spec(B)$, and $Y'$ with $\Spec(C)$, it is sufficient to prove that $\Spec(B\otimes_A C)\to \Spec(C)$ is surjective. By \ref{ffus}, it is sufficient to prove that $C\to B\otimes_A C$ is faithfully right flat. Hence it is sufficient to prove that $A\to B$ is faithfully right flat. We have $A\cong M_n(D)$ for some $n\geq 1$ and for some finite-dimensional central division algebra $D$ over a commutative field $k$. By \ref{Az}, the $A$-algebra $B$ is written as $M_n(D \otimes_k B')$ for some non-zero $k$-algebra $B'$ and hence is faithfully right flat over $A$. 
  \end{pf}
  
  \begin{para}\label{Noether} Let $X$ be a quasi-compact non-commutative scheme satisfying ({\bf F}). Then $X$ is an Noetherian space  (that is, every open subset of $X$ is quasi-compact). In fact,  $X$ is locally an open subspace of $\Spec(A)$ for some algebra $A$ such that $Z(A)$ is a Noetherian ring and $A$ is finitely generated as a $Z(A)$-module. Since all two-sided ideals of $A$ are finitely generated, $\Spec(A)$ is a Noetherian space. 
  
  Hence $X$ is a finite union of irreducible closed subsets. 
  
  \end{para}
  \begin{prop}\label{a1} Let $X\to Y$ be an  l.f.p.  morphism of non-commutative schemes. Assume $X$ and  $Y$ satisfy ({\bf F}) and $X$ and $Y$ are  quasi-compact. Then $f(X)$ is a constructible subset of $Y$.

    \end{prop}
  
  \begin{pf} 
  By Noetherian induction, we may assume that $X$ is a prime non-commutative scheme and that $f(C)$ is constructible for every closed subset $C\neq X$ of $X$. In fact, if $f(X)$ is not constructible, since $X$ is a Noetherian space, there is a minimal element in the set of all closed subsets $C$ of $X$ such that $f(C)$ are not constructible. If this minimal element $C$ is a finite union of closed subsets $C_i$, we have $C=C_i$ for some $i$ because otherwise, all $f(C_i)$ would be constructible and hence $f(C)$ would be constructible. Thus this minimal element $C$ is irreducible. It is sufficient to consider the case $X$ is this $C$.

 We may assume that $Y$ is the closure of the image of the generic point of $X$ in $Y$, and we may assume that $Y$ is a prime non-commutative scheme. If $X$ is  union of open subsets $U$ and $U'$ such that all $f(U)$ and $f(U')$ are constructible, then $f(X)=  f(U)\cup f(U')$ is constructible. 
  By this, we may assume that $Y$ is an open subspace of $\Spec(A)$ for some prime ring $A$ such that $Z(A)$ is Noetherian and $A$ is finitely generated as a $Z(A)$-module, 
$X$ is an open subspace of $\Spec(B)$ for some prime ring $B$ such that $Z(B)$ is Noetherian and $B$ is finitely generated as a $Z(B)$-module, 
    $B$ is an $A$-algebra of finite presentation, $A\to B$ is injective, and the morphism $\Spec(B)\to \Spec(A)$ is compatible with $f:X\to Y$. Note that then $Z(A)$ and $Z(B)$ are integral domains.  There are a non-zero element $a$ of $Z(A)$ and a non-zero element $b$ of $Z(B)$ such that $b/a\in Z(B)$ and such that $A[1/a]$ is an Azumaya algebra over $Z(A)[1/a]$, $V:= Y\times_{\Spec(Z(A))} \Spec(Z(A)[1/a])$ coincides with $\Spec(A[1/a])$, $B[1/b]$ is an Azumaya algebra over $Z(B)[1/b]$, and $U:= X \times_{\Spec(Z(B))} \Spec(Z(B)[1/b])$ coincides with $ \Spec(B[1/b])$. Then via the isomorphisms of topological spaces $V\overset{\cong}\to \Spec(Z(A)[1/a])$ and $U\overset{\cong}\to \Spec(Z(B)[1/b])$, the constructibility of $f(U)$ in $V$ is reduced to the constructibility of the image of the morphism $\Spec(Z(B)[1/b]) \to \Spec(Z(A)[1/a])$ of finite presentation in the theory of schemes. Since $f(X\smallsetminus U)$ is constructible by our Noetherian induction, we have that $f(X)=f(U)\cap f(X\smallsetminus U)$ is constructible. 
   \end{pf}

  \begin{prop}\label{a2} Let $X\to Y$ be a right flat morphism of non-commutative schemes. Assume both $X$ and $Y$ satisfy ({\bf F}). Then $f(X)$ is stable in $Y$ under generalizations.

  \end{prop}
  
  \begin{pf} 
 Let $x\in X$ with image $y$ in $Y$. Let $y'\in Y$ be a generalization of $y$ (that is, $y$ belongs to the closure of $y'$ in $Y$). We prove that there is an element $x'$ of $X$ with image $y'$ in $Y$.
 We may assume that $Y$ is an open subspace of $\Spec(A)$, where $A$ is an $R$-algebra over a commutative ring $R$ such that $A$ is finitely generated as an $R$-module.  Replacing $A$ by $A/\frak p$ where $\frak p$ is the prime ideal of $A$ given by $y'$ and replacing $Y$ by $Y\times_{\Spec(A)} \Spec(A/\frak p)$, we may assume $A$ is a prime ring and $y'$ corresponds to the prime ideal $(0)$ of $A$. By replacing $R$ by the image of $R$ in $A$, we may assume $R$ is an integral domain. 
 
 We may assume that $X$ is an open subscheme of $\Spec(B)$ where $B$ is an algebra  over an excellent commutative ring $R'$
  and is finitely generated as an $R'$-module. By replacing $R'$ with the strict henselization at the image of 
  $x$ in $\Spec(R')$, we may assume that $X=\cU(x)$ (\ref{cU(x)}). We prove that there is a point of $\cU(x)$ lying over the prime ideal 
  $(0)$ of 
  $A$, that is, lying over the prime ideal $(0)$ of $R$. If there is no such element, 
   $ \cU(x)\times_{\Spec(R)} \Spec(Q(R))$ is empty, where $Q(R)$ denotes 
  the field of fractions of $R$. Since every element in the intersection of all prime ideals of a ring is nilpotent (\cite{La0} Theorem 10.7),
   the image of $e$ in $Q(R)\otimes_R A$ is nilpotent and hence is $0$ because it is an idempotent. Hence there is a non-zero element $a$ of 
   $R$ such that $ae=0$. Since $A\to A\;;\;x\mapsto ax$ is injective and the direct summand $eB$ of $B$ is right flat over $A$, $eB\mapsto eB\;;\;x\mapsto ax$ 
   is injective. Since $ae=0$, we have   $e=0$, a contradiction. 
   \end{pf}
   
   \begin{prop}\label{a3} Let $X\to Y$ be an  l.f.p. right  flat morphism of non-commutative schemes. Assume both $X$ and $Y$ satisfy ({\bf F}).   Then $f$ is an open map. 
      \end{prop}

\begin{pf}
This follows from \ref{a1} and \ref{a2}. 
\end{pf}

\section{\'Etale cohomology}\label{sec5}

 In this section, we define the \'etale cohomology groups of non-commutative schemes imitating SGA 4 (\cite{SGA4}).

 We prove a finiteness theorem (\ref{fin2}) for non-commutative schemes satisfying ({\bf F}). We also prove a proper base change theorem (\ref{pbc2}) for non-commutative schemes satisfying ({\bf F}) and for relatively commutative morphisms between them.

\begin{para} Let $f: X\to Y$ be a morphism of non-commutative schemes. 

We say $f$ is \'etale if it satisfies the following two conditions.  

(i) $f$  is flat (\ref{flat2}).

(ii) $f$ is unramified (\ref{sep}, this tells that $f$ is r.c. and l.f.p). 
\end{para}

For schemes, this \'etaleness coincides with the usual \'etaleness.

\begin{lem}\label{lemet}  (1) If $X\to Y$ and $Y\to Z$ are  \'etale morphisms of non-commutative schemes, the composition $X\to Z$ is  \'etale. 

(2) Let  $X\to Y$ be an  \'etale morphism. Then for every morphism $Y'\to Y$ of non-commutative schemes, the induced morphism $X\times_Y Y'\to Y'$ is  \'etale.

(3) Let $X\to Z$ and $Y\to Z$ be \'etale morphisms of non-commutative schemes. Then every morphism $X\to Y$ over $Z$ is \'etale. 
\end{lem}

\begin{pf} 
(1) (resp. (2)) follows from \ref{lemflat} (1) (resp. (2))  for flat morphisms and from \ref{sep2} (1) (resp. (2)) for unramified morphisms.

(3)  We first prove that the morphism $X\to Y$ is flat. This morphism is the composition $X\to X\times_Z Y \to Y$. Here the first arrow is flat because it is the base change of the open immersion  $Y\to Y\times_Z Y$ by $X\times_Z Y \to Y\times_Z Y$. The second arrow is flat because it is the base change of the flat morphism $X\to Z$ by $Y\to Z$. The morphism  $X\to X\times_Y X$ is an open immersion because the composition $X\to X\times_Y X\to X\times_Z X$ is an open immersion and the morphism $X\times_Y X \to X\times_Z X$ is the base change of the open immersion $Y\to Y \times_Z Y$ by $X\times_Z X \to Y \times_Z Y$ and hence is an open immersion. 
\end{pf}

\begin{para} For a non-commutative scheme $X$, 
we define the  \'etale site $X_{et}$ as follows. 

An object of $X_{et}$ is a non-commutative scheme over $X$ which is \'etale over $X$. 

A morphism is a morphism of non-commutative schemes over $X$. (Note that a morphism in this site  is \'etale by \ref{lemet} (3).)
 
A morphism $f: Y\to Z$ in $X_{et}$ is a covering if it is universally surjective (that is, for every morphism $Z'\to Z$ of non-commutative schemes, the map $Y':=Y\times_Z Z'\to Z'$ induced by $f$ is surjective).

\end{para}

\begin{rem}\label{remsurj} By \ref{Fprop}, if  $X$ satisfies ({\bf F}) and $f: Y\to Z$ is a morphism in $X_{et}$, $f$ is a covering if it is surjective.   The author does not know any example of an r.c. l.f.p. morphism of non-commutative schemes which is surjective but not universally surjective.
\end{rem}

\begin{para} Functoriality. If $X\to Y$ is a morphism of non-commutative schemes, $U\mapsto X \times_Y U$ sends an object of $Y_{et}$ to an object of $X_{et}$  and a covering to a covering. 
Hence we have a morphism of the associated topoi $\tilde X_{et}\to \tilde Y_{et}$. 

\end{para}

\begin{prop}\label{Azet} Let $S$ be a scheme, let $\cB$ be an Azumaya algebra over $\cO_S$ on $S$, and let $X=\Spec(\cB)$. Then the pullback functor $X\times_S$ gives an equivalence of sites between $S_{et}$ and $X_{et}$. \end{prop}

\begin{pf} This is deduced  form \ref{AzX}.
\end{pf}

\begin{para}\label{Oet}  For a non-commutative scheme $X$, the presheaves $U\mapsto \cO(U)$ and $U\mapsto \Gamma(U, Z(\cO_U))$ on the \'etale site $X_{et}$ need not be a sheaf. Note that $\Gamma(U, Z(\cO_U))$ is understood as the set of morphisms  of non-commutative schemes $U\to \Spec(\Z[T])$ (\ref{A[T]}, \ref{specspi2}). Thus a representable functor is not necessarily a sheaf for the \'etale topology. 
 
  Let $\cO_{X_{et}}$ be the sheafification of the presheaf $U\mapsto \cO(U)$ on $X_{et}$. See \ref{Oet2} for an example of $X$ such that $\cO(X)=\Gamma(X, \cO_X)$ does not coincide with $\Gamma(X_{et}, \cO_{X_{et}})$ and $\Gamma(X, Z(\cO_X))$ does not coincide with $\Gamma(X_{et}, Z(\cO_{X_{et}}))$.
See also \ref{Oet4} (1). 
\end{para}

\begin{para}\label{lim0}  Let $X$ be a quasi-compact non-commutative scheme which satisfies ({\bf F}). 

Then $X$ is a Noetherian space 
(\ref{Noether}). 
By this fact and  by the fact that every morphism between objects of $X_{et}$ is an open map by \ref{a3},  and by the fact every surjective morphism between objects of $X_{et}$ is a covering (\ref{remsurj}), we have that the category of sheaves on $X_{et}$ is a Noetherian topos in the sense of \cite{SGA4} vol. 2, Exp. VI, 2.11. Hence by loc.cit Theorem 5.1 and Theorem. 8.7.3,  we have the following \ref{lim1} and its generalization \ref{lim3}.

\end{para}

\begin{prop}\label{lim1} Let $X$ be a quasi-compact non-commutative scheme satisfying ({\bf F}), and let $(\cF_{\lam})_{\lam}$ be a directed family of sheaves of abelian groups on $X_{et}$. Then for every $m$, we have an isomorphism 
$$\varinjlim_{\lam}\; H^m(X_{et}, \cF_{\lam}) \overset{\cong}\to H^m(X_{et}, \varinjlim_{\lam} \cF_{\lam}).$$

\end{prop}

\begin{para}\label{lim2} This is a preparation for \ref{lim3}.  Let $\La$ be a directed ordered set and assume that we are given an inverse system $(X_{\la})_{\la\in \La}$ of non-commutative schemes. Let $I$ be a finite set. Assume that for each $\la\in \La$, we have an open covering $(U_{\la,i})_{i\in I}$ of $X_{\la}$. Assume that for each $i\in I$, we have an inductive system  $(A_{\la,i})_{\la\in \La}$ in the category (Alg). Assume that we have an open immersion $U_{\la,i}\overset{\subset}\to \Spec(A_{\la,i})$ for each $\la\in \La$ and $i\in I$. We assume the following (i)--(iii).

(i) For $\la,\mu\in \La$ such that $\la\leq \mu$ and for each $i\in I$, the inverse image of $U_{\la,i}\subset X_{\la}$ in $X_{\mu}$ is $U_{\mu,i}$, the inverse image of $U_{\la,i}\subset \Spec(A_{\la,i})$ in $\Spec(A_{\la,i})$ in $\Spec(A_{\mu,i})$ is $U_{\mu, i}$,  and the morphism $U_{\mu,i}\to U_{\la,i}$ induced by $\Spec(A_{\mu,i})\to \Spec(A_{\la, i})$ coincides with the morphism induced by the transition morphism $p_{\mu,\la}: X_{\mu}\to X_{\la}$. 

(ii) For each $\la\in \La$, the center $Z(A_{\la,i})$ is an excellent ring and $A_{\la,i}$ is finitely generated as a $Z(A_{\la,i})$-module.

(iii) For each $i\in I$, if we write $A_i=\varinjlim_{\la}\; A_{\la,i}$, then $Z(A)$ is an excellent ring and $A_i$ is finitely generated as a $Z(A_i)$-module. 

Let  $X= \varprojlim_{\la}\; X_{\la}$ as a topological space, endow $X$ with the sheaf of rings $\varinjlim_{\la} \;p^{-1}_{\la}(\cO_{X_{\la}})$ where $p_{\la}: X\to X_{\la}$ is the projection, and define $\frak p(x):= \varinjlim_{\la} \frak p(x_{\la})$ where $x_{\la}$ is the image of $x$ in $X_{\la}$. Then $X= (X, \cO_X, (\frak p(x))_{x\in X})$ is a non-commutative scheme.

Assume that for each $\lam\in \La$, we are given a sheaf $\cF_{\lam}$ of abelian groups on $X_{\lam,et}$ and assume that for each $\lam, \mu\in \La$ such that $\la\leq \mu$, we have a homomorphism $h_{\lam,\mu}: p_{\mu, \lam}^{-1}(\cF_{\la})\to \cF_{\mu}$. We assume that $h_{\la,\la}$ is the identity homomorphism for $\la\in\La$ and that we have   $h_{\la, \nu}= h_{\mu,\nu}\circ p^{-1}_{\nu,\mu}(h_{\la, \mu})$ when $\la\leq \mu\leq \nu$. 

Let $\cF:=\varinjlim_{\la} p_{\la}^{-1}(\cF_{\la})$. 

\end{para}

\begin{lem}\label{lim3} For each $m$, we have an isomorphism
$$\varinjlim_{\la}\; H^m(X_{\la,et}, \cF_{\la}) \overset{\cong}\to H^m(X, \cF).$$

\end{lem}

\begin{prop}\label{small3} Let $R$, $S$, $A$, $x$, $\cU(x)$  be as in \ref{cU(x)}.  Let $U\to \cU(x)$ be an \'etale morphism such that there is $u\in U$ whose image in $\cU(x)$ is $x$.  Then there is an open immersion  $\cU(x)\to U$ over $\cU(x)$.
\end{prop}

\begin{pf}  (In this proof here,  we use the fact that $f$ is right flat and also left flat.)

 The integral closure $R_1$ of the image of $R$ in $A$ is a finitely generated $R$-module and hence it is a finite product of local rings because $R$ is henselian. Let $R_2$ be the component of $R_1$ in this direct product decomposition such that the image of $x$ in $\Spec(R_1)$ belongs to $\Spec(R_2)$. By replacing $R$ with $R_2$ and $A$ with $A\otimes_{R_1} R_2$, we may assume that $R$ is the center of $A$ and is a local ring. Then we have $A=\cO_{\cU(x), x}$. 
 
  Since $U\to \cU(x)$ is r.c. and l.f.p., 
we may assume that 
there are a finitely generated commutative ring $R'$ over $R$, an $A\otimes_R R'$-algebra $B$ such that $A\otimes_R R' \to B$ is surjective, and an open immersion $U\to \Spec(B)$ over $A$. Let $V$ be the set of points of $\Spec(R')$ at which the morphism $\Spec(R')\to \Spec(R)$ is quasi-finite. Then $V$ is open in $\Spec(R')$ by \cite{EGA4.3} Cor. 13.1.4. Furthermore, the image $t$ of $u$ in $\Spec(R')$ is contained in $V$ by \ref{unrlem}. By replacing $\Spec(R')$ by an affine open neighborhood of $t$  in $V$, we may assume that $\Spec(R')$ is quasi-finite over $S$. Then since $R$ is henselian, there is an open and closed neighborhood of $t$ in $\Spec(R')$ which is finite over $S$. By replacing $\Spec(R')$ by this open an closed subset, we may assume that $R'$ is  finitely generated as an $R$-module. Hence $B$ is  finitely generated as  a  left $A$-module.  Let $R'_1$ be the integral closure of the image of $R'$ in $B$. Since $R'$ is henselian, $R'_1$ is a finite product of local rings. Let $R'_2$ be the component of $R'_1$ in this direct product decomposition such that the image of $u$ in $\Spec(R'_1)$ is contained in $\Spec(R'_2)$. By replacing $R'$ with $R'_2$ and replacing $B$ with $B\otimes_{R'_1} R'_2$, we may assume that  $R'$ is the center of $B$ and is a local ring. We have $B=\cO_{U,u}$. 
By \ref{unrlem} and by the subjectivity of $A\otimes_R R' \to B$, the maps $A/\m \to B/B\m= B/\m B$ are bijective for all $\m\in \max(A)$. Hence the map $A\to B$ is surjective by Nakayama's lemma. For the proof of \ref{small3}, it is now sufficient to prove $A\overset{\cong}\to B$. In fact, then  we will have $U\overset{\cong}\to \cU(x)$ because $\cU(x)$ is the smallest neighborhood of $x$ in $\Spec(A)$.

 Let $I$ be the kernel of the surjection $A\to B$.  We prove $I=0$. Since $B$ is  left flat, by \cite{La} Chap. 2, \S4D Proposition 4.30, $B$ is projective as a left $A$-module and hence  $A\cong I \oplus B$ as a left $A$-module.

{\bf Claim}. If $M$ is a subquotient of the left $A$-module $I$ and $N$ is a subquotient of the left $A$-module $B$, and if $M\cong N$ as a left $A$-module, then $M=0$ and $N=0$. 

Proof of Claim.  Since $B$ is right flat over $A$, $B\otimes_A I=0$ and hence $B\otimes_A M=0$. On the other hand, $B\otimes_A N=N$. Hence $N=0$ and hence $M=0$.

Now by right multiplication, the opposite ring $A^{\circ}$ is isomorphic to the endomorphism ring of the left $A$-module $A\cong I \oplus B$. By Claim, this endomorphism ring is  the direct product (as a ring)  of the endomorphism ring of the left $A$-module $I$ and the endomorphism ring of the left $A$-module $B$. Since the center of $A$ is a local ring, the ring $A$ can not be a product of two non-zero rings.  Hence the endomorphism ring of the left $A$-module $I$ is zero, that is,  the right multiplication by $1$ on $I$ is the zero map. Hence $I=0$.
\end{pf}

\begin{para}\label{barx} Let $X$ be  a non-commutative scheme satisfying ({\bf F}). We consider stalks of a sheaf $\cF$ on $X_{et}$. 

Let  $x\in X$.
Then   $\kappa(x):=\cO_{X,x}/\frak p(x)$ is a finite-dimensional central simple algebra over its center $k$ which is a commutative field. Let $\kappa(\bar x)= \kappa(x) \otimes_k \bar k$, where $\bar k$ is the separable closure of $k$, and let $\bar x=\Spec(\kappa(\bar x))$. We have a canonical morphism $i_{\bar x}: \bar x \to X$ of non-commutative schemes. The topos of sheaves on $\bar x_{et}$ is equivalent to the category of sets via $\cG\mapsto \cG(\bar x)$. This is by the fact that $\kappa(\bar x)\cong M_n(\bar k)$ for some $n\geq 1$ and by \ref{Azet}.  

For a sheave $\cF$ on $X_{et}$, let $\cF_{\bar x}=(i_{\bar x}^{-1}\cF)(\bar x)$. In the case $k$ is separably closed, we denote $\cF_{\bar x}$ also by $\cF_x$. 
\end{para}

 \begin{lem}\label{stalk} Assume $X$ satisfies ({\bf F}).
 
(1)  For a morphism $h:\cF\to \cG$ of sheaves on $X_{et}$, $h$ is an isomorphism if and only if the induced maps $\cF_{\bar x} \to \cG_{\bar x}$ are bijective for all $x\in X$. 
 
 (2) A sequence 
$\cF'\to \cF \to \cF''$ of sheaves of abelian groups on $X_{et}$ is exact if and only if $\cF'_{\bar x}\to \cF_{\bar x}\to \cF''_{\bar x}$ are exact for all $x\in X$.

\end{lem}

\begin{pf} (1) Assume that $\cF_{\bar x}\to \cG_{\bar x}$ are bijective for all $x\in X$. We prove that $h$ is injective (resp. surjective).  Let $U\in X_{et}$. Let $a_1, a_2\in \cF(U)$ and assume that there images in $\cG(U)$ are the same (resp. let $a\in \cG(U)$). Let $u\in U$, and let $x\in X$ be the image of $u$ in $X$. Then the morphism $\bar x\to x$ factors as $\bar x\to u \to x$. Hence by the assumption, there are $U_u \in U_{et}$, and a morphism $\bar x\to U_u$ over $U$ such that the  image of $\bar x$ in $U$ is $u$ and such that the images of $a_1, a_2$ in $\cG(U_u)$ are the same (resp. there is  $b\in \cF(U_u)$ whose image in $\cG(U_u)$ is the pullback of $a$). The morphism $\coprod_{u\in U} U_u \to U$ is surjective and hence universally surjective by \ref{Fprop}, and hence $\coprod_{u\in U} U_u\to U$ is a covering in $X_{et}$. This proves that $h$ is injective (resp. surjective).

 (2) follows from (1). 
\end{pf}

\begin{prop}\label{small4} Let $S$, $A$, $x$, $\cU(x)$ be as in \ref{cU(x)} and let $s$ be the closed point of $S$. Let $\cF$ be a sheaf of abelian groups on $\cU(x)_{et}$. 

(1)  We have $R\Gamma(\cU(x)_{et}, \cF)\overset{\cong}\to \cF_x$.

(2) Let $\epsilon: \cU(x)\to S$ be the canonical morphism. Then $(R\epsilon_*\cF)_s=(\epsilon_*\cF)_s=\cF_x$.

\end{prop}
\begin{pf} (1) follows from \ref{small3}. (2) follows from (1). 
\end{pf}

\begin{lem}\label{imet1}

Let $i: Y\to X$ be an l.f.p. closed immersion of non-commutative schemes. Assume $X$ satisfies ({\bf F}).

(1) The functor $i_*$ is exact for sheaves $\cF$ of abelian groups on $Y_{et}$. For $x\in X$, we have $(i_*\cF)_{\bar x}=\cF_{\bar x}$ if $x\in Y$, and $(i_*\cF)_{\bar x}=0$ if $x\notin Y$. 

(2)  The functor $i_*$ gives  an equivalence from the category of  sheaves of abelian groups on $Y_{et}$ to  the category of sheaves of abelian groups on $X_{et}$ whose restrictions to $U=X\smallsetminus Y$ are zero. The converse functor is given by $i^{-1}$. In particular, the former category depends only on the closed set $Y$ of $X$, not on the non-commutative scheme structure of $Y$. 
 \end{lem}

\begin{pf}
(1)  follows from \ref{small4}. 

(2) If $\cF$ is a sheaf of abelian groups on $Y_{et}$, $i^{-1}i_*\cF \to \cF$ is an isomorphism as is seen by checking the stalks. If $\cG$ is a sheaf of abelian groups on $X_{et}$ such that $\cG|_U=0$, $\cG\to i_*i^{-1}\cG$ is an isomorphism as is seen by checking the stalks.
\end{pf}
\begin{lem}\label{imet2}  Let $j:U\to X$ be an open immersion of non-commutative schemes. Assume $X$ satisfies ({\bf F}).

(1) The inverse image functor $j^{-1}$ from the category of sheaves of abelian groups on $X_{et}$ to the category of abelian groups on $U_{et}$ has a left adjoint functor $j_!$.

(2) $j_!$ is an exact  functor for sheaves $\cF$ of abelian groups on $U_{et}$. For $x\in X$, we have $(j_!\cF)_{\bar x}=\cF_{\bar x}$ if $x\in U$, and $(j_!\cF)_{\bar x}=0$ if $x\notin U$.

(3) If $i:Y\to X$ is an l.f.p. closed immersion and if $j:U\to X$ is the complement, 
we have an exact sequence $0\to j_!j^{-1}F\to F \to i_*i^{-1}F \to 0$. 

\end{lem}

\begin{pf} For a sheaf of abelian groups $\cF$ on $U_{et}$, define $j_!\cF$ to be the sheaf associated to the following  presheaf $\cG$ on $X_{et}$. For $V\in X_{et}$, $\cG(V)$ is $\cF(V)$ if the image of $V$ in $X$ is contained in $U$ (then $V\to X$ factors uniquely as $V\to U\to X$ and $V$ is regarded as an object of $U_{et}$ and hence $\cF(V)$ is defined)  and $\cG(V)=0$ otherwise. (1) and (2) follow from this definition. (3) is shown by checking the stalks by \ref{stalk} (2). 
\end{pf}

\begin{para}  Let $f: Y\to X$ be an  l.f.p. immersion of non-commutative schemes. Assume $X$ satisfies ({\bf F}). We define $f_!$. 

For a sheaf $\cF$ of abelian groups on $Y_{et}$, by using a factorization 
$Y\overset{i}\to U \overset{j}\to X$ where $i$ is an l.f.p. closed immersion  and $j$ is an open immersion, we define $f_!\cF:=j_!i_*\cF$. This is independent  of the choice of the factorization. In fact, if $Y\to U'\to Y$ is another factorization, we have the third  factorization $Y\overset{k}\to U\cap U'\overset{l}\to X$, and we have a canonical morphism $l_!k_*(\cF) \to j_!i_*(\cF)$. 
We see that this is an isomorphism by comparing the stalks. 

\end{para}

\begin{para}\label{const1} Let $X$ be a non-commutative scheme satisfying ({\bf F}).

We say sheaf $\cF$ on $X_{et}$ is  constructible if locally on $X$, there is a finite family of l.f.p. immersions $Y_i\to X$ such that the set $X$ is the disjoint union of these $Y_i$ and such that the pullback of $\cF$ to each $Y_i$ is locally constant and finite. 

By \ref{imet2} (3), we see that a sheaf $\cF$ of abelian groups on $X_{et}$ is constructible if and only if locally on $X$, there is a finite filtration on $\cF$ whose each graded quotient is isomorphic to $f_!\cG$ for some l.f.p. immersion $f:Y\to X$ and for some locally constant finite sheaf $\cG$ of abelian groups on $Y_{et}$.

\end{para}

\begin{lem}\label{const2} Let $X$ be as in \ref{const1}. 

(1) For a homomorphism $h: \cF\to \cG$ between constructible (resp. locally constant and finite) sheaves of abelian groups on $X_{et}$, the kernel and the cokernel of $h$ are smooth (resp. constructible).

(2) For an exact sequence $0\to \cF \to \cG \to \cH \to 0$ of sheaves of abelian groups on $X_{et}$ such that $\cF$ and $\cH$ are constructible (resp. locally constant and finite),  $\cG$ is also constructible (resp. locally constant and finite)

\end{lem}

\begin{pf} These are proved as in the case of schemes treated in Lemma 2.1 and Proposition 2.6 in \cite{SGA4} vol. 3, Exp. IX.  \end{pf}

We consider a non-commutative version of the finiteness theorem of Gabber (\cite{ILO} Exp. XIII, Exp. XXI; it is a generalization the finiteness theorem of Deligne in the chapter on finiteness in \cite{De2}). 

\begin{thm}\label{fin2}  Let $X$ and $Y$ be quasi-compact non-commutative schemes satisfying ({\bf F}), and let $f:X\to Y$ be an l.f.p. morphism.

Let $\cF$ be a constructible sheaf of abelian groups on $X_{et}$ such that the orders of all stalks of $\cF$ are invertible on $Y$.  Then the sheaf  $R^mf_*\cF$ on $Y_{et}$  is constructible  for every $m$ and is zero for $m\gg 0$. 

(2) Let $\cF$ be a constructible sheaf of sets on $X_{et}$.  Then the sheaf  $f_*\cF$ on $Y_{et}$  is constructible.

(3) Let $\cF$ be a constructible sheaf of  groups on $X_{et}$ such that the orders of all stalks of $\cF$ are invertible on $Y$.  Then the sheaf  $R^1f_*\cF$ on $Y_{et}$  is constructible.

\end{thm}

The proof of \ref{fin2} will be completed in \ref{pffin2}.

\begin{para}\label{spcase1} To prove \ref{fin2}, we first consider the following situation. 

Let $A$ be a ring and let $B$ be an $A$-algebra of finite presentation. Assume that $A$ and $B$ are prime rings, that they are finitely generated as modules over their centers, and that the map $A\to B$ is injective. Note that $Z(A)$ and $Z(B)$ are integral domains. We assume that $Z(A)$ and $Z(B)$ are excellent rings and $Z(A)$ is regular. 

Let $Y$ be an open subspace of $\Spec(B)$, let $X$ be an open subspace of $\Spec(B)$, and assume that $X$ is contained in the inverse image of $Y$ in $\Spec(B)$. 

Let $a$ be a non-zero element of $Z(A)$, let $b$ be a non-zero element of $Z(B)$ such that $b/a\in Z(B)$, 
$A[1/a]$ is an Azumaya algebra over $Z(A)[1/a]$, $B[1/b]$ is an Azumaya algebra over $Z(B)[1/b]$, $Y\times_{\Spec(Z(A))} \Spec(Z(A)[1/a])=\Spec(A[1/a])$, and $X\times_{\Spec(Z(B)} \Spec(Z(B)[1/b])= \Spec(B[1/b])$. 

Let $f: X\to Y$, $\pi:Y\to \Spec(Z(A))$, $j: U:= X \times_{\Spec(Z(B))} \Spec(Z(B)[1/b]) = \Spec(B[1/b]) \to X$, $g=\pi\circ f\circ j: U
 \to Y$,  $h: V:= \Spec(Z(B)[1/b])\to \Spec(Z(A))$ 
 be the canonical morphisms.  Note that $j$ is an open immersion. 

By \ref{Azet}, the topos of sheaves on $U_{et}$ and the topos of sheaves on $V_{et}$ are equivalent. We will identify a sheaf on $U_{et}$ with the corresponding sheaf on $V_{et}$.

\end{para}

We use the following elementary lemma. 

\begin{lem}\label{D(e)big}  Let $A$ be an Azumaya algebra over an integral domain $R$, and let $e$ be a non-zero idempotent of $A$. Then $D(e)=\{\frak p\in \Spec(A)\;|\; e\notin \frak p\}$ coincides with the whole  $\Spec(A)$. 

\end{lem}

\begin{pf} Let $Q(R)$ be the field of fractions of $R$. Then $Ae$ is a direct summand of the $R$-module $A$, and $Q(R)\otimes_R Ae\neq 0$. Hence $Ae$ is a finitely generated projective $R$-module rank $\geq 1$. Let $\frak p$ be a prime ideal of $A$. Then $\frak p=A\frak q$ for a prime ideal $\frak q$  of $R$. We have that $0\neq R/\frak q\otimes_R Ae\subset R/\frak q\otimes_R A$. This shows $e\notin \frak p$. 
\end{pf}

\begin{lem}\label{spcase2} Let the situation be as in \ref{spcase1}. 
Then for a sheaf $\cF$ of abelian groups (resp. sets, resp. groups) on $U_{et}$, the canonical morphism $\pi^{-1}R^mh_*\cF\to R^mg_*\cF$ is an isomorphism for every $m$ (resp. for $m=0$, resp. for $m=1$).

\end{lem}

\begin{pf}  Let $y\in Y$ and let $s$ be the image of $y$ in $S:=\Spec(Z(A))$. It is sufficient to prove that  the map  $(R^mh_*\cF)_{\bar s}\to (R^mg_*\cF)_{\bar y}$ is an isomorphism. Let $S_{\bar s}$ be the strict henselization of $S$ at $\bar s$. We consider the smallest neighborhood $\cU(\bar y)$ of $\bar y$ in $Y \times_S S_{\bar s}$ (\ref{cU(x)}). Since $R$ is regular, $\cO_{S, \bar s}$ is an integral domain. Hence by \ref{D(e)big} applied to the Azumaya algebra 
$A\otimes_{Z(A)} \cO_{S,\bar s}[1/a]$ over $\cO_{S,\bar s}[1/a]$, we have $\cU(\bar y)\times_S \Spec(Z(A)[1/a])=  Y \times_S \Spec(\cO_{S, \bar s}[1/a])$. From this, we have
$U\times_S S_{\bar s}=U\times_Y \cU(\bar y)$. Hence we have
$$(R^mg_*\cF)_{\bar y}= H^m((U\times_Y \cU(\bar y))_{et}, \cF)=H^m((U\times_S S_{\bar s})_{et}, \cF),$$ where the first $=$ follows from  \ref{small4} (1). 
On the other hand, $$(R^mh_*\cF)_{\bar s}= H^m((V \times_S S_{\bar s})_{et}, \cF).$$ By \ref{Azet}, we have 
$H^m((U\times_S S_{\bar s})_{et}, \cF)\cong H^m((V \times_S S_{\bar s})_{et}, \cF)$.
\end{pf}

\begin{lem}\label{spcase2.5} Let the situation be as in \ref{spcase1}. 
Then for a constructible sheaf $\cF$ of abelian groups  on $U_{et}$, the sheaf  $R^mg_*\cF$ is constructible for every $m$  and is zero for $m\gg 0$. For a constructible sheaf $\cF$ of sets (resp. groups) on $U_{et}$, $f_*\cF$ (resp. $R^1f_*\cF$) is constructible.

\end{lem}

\begin{pf} By the finite theorem of Gabber (\cite{ILO} Exp. XIII and Exp XXI), for a constructible sheaf $\cF$ of abelian groups $\cF$ on $V_{et}$, $R^mh_*\cF$ is constructible for every $m$ and is zero for $m\gg 0$, and for a constructible sheaf $\cF$ of sets (resp. groups) on $V_{et}$, $h_*\cF$ (resp. $R^1h_*\cF$) is constructible. This lemma follows from it by  \ref{spcase2}. 
\end{pf}
\begin{para}\label{pffin2} We prove \ref{fin2}.

For a sheaf $\cF$ on $X_{et}$ and for a closed subset $C$ of $X$, there is a  unique quotient sheaf $\cF$ of $\cF$ such that $\cF_{\bar x}\overset{\cong}\to (\cF_C)_{\bar x}$ if $x\in C$ and $(\cF_C)_{\bar x}$ is a one point set if $x \in X\smallsetminus C$. The uniqueness is clear. The existence of  $\cF_C$ is seen as follows. In the case  $C$ is irreducible, $C$ has the canonical structure of a prime  non-commutative scheme and we have the l.f.p. closed immersion $i:C\to X$. We have  $\cF_C:=i_*i^{-1}\cF$. In general, 
let $C_k$ ($1\leq k\leq r$) be all the irreducible components of $C$. Then we obtain $\cF_C$ as the image of $\cF\to \prod_{k=1}^r \cF_{C_k}$. If $\cF$ is a sheaf of  groups (resp. abelian groups), $\cF_C$ has the unique structure of a sheaf of  groups (resp. abelian groups) which is compatible with that of $\cF$. 

We say that $\cF$  has supports in  $C$ if $\cF=\cF_C$, that is, if the restriction of $\cF$ to $X\setminus C$ is the constant sheaf associated to a one point set. 

We prove (1) of \ref{fin2}. 
For closed subsets  $C$ and $C'$ of $X$, we have an exact sequence 

(1) $0\to \cF_{C\cup C'} \to \cF_C \oplus \cF_{C'} \to \cF_{C\cap C'}\to 0$.

Thus if \ref{fin2} (1) is true for $\cF$ with supports in $C$ and  for $\cF$ with supports in $C'$, then it is  true for $\cF$ with supports in $C\cup C'$. 
Note that $X$ is a Noetherian space. Hence by Noetherian induction, we may assume that $X$ is irreducible and \ref{fin2} (1)  is true if $\cF$ has supports in  some closed subset $C\neq X$ of $X$. Hence we assume that $X$ is a prime non-commutative scheme. Let $\eta$ be the generic point of $X$ and let $Y'$ be the closure of $f(\eta)$ in $Y$ with the prime non-commutative scheme structure. We can replace $Y$ by $Y'$ and hence we assume that $Y$ is also a prime non-commutative scheme.

 For open immersions $U\to X$ and $U'\to X$ such that $X$ is the union of $U$ and $U'$, we have an exact sequence (let $U''=U\cap U'$)

(2)  $\dots  \to R^mf_*\cF\to R^mf_{U, *}(\cF_U)\oplus R^mf_{U', *}(\cF_{U'})\to R^mf_{U'',*}(\cF_{U''}) \to R^{m+1}f_*\cF\to \dots,$

\noindent
where $f_U$ denoted the composition $U\to X \overset{f}\to Y$ and $f_{U'}$, $f_{U''}$ are defined similarly, $\cF_U$ denotes the restriction of $\cF$ to $U$, and $\cF_{U'}$, $\cF_{U''}$ are defined similarly.  By this, we are reduced to the situation of \ref{spcase1}. (The regularity of $Z(A)$ in \ref{spcase1} is attained since the set of regular points in an excellent scheme is open.)

Let $C=X\smallsetminus U$ and let  $j:U\to X$ be the inclusion morphism. 
We have a distinguished triangle

(3) $\cF \to Rj_*(\cF_U)\oplus \cF_C \to (Rj_*(\cF_U))_C \to.$ 

Note that $Rj_*(\cF_U)$ is   constructible by the case $X=Y$ of \ref{spcase2.5}, and hence  $(Rj_*(\cF_U))_C$ is constructible. 
By Noetherian induction,  $Rf_*(\cF_C)$ and $Rf_*((Rj_*(\cF_U))_C)$ are constructible. On the other hand, $Rf_*Rj_*(\cF_U)=Rg_*(\cF_U)$ is constructible by \ref{spcase2.5}. Hence $Rf_*\cF$ is constructible.  

The proof of \ref{fin2} (2) is similar to the above proof of \ref{fin2} (1). We just replace the above exact sequences (1) and (2) and the distinguished triangle (3), by the following facts (4), (5), and (6), respectively.

(4) The sheaf $\cF_{C\cup C'}$ is the fiber product of $\cF_C \to \cF_{C\cap C'} \leftarrow \cF_{C'}$. 

(5) The sheaf $\cF$ is the fiber product of $f_{U,*}(\cF_U)\to f_{U'',*}(\cF_{U''}) \leftarrow f_{U',*}(\cF_{U'})$. 

(6)  The sheaf $\cF$ is the fiber product of $\cF_C\to (j_*(\cF_U))_C  \leftarrow j_*(\cF_U)$.

We prove (3) of \ref{fin2}. 
Because we have  proved (2) of \ref{fin2}, we have

{\bf Claim}.  Let $\cF$ and $\cF'$ be constructible sheaves of groups on $X_{et}$ and let  $\cF\to \cF'$ be an injective homomorphism. If $R^1f_*(\cF')$ is constructible, then $R^1f_*\cF$ is constructible.

In the case of schemes, this is Lemma 3.1.1 of Exp. XXI of \cite{ILO}. The same proof works. 

For a quasi-compact non-commutative scheme $X$ satisfying ({\bf F}), let  $S_X$ (resp. $S_{X,C}$ for a closed subset $C$ of $X$) be the statement that \ref{fin2} (3) is true for  every $Y$ and every $f:X\to Y$ and every $\cF$ (resp. every $\cF$ with supports in $C$)  as in the hypothesis of \ref{fin2}. 
For closed subsets $C_i$ ($1\leq i\leq n$) of $X$, if $S_{X,C_i}$ is true for $1\leq i\leq n$, then  $S_{X,C}$ for $C=\cup_{i=1}^n C_i$ is true. This is seen by the above Claim  applied to the injective homomorphism $\cF_C\to \prod_{i=1}^n \cF_{C_i}$.
By this, for the proof of $S_X$, we may assume that $X$ is a prime non-commutative scheme and that $S_{X,C}$ is true for every closed subset $C\neq X$ of $X$. For open subsets $U_i$ ($1\leq i\leq n$) of $X$ such that $X=\cup_{i=1}^n U_i$, if $S_{U_i}$ are true for $1\leq i\leq n$, then $S_X$ is true. In fact, for each $i$,  if we denote by $j_i$ the inclusion morphism $U_i\to X$, $R^1(f\circ j_i)_*(\cF_{U_i})$ and $R^1j_{i,*}(\cF_{U_i})$ are constructible. Hence $f_*R^1j_{i,*}(\cF_{U_i})$ is constructible by (2) of \ref{fin2} and hence 
the kernel $R^1f_*j_{i,*}(\cF_{U_i})$ of $R^1(f\circ j_i)_*(\cF_{U_i})\to f_*R^1j_{i,*}(\cF_{U_i})$  is constructible. By Claim  applied to the injective homomorphism $\cF \to \prod_{i=1}^n j_{i,*}(\cF_{U_i})$, $R^1f_*\cF$ is constructible. By this, we are reduced to the situation of \ref{spcase1}. 
Now we consider  the injective homomorphism $\cF\to  j_*(\cF_U)\times \cF_C$. The sheaf $R^1f_*j_*(\cF_U)$ is constructible because it is the kernel of
$R^1g_*(\cF_U) \to f_*R^1j_*(\cF_U)$ and $R^1g_*(\cF_U)$ (resp. $ f_*R^1j_*(\cF_U)$) is  constructible by \ref{spcase2} (resp. by the case $Y=X$ of  \ref{spcase2.5} and by (2) of \ref{fin2}). The sheaf $R^1f_*(\cF_C)$ is constructible by the Noetherian induction. Hence $R^1f_*\cF$ is constructible by Claim. This completes the proof of \ref{fin2}.

\end{para}

The following \ref{pbc} is a non-commutative version of  the proper base change theorem in the usual \'etale cohomology theory. This will be further generalized in \ref{pbc2}. 

\begin{thm}\label{pbc} Let $Y$ be a non-commutative scheme satisfying ({\bf F}). 
 Let $T\to S$ be a proper morphism of schemes, let $Y\to S$ be a morphism of non-commutative schemes, and let $X=Y\times_S T$. Then we have the following proper base change theorem for $f: X\to Y$. Let $Y'$ be a non-commutative scheme satisfying ({\bf F}) and let $h:Y'\to Y$ be a morphism, and consider the following commutative diagram with $X'=X\times_Y Y'$ (note that $f$ is relatively commutative by our assumption and hence the fiber product is defined). $$\begin{matrix}  X' &\overset{g}\to& X\\ {f'}\downarrow\;\;\;\; &&\;\;\;\downarrow{f}\\
Y'&\overset{h}\to& Y\end{matrix}$$
Then for any sheaf of torsion abelian groups $\cF$ on $X_{et}$, the canonical morphism $h^{-1}Rf_*\cF\to Rf'_*g^{-1}\cF$ is an isomorphism.

\end{thm}

\begin{pf}   
We may assume that $S=\Spec(R)$ for an excellent commutative ring $R$ and  $X=\cU(x)$ for some $R$-algebra $A$ 
which is finitely generated as an $R$-module and for some closed point $x$ of  $\Spec(A)$. 
It is sufficient to prove that $(h^{-1}Rf_*\cF)_{\bar y'}\to (Rf'_*g^{-1}\cF)_{\bar y'}$ is an isomorphism 
for every $y'\in Y'$.  Let $y$ be the image of $y'$ in $Y$. Let $X(y)=X\times_Y y$, $X'(\bar y')= X'\times_{Y'} \bar y'$ and consider the cartesian  squares 
$$\begin{matrix}   X(y) & \to & X & X'(\bar y')  &\to& X' &  X'(\bar y') &\to & X(y) \\
\downarrow && \downarrow   & \downarrow && \downarrow &\downarrow &&\downarrow\\
y& \to& Y & \bar y' & \to& Y' & \bar y' & \to & y.  \end{matrix}$$
It is sufficient to prove \ref{pbc} in the cases of these squares (that is, in the cases the cartesian square in \ref{pbc} is replaced by these cartesian squares).
In the second square, by localizing $Y'$, we may assume that $y'=\bar y'$. 
It is sufficient to consider the first and the third squares. Note that $X(y)= T(s)\times_s y$,  $X'(y')= T(s) \times_s y'$, where $T(s)=T\times_S s$.

Now we  consider the first square. Let  the notation be $\epsilon^X: X\to T$, $\epsilon^Y: Y\to S$,  $\epsilon^{X(y)}: X(y)  \to T(s)$, $\epsilon^y: y\to s$, $f^y: X(y) \to y$,  $p: T\to S$,  $p^s: T(s)\to s$. 
Now \ref{pbc} for the first square is proved by 
$$(Rf_*\cF)_x= (\epsilon^Y_*Rf_*\cF)_s= (Rp_*\epsilon^X_*\cF)_s= (Rp^s_*((\epsilon^X_*\cF)|_{T(s)}))_s= (Rp^s_*\epsilon^{X(y)}_*(\cF|_{X(y)}))_s$$ 
$$
=(\epsilon^y_*Rf^y_*(\cF|_{X(y)}))_s= (Rf^y_*(\cF|_{X(y)}))_x.$$
Here, the second and the fifth $=$ are the evident ones. The first and the last (sixth) $=$ are by \ref{small4} (2).  In the third $=$, we use the usual proper base change theorem of schemes for the proper morphism $p:T\to S$. The fourth $=$ is by 

{\bf Claim.} $(\epsilon^X_*\cF)|_{T(s)}= \epsilon^{X(y)}_*(\cF|_{X(y)})$.

We prove Claim. 
Let  $t\in T(s)$ and let $x$ the unique  point of $X(y)$ lying over $t$. Then if we denote the the strict localization of $T$ at $t$ by $T_{\bar t}$, the morphism $X\times_T T_{\bar t}\to T_{\bar t}$ satisfies the present condition of $X\to S$ and $X\times_T T_{\bar t}=\cU(x)$. Hence we have 
$(\epsilon^X_*\cF)|_{T(s)})_{\bar t}= \cF_{\bar x}$ and 
$(\epsilon^{X(y)}_*(\cF|_{X(y)}))_{\bar t}= (\cF|_{X(y)})_{\bar x}=\cF_{\bar x}$ by \ref{small4} (2). 
This proves Claim and hence \ref{pbc} for the first square.

Now we consider the third square. Note that $s=\Spec(k)$, $y=\Spec(M_n(k))$, $y'=\Spec(M_{mn}(k'))$ for some separably closed commutative fields $k$ and $k'$ and for some $m,n\geq 1$. Since $X(y)=T(s) \times_s y$, a sheaf  of torsion abelian groups on $X(y)_{et}$ comes form a sheaf  of torsion abelian groups on $T(s)_{et}$ by \ref{Azet}. 
Let $s'=\Spec(k')$. Then $X'(y')= (T(s) \times_s s')\times_{s'} y'$. Hence by \ref{Azet}, \ref{pbc} for the third square follows from the usual proper base change theorem for the proper morphism $T(s)\to s$ and the base change by  $s'\to s$. 
\end{pf}

\begin{para}\label{indep} Let $X$ and $Y$ be non-commutative schemes satisfying ({\bf F}) and let $f:X\to Y$ be a morphism satisfying the following condition \ref{indep}.1.

\medskip

\ref{indep}.1.  Locally on $Y$, there are schemes $S$ and $T$, a proper morphism $T\to S$,  a morphism $Y\to S$, and an l.f.p.  immersion $j: X\to Y\times_S T$ such that $f$ is the composition $X\overset{j}\to Y\times_ST \to Y$. 

\medskip

We define $Rf_!\cF$ for a sheaf $\cF$ of torsion abelian groups on $X_{et}$. Assume we have $T\to S$, $Y\to S$, and $X\overset{j}\to Y\times_S T$ as above. We define 
 $$Rf_!:=Rp_*\circ  j_!$$ where  $p:Y\times_S T \to Y$. This is independent  of the choices. In fact, for other $T'\to S'$, $Y\to S'$, and $X\overset{j'}\to Y\times_{S'} T'$, we have the third choices 
 $T''=T\times T'\to S''=S\times S'$ (here $\times$ is the fiber product over $\Spec(\Z)$), $Y\to S\times S'$, and  $X\overset{j''}\to Y\times_{S''} T''$, and we have 
 a commutative diagram
 $$\begin{matrix} X &\overset{j''}\to & Y\times_{S''} T''&\overset{\pi}\to & Y\times_{S''} (T \times S')  = Y \times_S T & \to & Y \\
 \Vert&&&&\Vert&&\Vert\\
 X&& \overset{j}\to && Y\times_S T & \to & Y\end{matrix}$$
in which $T'' \to T\times S'$ is proper. We show that $j''$ is an l.f.p. immersion and that $j_!= R\pi_*\circ  j''_!$ for sheaves of torsion abelian groups. Then we will have $Rp_*\circ j_!= Rp_*\circ R\pi_* \circ j''_!= Rp''_*\circ  j''_!$ where $p'': Y \times_{S''} T''\to Y$, and we have similarly $Rp'_*\circ j'_!= Rp''_*\circ j''_!$ and hence consequently $Rp_*\circ j_!=Rp'_*\circ j'_!$. 
 
 This $j_!=R\pi_*\circ j''_!$ is  reduced to the following. Assume we are given non-commutative schemes $X$ and $Y$ satisfying ({\bf F}), an l.f.p. immersion $j:X\to Y$, a proper morphism of schemes $T\to S$, a morphism $Y\to S$, and a morphism $j': X\to Y \times_S T$ over $Y$. Then $j'$ is an l.f.p. immersion and $j_!=Rp_* \circ j'_!$ for sheaves of torsion abelian groups, where $p:Y\times_S T\to Y$. 
 The fact $j'$ is an l.f.p. immersion follows from \ref{sep2} (3). We have a canonical morphism $j_!F\to Rp_*\circ j'_!F$. We prove $(j_!F)_{\bar y}\to (Rp_*\circ j'_!F)_{\bar y}$ is an isomorphism for every $y\in Y$. By \ref{pbc}, this is reduced to the case $Y= \bar y$. In this case,  $X$ is either $\bar y$ or $\emptyset$. In the former case, both $j_!\cF$ and $Rp_*\circ j_!\cF$ are $\cF$. The latter case is clear. 
 
 Thus $Rf_!\cF$ is well defined for sheaves of torsion abelian groups on $X_{et}$. 

 \medskip

If $X$, $Y$, $Z$ are non-commutative schemes satisfying ({\bf F}) and if $f:X\to Y$ and $g:Y\to Z$ are morphisms satisfying the  following condition \ref{indep}.2 for a fixed scheme $S$, we have
$$R(g\circ f)_!=Rg_!\circ Rf_!$$
for sheaves of torsion abelian groups. 

\medskip
\ref{indep}.2. There are  proper schemes $T$ and $T'$ over $S$ and  l.f.p.  immersions $j: X\to Y\times_S T$ and $ j':Y\to Z\times_S T'$ such that $f$ is the composition $X\overset{j}\to Y\times_S T \to Y$ and $g$ is the composition $Y\overset{j'}\to Z\times_S T'\to Z$. 

\medskip

Note that then $g\circ f$ is the composition  $X\to Z\times_S (T\times_S T')\to Z$ in which the first arrow is an l.f.p. immersion and the $T\times_ST'\to S$ is proper, and hence $R(g\circ f)_!$ is defined.

\end{para}

\begin{thm}\label{pbc2} Let $X$ and $Y$ be non-commutative schemes satisfying ({\bf F}), and let  $f:X\to Y$ be a morphism  satisfying the condition \ref{indep}.1.

(1) If $Y'$ is a non-commutative scheme satisfying ({\bf F}) and if $Y'\to Y$ is a morphism, $R^mf_!$ and $R^mf'_!$ with $f':X\times_Y Y'\to Y'$ commute with pullbacks of sheaves of torsion  abelian groups.

(2) $R^mf_!$ sends constructible sheaves of abelian groups  on $X_{et}$ to constructible sheaves of abelian groups on $Y_{et}$. If $m\gg 0$, $R^mf_!\cF=0$ for all constructible sheaves $\cF$ of abelian groups on $X_{et}$. 

Remark. For a constructible sheaf $\cF$ of abelian groups on $X_{et}$ such that the orders of all stalks of $\cF$ are invertible on $Y$, the constructibility of $R^mf_!\cF$  in (2) follows from \ref{fin2}.
\end{thm}

\begin{pf}

(1) follows from  \ref{pbc}. ((1) for the case $f$ is an l.f.p. closed immersion and the case $f$ is an open immersion are clear by checking stalks.)

(2) By Noetherian induction, we may assume that $Y$ is a prime non-commutative scheme and $(Rf_!\cF)_C$ is constructible for every closed subset $C\neq Y$ of $Y$ and every constructible sheaf of abelian groups $\cF$ on $X_{et}$. 
In fact, if $Rf_!\cF$ is not constructible for some  constructible $\cF$, there is a minimal element in the set of all closed subsets $C$ of $Y$ satisfying the following condition:  $(Rf_!\cF)_{C'}$ is not constructible for some closed subset $C'$ of $C$ and for some constructible sheaf $\cF$ of abelian groups on $X_{et}$ . If $C=C_1\cup C_2$ for closed subsets $C_i$, we would have $C=C_1$ or $C=C_2$ because otherwise, for every closed subset $C'$ of $C$, the exact sequence 
 (1) in the proof of \ref{fin2},  which we apply by replacing $X$, $\cF$, $C$ and $C'$ there with $Y$, $R^mf_!(\cF)$, $C_1\cap C'$ and $C_2\cap C'$ here,   would show that $(Rf_!\cF)_{C'}$ is constructible.  Hence this minimal $C$ is irreducible. 
We have $(Rf_!\cF)_C= Rf'_!(\cF_{X'})$ for this minimal $C$ which we endow with the prime non-commutative scheme structure, where $f': X'=X\times_Y C\to Y$. We may assume that $Y$ is this minimal $C$. 

By Noetherian induction by using the exact sequence (1) in the proof of \ref{fin}, we may assume that  $X$ is a prime non-commutative scheme and $Rf_!(\cF_C)$ is constructible for  every  closed subset $C\neq X$ of $X$. 
 We may assume that $Y$ is the closure of the image of the generic point of $X$ in $Y$ and that $Y$ has the prime non-commutative scheme structure. 
 If $X=U\cup U'$ for open subsets of $X$, we have an exact sequence $0\to j''_!(\cF_{U''}) \to j_!(\cF_U) \oplus j_!\cF_{U'} \to \cF\to 0$, where $U''=U\cap U'$ and $j$, $j'$, $j''$ are inclusion morphisms from $U$, $U'$, $U''$ to $X$, respectively.
  By this, we may assume that $Y$ is an open subspace of $\Spec(A)$ for some prime ring $A$ such that $Z(A)$ is an excellent ring and $A$ is finitely generated as a $Z(A)$-module, 
$X$ is an open subspace of $\Spec(B)$ for some prime ring $B$ such that $Z(B)$ is an excellent ring and $B$ is finitely generated as a $Z(B)$-module, 
    $B$ is an $A$-algebra of finite presentation, $A\to B$ is injective, and the morphism $\Spec(B)\to \Spec(A)$ is compatible 
    with $f:X\to Y$.  There are a non-zero element $a$ of $Z(A)$ and a non-zero element $b$ of $Z(B)$ such that $b/a\in Z(B)$ and such that $A[1/a]$ is 
    an Azumaya algebra over $Z(A)[1/a]$, $V:= Y\times_{\Spec(Z(A))} \Spec(Z(A)[1/a])$ coincides with $\Spec(A[1/a])$, $B[1/b]$ is an Azumaya algebra over $Z(B)[1/b]$, and $U:= X \times_{\Spec(Z(B))} \Spec(Z(B)[1/b])$ coincides with $ \Spec(B[1/b])$. 
    Consider the exact sequence $0\to j_!(\cF_U) \to \cF \to \cF_C\to 0$ on $X_{et}$,  where $C=X\smallsetminus U$ and $j$ is the inclusion morphism $U\to X$. By our Noetherian induction on $X$, $Rf_!(\cF_C)$ is constructible. Let $D=Y\smallsetminus V$. By our Noetherian induction on $Y$, $(Rf_!j_!(\cF_U))_D$ is constructible. It remains to prove that the pullback of 
    $Rf_!j_!(\cF_U)$ on $V_{et}$ is constructible, that is, $Rg_!(\cF_U)$ is constructible, where $g:U\to V$. Let $h: \Spec(Z(B)[1/b])\to \Spec(Z(A)[1/a])$ and 
    let $\cG$ be the sheaf on $(\Spec(Z(B)[1/b])_{et}$ corresponding to $\cF_U$ via \ref{Azet}. Then 
    $Rg_!(\cF_U)$ corresponds to $Rh_!\cG$ via \ref{Azet} and the constructibility of the former is reduced to the constructibility of the latter in the scheme theory (\cite{SGA4}, vol. 3, Exp. XIV, Theorem 1.1).
      \end{pf}

\begin{para}\label{ladic}  Let $X$ be a non-commutative scheme satisfying ({\bf F}) and let $\ell$ be a prime number which is invertible on $X$. We consider constructible $\La$-sheaves on $X$ for $\La=O_E, E, \bar \Q_{\ell}$ for a finite extension $E$ of $\Q_{\ell}$ ($O_E$ denotes the integer ring of $E$) and for an algebraic closure $\bar \Q_{\ell}$ of $\Q_{\ell}$,  following  \cite{SGA4}, \cite{De}, \cite{SGA5}.

An $O_E$-sheaf on $X$  is a family ($\cF_n)_{n\geq 1}$ of sheaves $\cF_n$  of $O_E/\ell^nO_E$ -modules on $X_{et}$ endowed with an isomorphism $\cF_{n+1}\otimes_{O_E/\ell^{n+1}O_E} O_E/\ell^nO_E\cong \cF_n$ for each $n$. We say a  $\La$-sheaf is smooth if $\cF_n$ are locally constant and finite. We say a $\La$-sheaf is constructible if locally on $X$, there is a finite family of l.f.p. immersions $Y_i\to X$ such that the set $X$ is the disjoint union of these $Y_i$ and such that the pullback of $\cF$ to each $Y_i$ is smooth. 

The category of constructible $E$-sheaves on $X$ is defined as the quotient category of the category of constructible $O_E$-sheaves by the full subcategory consisting of objects which are killed by some powers of $\ell^n$. 

The category of constructible $\bar \Q_{\ell}$-sheaves on $X$  is defined to be the inductive limit of categories of constructible $E$-sheaves, where $E$ ranges over all finite extensions of $\Q_{\ell}$ in $\bar \Q_{\ell}$.

In the case $X$ satisfies (F$_\Z$), we consider mixed sheaves following \cite{De}. Fix an isomorphism $\bar \Q_{\ell}\cong \C$ of commutative fields. A mixed $\bar \Q_{\ell}$-sheaf on $X$ is a constructible $\bar \Q_{\ell}$-sheaf $\cF$ satisfying the following condition. There is a finite family $\cF^{(i)}$ ($0\leq i\leq m+1$) of constructible $\bar \Q_{\ell}$-subsheaves $\cF^{(i)}$ of $\cF$  such that $\cF^{(0)}\supset \cF^{(1)}\supset \dots \supset \cF^{(m+1)}=0$ and such that for $0\leq i\leq m$, there is an integer $w(i)$ satisfying the following condition (*). For each closed point $x$ of $X$, since $\cO_{X,x}/\frak p(x)\cong M_r(\F_q)$ for some $r\geq 1$ and for some finite field $\F_q$, the stalk $(\cF^{(i)}/\cF^{(i+1)})_{\bar x}$, which is a finite-dimensional $\bar \Q_{\ell}$-vector space,  has a linear action of $\Gal(\bar \F_q/\F_q)$ and in particular, the action of the Frobenius $\varphi_x\in\Gal(\bar \F_q/\F_q)$, where $\varphi_x(a)=a^q$ for $a\in \bar \F_q$. 

(*) For each $x\in X$, all the eigen values of $\varphi_x^{-1}$ on $(\cF^{(i)}/\cF^{(i+1)})_{\bar x}$ regarded as  elements of $\C^\times$ have absolute value $q^{w(i)/2}$.

\end{para}

\begin{thm}\label{lfin} Let $X$, $Y$ and $f:X\to Y$ be as in \ref{fin2}, let $\ell$ be a prime number which is invertible on $Y$, and let $\cF$ be a constructible $\bar \Q_{\ell}$-sheaf on $X_{et}$. Then $R^mf_*\cF$ is a constructible $\bar \Q_{\ell}$-sheaf on $Y$ for every $m$ and is $0$ for $m\gg 0$. 

\end{thm}

\begin{pf} This follows from \ref{fin2} as in the case of schemes treated in \cite{SGA5}.  \end{pf}

\begin{thm}\label{pbc3} Let $X$, $Y$ and $f: X\to Y$ be as in \ref{pbc2}. Let $\ell$ be a prime number which is invertible on $Y$. 

(1) $R^mf_!$ sends constructible $\bar \Q_{\ell}$ sheaves on $X$ to constructible $\bar \Q_{\ell}$-sheaves on $Y$. If $m\gg 0$, $R^mf_!\cF=0$ for all constructible $\bar \Q_{\ell}$-sheaves $\cF$ on $X$. 

(2) For a non-commutative scheme $Y'$ satisfying ({\bf F}) and for a morphism $Y'\to Y$, $R^mf_!$ and $R^mf'_!$ with $f':X\times_Y Y'\to Y'$ commute with pullbacks of constructible $\bar \Q_{\ell}$-sheaves.

(3) If $Y$ satisfies (F$_\Z$), $R^mf_!$ sends mixed $\bar \Q_{\ell}$ sheaves on $X$ to mixed $\bar \Q_{\ell}$ sheaves on $Y$.

\end{thm}

\begin{pf}
(1) and (2) are deduced from (2) and (1) of \ref{pbc2}, respectively, as in the case of schemes treated in \cite{SGA5}. (Actually, (1) follows from \ref{lfin}.) 
  
(3) By the method of the proof of (2) of \ref{pbc2}, we are reduced to the case $X$ and $Y$ are schemes, and to the theory of mixed sheaves of Deligne in  \cite{De}.
\end{pf}

\section{Betti cohomology}\label{sec6}

 For a scheme $S$ over $\C$ locally of finite type, the set $S(\C)$ of $\C$-valued points has the topology called the classical topology. 
Recall that as a set,  $S(\C)$ is identified with the set of all closed points of $S$.

In this Section \ref{sec6}, we consider a non-commutative version of this classical topology. For a non-commutative scheme $X$ over $\C$ satisfying the finiteness condition ({\bf F}$_\C$) (\ref{bfF}), we endow the set $X_{\cl}$ of all closed points of $X$ with a topology which we call the classical topology. (The notation $\cl$ in $X_{\cl}$ stands for classical.) We will call the cohomology of $X_{\cl}$ the Betti cohomology of $X$. 
 We will prove a comparison theorem (\ref{e=B}) between the \'etale cohomology of $X$ and the Betti cohomology of $X$.

\begin{lem} Assume $X$ satisfies ({\bf F}$_\C$), and let $x\in X$. Then $x\in X_{cl}$ if and only if the residue field of the local ring $Z(\cO_{X,x})$ is $\C$.

\end{lem}

\begin{pf}
This follows from \ref{fin} (1). 
\end{pf}

For $X$ satisfying ({\bf F}$_{\C}$), the classical topology on $X_{\cl}$ is defined in \ref{cl1} (1). The style of this definition follows the following  \ref{old1} and \ref{old2} concerning the classical topologies of $S(\C)$  for schemes $S$ locally of finite type over $\C$. 

\begin{para}\label{old1}

Recall that there is a unique way to endow the set $S(\C)$ for every scheme $S$ over $\C$ locally of finite type with a topology (the classical topology)  satisfying the following (i) and (ii).

(i) If $U$ is an open set of $S$, $U(\C)$ is open in $S(\C)$ and the classical topology of $U(\C)$ is the restriction of the classical topology of $S(\C)$.

(ii) If $S$ is $\Spec(R)$ for a 
 finitely generated commutative ring  $R$ over $\C$, the classical topology of $S(\C)$ is the topology of simple convergence of $\Hom_\C(R, \C)=S(\C)$, and it is also understood as in \ref{old2} below.

\end{para}

\begin{lem}\label{old2}  Let $R$ be a commutative finitely generated algebra over $\C$, let $S=\Spec(R)$, and let $R^{\ct}$ be the ring of all $\C$-valued continuous functions on $S(\C)$ for the classical topology. Let $\iota: S(\C) \to \max(R^{\ct})$ be the map which sends $s\in S(\C)$ to the kernel of the evaluation $R^{\ct} \to \C$ at $s$. Then $\iota$ is injective and the classical topology of $S(\C)$ coincides with the pullback of the Zariski topology (Jacobson topology) of $\max(R^{\ct})$ via $\iota$. 

\end{lem}

\begin{pf} This is well known. 
We are reduced to the case $R=\C[T_1, \dots, T_n]$. In this case, for $\alpha= (\alpha_1, \dots, \alpha_n)\in S(\C)=\C^n$, for $h_r\in R^{\ct}$ ($r>0$) defined by $h_r(z):= \prod_{i=1}^n \max(r- |z_i-\alpha_i|, 0)$ ($z\in \C^n$),  the sets $\{z \in\C^n\;|\; h_r(z)\neq 0\}=\{z\in \C^n\;|\; |z_i-\alpha_i|<1 \;(1\leq i\leq n)\}$ for $r>0$ form a base of neighborhoods of $\alpha$. 
\end{pf}

\begin{prop}\label{cl1} (1) There is a unique way to endow $X_{\cl}$ with a topology for every non-commutative scheme $X$ over $\C$ satisfying the condition ({\bf F}$_{\C}$), satisfying the following conditions (i) and (ii).

(i) If $U$ is an open subspace of $X$, $U_{\cl}$ is open in $X_{\cl}$ and the topology of $U_{\cl}$ is the restriction of the topology of $X_{\cl}$. 

(ii) Let $R$ be a finitely generated commutative ring over $\C$, let $A$ be an $R$-algebra which is finitely generated as an $R$-module, and let $X=\Spec(A)$. 
Let $S=\Spec(R)$,  let $R^{\ct}$ be the  ring of all $\C$-valued continuous functions on $S(\C)$, and consider the injection  $\iota: 
\max(A) \to \max(R^{\ct}\otimes_R A)$ which sends $\m\in \max(A)$ to the kernel of 
  $R^{\ct} \otimes_R A\to \C \otimes_R A= A/\m'A\to A/\m$, where $\m'\in \max(R)$ is the inverse image of $\m$ under $R\to A$, the first arrow is defined by the evaluation  $R^{\ct}\to \C$ at $\m'$,  and the map $R\to \C$ in $\C\otimes_R$ is the evaluation at  $\m'$. 
Then the topology of $X_{\cl}$ coincides with the pullback of the Zariski topology (Jacobson topology) of $\max(R^{\ct}\otimes_R A)$ via $\iota$. 

We will call this topology of $X_{\cl}$ the classical topology.
\medskip

(2) The classical  topology of $X_{\cl}$ is the weakest topology such that the following sets $D_{U,A}(f)$ are open. Here $A$ is a $\C$-algebra such that $Z(A)$ is a finitely generated $\C$-algebra and $A$ is finitely generated as a $Z(A)$-module, $U$ is an open set of $X$ endowed with a morphism $U\to \Spec(A)$ over $\C$, $f$ is an element of $Z(A)^{\ct} \otimes_{Z(A)} A$ ($Z(A)^{\ct}$ is as in (1) (ii)), and $D_{U,A}(f)\subset U_{\cl}$ denotes the inverse image of
$D(f):=\{\m\in \max(Z(A)^{\ct} \otimes_{Z(A)} A)\;|\;f\notin \m\}$ under $U_{\cl}\to \max(A)\overset{\iota}\to \max(Z(A)^{\ct} \otimes_{Z(A)} A)$. 

\medskip

The following (3) (resp. (4)) is similar to (2) (resp. (1) (ii)) but uses analytic functions instead of continuous functions. 

\medskip

(3) The classical  topology of $X_{\cl}$ is the weakest topology such that the following sets $D_{U,A,V}(f)$ are open. Here $A$ and $U$ are as in (2),  $V$ is an open set of $S(\C)$, where $S:=\Spec(Z(A))$, $f$ is an element of $A_V:=\cO^{\an}_S(V) \otimes_{Z(A)} A$, 
where $\cO_S^{\an}$ is the structure sheaf of the complex analytic space associated to $S$,  and $D_{U,A,V}(f)\subset U_{\cl}\times_{S(\C)} V$ denotes the inverse image of $D(f):= \{\m\in \max(A_V)\;|\; f\notin \m\}$ under the map $U_{\cl}\times_{S(\C)} V \to \max(A_V)$ which sends $(u, v) \in U_{\cl}\times_{S(\C)} V$ with image $s$ in $S(\C)$ to the kernel of $A_V\to \C \otimes_{\cO(S)} \cO(U)\to \cO_{U, u}/\frak p(u)$. Here the first arrow is induced by the evaluation $\cO_S^{\an}(V)\to \C$ at $v$ and $\cO(S)\to \C$ is the evaluation at $s$. 

(4) Let $R$ be a finitely generated commutative ring over $\C$, let $A$ be an $R$-algebra which is finitely generated as an $R$-module, and  assume that $X$ is an open subspace of $\Spec(A)$. Then  the classical topology of $X_{\cl}$ is the weakest topology such that the following sets $D_V(f)$ are open. 
Here $V$ is an open set of $S(\C)$, where $S:=\Spec(R)$, $f$ is an element of $A_V:=\cO^{\an}_S(V) \otimes_R A$, 
and $D_V(f)\subset U_{\cl}\times_{S(\C)} V$ denotes the inverse image of $D(f)\subset  \max(A_V)$.

(5) Let $x\in X_{\cl}$.  Take an open neighborhood $U$ of $x$, a finitely generated commutative ring $R$ over $\C$, an $R$-algebra $A$ which is finitely generated as an $R$-module, and an open immersion $U\to \Spec(A)$ over $\C$. 
Then a base of neighborhoods of $x$ in $X_{\cl}$ (for the classical topology) is given as follows.
Denote $\Spec(R)$ by $S$, the canonical map $U_{\cl}\to S(\C)$ by $\pi$, and let $s=\pi(x)$. 
   Take an open neighborhood $W$ of $s$ in $S(\C)$ and  an element $e\in \cO^{an}_S(W)\otimes_R  A$ such that the image of $e$ in $\cO^{\an}_{S,s} \otimes_R A$ is an idempotent for $(\cO^{\an}_{S,s} \otimes_R A, \m)$, where $\m$ is the maximal two-sided ideal corresponding to $x$ (\ref{idem}, such $e$ exists by \ref{idem2} (1)).
 Let $D(e)$ be the open subset of $\pi^{-1}(W)\subset X_{\cl}$ consisting of all points $y$ such that the image of $e$ in $\cO_{X,y}/\frak p(y)$ is not zero.
Then when $W'$ ranges over all open neighborhoods of $s$ in $W$,   $D(e)\cap \pi^{-1}(W')$ form a base of neighborhoods of $x$ in $X_{\cl}$. 
\end{prop}

\begin{pf} Let  $T_2$ (resp. $T_3$)  be the topology of $X_{\cl}$ defined as in the above (2) (resp. (3))  in \ref{cl1}. Assuming we are given $R$ and $A$ and an open immersion $X\to \Spec(A)$ as in (ii) of (1)  (resp. as in (4)), let $T_1$ (resp. $T_4$) be the topology of $X_{\cl}$ defined as in (ii) of (1) (resp. as in (4)). Let $x\in X$. For the proof of \ref{cl1}, it is sufficient to prove that for each $k=1,2,3,4$ and for the topology $T_k$, $D(e)\cap \pi^{-1}(W')$ in (5) form a base of neighborhoods of $x$ in $X_{\cl}$. 

We first prove that $D(e)\cap \pi^{-1}(W')$ is a neighborhood of $x$ for $T_k$. The cases $k=3, 4$ are clear. The case $k=2$ follows from the case $k=1$. In the case $k=1$, take an open set $V$ and a closed set $C$ of $S(\C)$ 
 for the classical topology such that $s\in V\subset C\subset W$.   Write $e=\sum_{i=1}^n  b_i \otimes a_i$ with $b_i\in \cO_S^{an}(W)$ and $a_i\in A$, and let $c_i$ be a $\C$-valued continuous function on $S(\C)$ 
 which coincides with $b_i$ on $V$ and with $0$ on $S(\C)\smallsetminus C$. Let $e'= \sum_{i=1}^n c_i\otimes a_i\in R^{\ct} \otimes_R A$. Then $e$ and $e'$ coincide on $\pi^{-1}(V)$. This shows that for an open neighborhood $W'$ of $s$ in $V$, $D(e) \cap \pi^{-1}(W')$ is a neighborhood of $x$ for $T_1$.

We next prove that each neighborhood of $x$ in $X_{\cl}$ for $T_k$ contains $D(e)\cap\pi^{-1}(W')$ for some $W'$. For a $\C$-algebra $A$ such that $Z(A)$ is a finitely generated $\C$-algebra and $A$ is finitely generated as a $Z(A)$-module and for a morphism $X\to \Spec(A)$ over $\C$, locally on $X$, we have a  $\C$-algebra $B$ such that $Z(B)$ is a finitely generated $\C$-algebra and $B$ is finitely generated a $Z(B)$-module and 
an open immersion $X\subset \Spec(B)$ over $\C$ with $B$, and a morphism $A\to B$ in (Alg$/\C$) such that $X\to \Spec(A)$ factors as $X\to \Spec(B)\to \Spec(A)$. From this, the case $T_2$ is reduced to the case $T_1$, and the case $T_3$ is reduced to the case $T_4$. We consider $T_1$ (resp. $T_4$).
Let $f$ be an element as in the explanation of $T_1$ (resp. $T_4$). It is sufficient to prove that if  $x\notin D_{X,A, S}(f)$ (resp. $D_V(f)$), $D_{X, A, S}(f)$ (resp. $D_V(f)$) contains $D(e)\cap \pi^{-1}(W')$ for some open neighborhood $W'$ of $s$ in $W$. Let $\cO_S^{\ct}$ be the sheaf of $\C$-valued continuous functions on $S(\C)$. For $(*)=\an, \ct$, let $f_x$ be the image of $f$ in $\cO^{(*)}_{S,s} \otimes_R A$. By \ref{idem2} (2)  applied to $I=(f_x)$,  $e\in (f_x)$ in this ring. Hence for some open neighborhood 
$W'$ of $s$ in $W$ (resp. $W\cap V$), we have $e\in (f)$ in $\cO_S^{(*)}(W) \otimes_R A$. Hence  $D(e)\cap \pi^{-1}(W')\subset D_{X,A,S}(f)$ (resp. $D_V(f)$). 
\end{pf}

\begin{lem} For a morphism $Y\to X$ of non-commutative schemes over $\C$ satisfying ({\bf F}$_\C$), the induced map  $Y_{cl} \to X_{cl}$ is continuous. \end{lem}
\begin{pf} This is seen from (2) or (3) in \ref{cl1}. 
\end{pf}

\begin{lem}\label{AzB} Let $S$ be a scheme over $\C$ locally of finite type, and let $X=\Spec(\cB)$ for an Azumaya algebra $\cB$ over $\cO_S$. Then the canonical map $X_{cl}\to S_{cl}$
is homeomorphism. 
\end{lem}

\begin{pf} We may assume that $S=\Spec(R)$ for a finitely generated commutative $\C$-algebra $R$ and $X=\Spec(A)$ for an Azumaya algebra $A$ over $R$. Then $\max(R^{\cont}\otimes_R A)\to \max(R^{\cont})$ is a homeomorphism for the Zariski topology because $R^{\cont}\otimes_R A$ is an Azumaya algebra over $R^{\cont}$. 
\end{pf}

\begin{para}\label{antop} For a scheme $S$  over $\C$ of finite type, we have the sheaf  $\cO_S^{\an}$ on $S(\C)$. This is  generalized to a non-commutative version as follows. 

For a non-commutative scheme $X$ over $\C$ satisfying ({\bf F}$_{\C}$), we have a sheaf  $\cO_X^{\an}$ of rings on $X_{\cl}$ defined as follows. We will show that we have a non-commutative spaces of prime ideal $X^{\an}=(X_{\cl}, \cO_X^{\an}, (\frak p(x))_x)$.

Let $X$ be a non-commutative scheme over $\C$ satisfying ({\bf F}$_{\C}$). For an open set $U$ of $X_{\cl}$, let $\cF(U):=\prod_{x\in U}\; \hat \cO_{X,x}$ where $\hat \cO_{X,x}$ denotes the completion of $\cO_{X,x}$ by the ideal generated by the maximal ideal of $Z(\cO_{X,x})$. Then $\cO_X^{\an}$ is the   subsheaf of $\cF$ characterized as follows. If $A$ is a $\C$-algebra such that $Z(A)$ is a finitely generated $\C$-algebra and $A$ is a finitely generated $Z(A)$-module, and if $U$ is an open set of $X$ and if $U\to \Spec(A)$ is an open immersion over $\C$, the restriction of $\cO_X^{\an}$  to $U_{\cl}$ coincides with the image of the canonical injection $f^{-1}(\cO_S^{\an} \otimes_{Z(A)} A)\to \cF$, where $S=\Spec(Z(A))$ and $f$ is the canonical map $U_{\cl}\to S(\C)$. 
For $x\in X_{\cl}$, the stalk $\cO^{\an}_{X,x}$ is identified with $\cO^{\an}_{S,s}\otimes_{Z(\cO_{X,x})} \cO_{X,x}$ where $s$ is the image of $x$ in $S(\C)$ and its center is identified with $\cO^{\an}_{S,s}$.

 We have the non-commutative space of prime ideals  $X^{\an}=(X_{\cl}, \cO^{\an}_X, (\frak p(x))_x)$, where $\frak p(x)$ for $x\in X_{\cl}$ is the kernel of 
 the composite map $Z(\cO^{\an}_{X,x}) \otimes_{Z(\cO_{X,x})} \cO_{X,x}\to \C \otimes_{Z(\cO_{X,x})} \cO_{X,x} \to \cO_{X,x}/\frak p(x)$, in which the first arrow is given by the evaluation $Z(\cO^{\an}_{X,x})\to \C$ at $x$. 
We have a functor $X\mapsto X^{\an}$ from the category of non-commutative schemes over $\C$ satisfying ({\bf F}$_{\C}$) to the category of non-commutative spaces of prime ideals.

\end{para}

\begin{thm}\label{lprcl} Let $X$, $Y$, $Z$ be non-commutative schemes over $\C$ satisfying ({\bf F}$_{\C}$), and assume that we are given morphisms $Y\to X$ and $Z\to X$ over $\C$. Then   
the canonical map 
$$(Y\otimes _X Z)_{\cl} \to Y_{\cl} \times_{X_{\cl}} Z_{\cl}$$ 
is a homeomorphism.

\end{thm}

\begin{pf} First we consider the case  $X=\Spec(\C)$. Write $Y\otimes_{\Spec(\C)} Z$ by $Y\otimes_\C Z$. We prove that the canonical map $(Y \otimes_\C Z)_{\cl} \to Y_{\cl} \times Z_{\cl}$ is a homeomorphism. We are reduced to the case where $Y=\Spec(B)$, $Z=\Spec(C)$ for $\C$-algebras.  ($B$ and $C$ are $\C$-algebras whose centers are finitely generated $\C$-algebras and which are finitely generated as modules over the centers.)  Then this canonical map is bijective because we have the converse map $Y_{\cl} \times Z_{\cl} \to (Y \otimes_\C Z)_{\cl}= (B\otimes_\C C)_{\cl}$ which sends a pair $(\frak p, \frak q)$ for $\frak p\in Y_{\cl}$ and $\frak q\in Z_{\cl}$ to the kernel of $B\otimes_\C C \to B/\frak p \otimes_{\C} C/\frak q\cong M_m(\C)\otimes_{\C} M_n(\C)\cong M_{mn}(\C)$. 

We show that this bijection is a homeomorphism. Note that $B \otimes_\C C$ is a finitely generated module over $Z(B)\otimes_\C Z(C)$.  Let $S=\Spec(Z(B))$, $T=\Spec(Z(C))$. 
Let  $R$ be the image of $Z(B)\otimes_\C Z(C) \to A \otimes_\C B$ and let $P=\Spec(R)$. Then the composition $(Y\otimes_\C Z)_{\cl} \to P(\C) \overset{\subset}\to S(\C) \times T(\C)$ is identified with the product of the surjections $Y_{\cl}\to S(\C)$ and $Z_{\cl} \to T(\C)$, and hence the closed immersion $P(\C) \to S(\C) \times T(\C)$ is bijective and hence is a homeomorphism.
Let $(\frak p, \frak q)\in Y_{\cl} \times Z_{\cl}$ and let $x$ be the corresponding point of $(Y\otimes_\C Z)_{\cl}$.  Let $e$ be as in (5) of \ref{cl1} for $(Y, \frak p)$ and let $e'$ be one for $(Z, \frak q)$. Then $e\otimes e'$ is as in (5) of \ref{cl1} for $(Y \otimes_\C Z, x)$. Let $\pi: Y_{\cl} \times Z_{\cl} \to S(\C) \times T(\C)$ be the projection, and let $(s, t)= \pi(\frak p, \frak q)$. Then by (5) of \ref{cl1}, when $U$ ranges over sufficiently small neighborhoods of $s$ in $S(\C)$ and $U'$ ranges over sufficiently small neighborhoods of $t$ in $T(\C)$, $(D(e) \times D(e')) \cap \pi^{-1}(U \times U')$ form a base of neighborhoods of $(\frak p, \frak q)$ in $Y_{\cl} \times Z_{\cl}$. Denote the projection $(Y\otimes_\C Z)_{\cl}\to P(\C)$ also by $\pi$ and let $y=\pi(x)$. By (5) of \ref{cl1}, when $U$ ranges over sufficiently small neighborhoods of $y$ in $P(\C)$, $D(e\otimes e')\cap \pi^{-1}(U)$ form a base of neighborhoods of $x$ in $(Y\otimes_\C Z)_{\cl}$. We have $D(e\otimes e')=D(e) \times D(e')$. This shows the coincidence of the topology of $(Y\otimes_\C Z)_{\cl}$ and that of $Y_{\cl} \times Z_{\cl}$.

Now consider a general $X$. We may assume that  $X=\Spec(A)$, and $Y=\Spec(B)$, $Z=\Spec(C)$ for $A$-algebras $B$, $C$, respectively. We first show that the map $(Y\otimes_X Z)_{\cl} \to Y_{\cl} \times_{X_{\cl}} Z_{\cl}$ is bijective. This is a map induced  from the bijection  $(Y\otimes_\C Z)_{\cl} \to Y_{\cl}\times Z_{\cl}$. Hence it is enough to prove the surjectivity of the map in problem. Assume $\frak p\in Y_{\cl}$ and $\frak q\in Z_{\cl}$ have the same image $\frak r$ in $X_{\cl}$. Then $A/\frak r\cong M_r(\C)$ for some $r$. Since the category of $\C$-algebras and the category of $M_r(\C)$-algebras are equivalent by $E\mapsto M_r(E)$ (\ref{Az}),  the $A/\frak r$-algebra $B/\frak p$ corresponds to a $\C$-algebra $E$ and $B/\frak p\cong M_r(E)$. Since $B/\frak p\cong M_m(\C)$ for some $m$, $E$ must be isomorphic to $M_b(\C)$ for some $b$ and hence $B/\frak p\cong M_b(M_r(\C))$ for some $b$. Similarly,  $C/\frak q\cong M_c(M_r(\C))$ for some $c$. Hence $B/\frak p \otimes_{A/\frak r} C/\frak q\cong M_{bc}(M_r(\C))\cong M_{bcr}(\C)$. Let $\frak s\in (Y \otimes_X Z)_{\cl}$ be the kernel of $B\otimes_A C \to M_{bcr}(\C)$. Then the map in problem sends $\frak s$ to $(\frak p, \frak q)$.

Since the topology of $(Y\otimes_X Z)_{\cl}$ is the restriction of the topology of $(Y\otimes_\C Z)_{\cl}$ and the topology of $Y_{\cl} \times_{X_{\cl}} Z_{\cl}$ is the restriction of the topology of $Y_{\cl} \times Z_{\cl}$, this bijection $(Y\otimes_X Z)_{\cl}\to Y_{\cl} \times_{X_{\cl}} Z_{\cl}$ is a homeomorphism.
\end{pf}

\begin{prop}\label{etcl} Let $X, Y$ be non-commutative schemes over $\C$ satisfying ({\bf F}$_{\C}$), and let $f:X\to Y$ be an \'etale morphism over $\C$. Then the induced map $f_{\cl}: X_{\cl} \to Y_{\cl}$ is 
a local homeomorphism. 
\end{prop}

\begin{pf}
Let $x\in X_{cl}$, $y:=f_{cl}(x)\in Y_{cl}$. By  \ref{small3}, there are morphisms of non-commutative schemes $U\overset{\subset}\to V'\to V\overset{\subset}\to  Y$, where $V\to Y$ and $U\to V'$ are open immersions and $V'\to V$ is a morphism of the type explained below, an open immersion $U\to X\times_Y U$ over $U$, and a point $u$ of $U$ whose image under $U\to X\times_Y U \to X$ is $x$. Here $V'\to V$ is such that there are a scheme $S$ over $\C$ of finite type, a scheme $S'$ over $S$ which is \'etale over $S$, a morphism $V\to S$ over $\C$, and $V'=V\times_S S$. By \ref{lprcl}, $(V\times_S S')_{cl}= V_{\cl}\times_{S(\C)} S'(\C)$ and $(X\times_Y U)=X_{cl} \times_{Y_{cl}} U_{cl}$. The maps $V_{cl}\to Y_{cl}$, $U_{cl}\to (V')_{cl}$, $U_{cl}\to (X\times_Y U)_{cl}$ are open immersions of topological spaces. The map $S'(\C)\to S(\C)$ is a local homeomorphism.  Hence the maps $U_{cl}\to X_{cl}$ and $U_{cl}\to Y_{cl}$ are local homeomorphisms. This shows that there are an open neighborhood of $x$ in $X_{cl}$ and an open neighborhood of $y$ in $Y_{cl}$ such that $f_{cl}$ induces a homeomorphism between these open neighborhoods. 
\end{pf}

\begin{para} By \ref{etcl}, we have a morphism from the topos $\tilde X_{cl}$ of sheaves on $X_{cl}$ to the topos $\tilde X_{et}$ of sheaves on $X_{et}$. 

\end{para}

The following is an analogue of Artin's comparison theorem (\cite{SGA4} vol. 3, Exp. XVI, 4).

\begin{thm}\label{e=B} 
Let $X$ and $Y$ be quasi-compact non-commutative schemes over $\C$ satisfying ({\bf F$_\C$}), and let $f:X\to Y$ be an l.f.p. morphism over $\C$.

Let $\cF$ be a constructible sheaf of abelian groups (resp. sets, resp. groups) on $X_{et}$.  Then we have an isomorphism $(R^mf_*\cF)_{cl}\overset{\cong}\to R^mf_{cl,*}(\cF_{cl})$ for every $m$ (resp. for $m=0$, resp. for $m=1$). Here $(\;)_{cl}$ denotes the pullback to $Y_{cl}$, $X_{cl}$, respectively.

\end{thm}
The proof of \ref{e=B} will be completed in \ref{pfe=B}, 

\begin{lem}\label{spcase3} Let the situation be as in \ref{spcase1} and assume that the center of $Z(A)$ is a finitely generated $\C$-algebra.  Via the homeomorphism $U_{cl}\to V_{cl}$(\ref{AzB}), identify a sheaf on $U_{cl}$ with the corresponding sheaf on $V_{cl}$. 

Then for a sheaf $\cF$ of abelian groups (resp. sets, resp. groups) on $U_{cl}$, the canonical morphism $\pi_{cl}^{-1}R^mh_{cl, *}\cF\to R^mg_{cl,*}\cF$ is an isomorphism for every $m$ (resp. for $m=0$, resp. for $m=1$).

\end{lem}

\begin{pf}
Let $S=\Spec(Z(A))$. Let $y\in Y_{\cl}$ and let $s=\pi_{cl}(y)\in S(\C)$ ($\pi_{cl}$ is the map $Y_{cl}\to S(\C)$ induced by $\pi:Y\to S$).  It is sufficient to prove that the map $(R^mh_{cl,*}\cF)_s\to (R^mg_{cl,*}\cF)_y$ is an isomorphism. 
Recall that by \ref{cl1} (5), a base of neighborhoods of $y$ in $Y_{cl}$ is given as follows.   Take an open neighborhood $W$ of $s$ in $S(\C)$ and  an idempotent  $e\in \cO^{an}_S(W)\otimes_{Z(A)}  A$ such that the image of $e$ in $\cO^{\an}_{S,s} \otimes_{Z(A)} A$ is an idempotent for $(\cO^{\an}_{S,s} \otimes_{Z(A)} A, \m)$ where $\m$ is the maximal two-sided ideal corresponding to $y$ (\ref{idem}).
 Let $D(e)$ be the open subset of $\pi_{cl}^{-1}(W)\subset Y_{\cl}$ consisting of all points $y'$ such that the image of $e$ in $\cO_{Y,y'}/\frak p(y')$ is not zero.
Then when $W'$ ranges over all open neighborhoods of $s$ in $W$,   $D(e)\cap \pi_{cl}^{-1}(W')$ form a base of neighborhoods of $y$ in $Y_{\cl}$. Let $V_{cl}(W')\subset V_{cl}$ be the inverse image of $W'$ under $V_{cl}\to S(\C)$ and let $U_{cl}(W')\subset U_{cl}$ be the inverse image of $W'$ under $U_{cl}\to S(\C)$, so $U_{cl}(W')$ is the inverse image of $\pi_{cl}^{-1}(W')$ in $U_{cl}$ under $U_{cl}\to Y_{cl}$. By \ref{D(e)big}, the map $U_{cl}\to Y_{cl}$ factors through $D(e)\cap \pi_{cl}^{-1}(W')$. From this, we have
$$(R^mh_{cl,*}\cF)_{y}\overset{\cong}\to \varinjlim_{W'} H^m(V_{cl}(W'), \cF), \quad (R^mg_{cl,*}\cF)_y\overset{\cong}\to \varinjlim_{W'} H^m(U_{cl}(W'), \cF).$$
By \ref{AzB}, we have $H^m(V_{cl}(W'), \cF)\overset{\cong}\to H^m(U_{cl}(W'), \cF)$. 
\end{pf}

\begin{lem}\label{spcase4} Let the situation be as in \ref{spcase3}, Then 
for a constructible sheaf $\cF$ of abelian groups (resp. sets, resp. groups) on $U_{et}$, we have an isomorphism $(R^mg_*\cF)_{cl}\overset{\cong}\to R^mg_{cl,*}(\cF_{cl})$ for every $m$ (resp. for $m=0$, resp. for $m=1$). .

\end{lem}

\begin{pf}
 By the comparison theorem of Artin, $(R^mh_*\cF)_{cl}\overset{\cong}\to R^mh_{cl.*}\cF_{cl}$ is an isomorphism. This comparison for $R^mg_*\cF$ follows from  it by  \ref{spcase2} and \ref{spcase3}. 
\end{pf}

The following \ref{group1} and \ref{group2} are preparations for the proof of \ref{e=B} (3). 

\begin{para}\label{group1} We give basic things about torsors, We work on a site $\cS$.

(1)  Let $\cG$, $\cH$ be sheaves of groups,  let $p,q:\cG\to \cH$ be  homomorphisms, and let $\cF\subset \cG$ be the equalizer of $(p,q)$. Assume that the morphism  $\nu:  \cG\to \cH\;;\; g\mapsto p(g)q(g)^{-1}$ is surjective (equivalently, we have an isomorphism of sheaves of sets $\cG/\cF \overset{\cong}\to \cH\;;\; g\mapsto p(g)q(g)^{-1}$).  Then we have the following (1.1) and (1.2). 

(1.1) The category of $\cF$-torsors  is equivalent to   the categories of pairs $(T, \iota)$ where  $T$ a $\cG$-torsor and $\iota$ is an isomorphism of $\cH$-torsors $p_*(T)\overset{\cong}\to q_*(T)$. Here $p_*(T)$ (resp. $q_*(T)$) denotes the $\cH$-torsor $\cH \times^{\cG} T$ where $\cG\to \cH$ is $p$ (resp. $q$). The equivalence is given by $T\mapsto (T', \iota)$ where $T'=\cG\times^{\cF} T$ and  $\iota$ is the evident isomorphism. 
The converse functor sends  $(T, \iota)$ to the equalizer of the two maps $T\to q_*(T)$, one is the evident map and the other is the composition $T\to p_*(T)\overset{\iota}\to q_*(T)$.

(1.2) The canonical map from $H^1(\cS, \cF)$ to the equalizer of the two maps $H^1(\cS,\cG)\to H^1(\cS, \cH)$ induced by $p$ and $q$, respectively, is surjective. 

This follows from (1.1).

(2) Let $\cF$ be a sheaf of group and let $E$ an $\cF$-torsor. Then we have a sheaf $\cF'$ of groups which we call the twist of $\cF$ by $T$ defined as follows.  It is the sheaf of morphisms $h: E \to \cF$ such that $h(fx)=fh(x)f^{-1}$ for  $f\in \cF$, $x\in E$. The group structure of $\cF'$ is as follows. For $h_1, h_2\in \cF'$, $(h_1h_2)(x)= h_1(x)h_2(x)$ ($x\in T$). 

 The category of $\cF$-torsors and the category of $\cF'$-torsors are equivalent by $T \mapsto T'$ where $T'$ is the sheaf of morphisms  $t: E\to T$ such that $t(fx)= ft(x)$ for $f\in \cF$, $x\in E$. The $\cF'$-torsor structure on $T'$ is as follows. For $h\in \cF'$ and  $t\in T'$, $(ht)(x)= h(x)t(x)$. 
This functor sends the $\cF$-torsor $E$ to the the trivial $\cF'$-torsor $\cF'$.
The converse functor sends an $\cF'$-torsor $T'$ to the $\cF$-torsor $T$ of sheaves of morphisms $s:E\to T'$ which satisfy the following condition. For $x\in E$, $f\in \cF$, $h\in \cF'$, $s(fx)=hs(x)$ if $h(x)=f^{-1}$. The action of $\cF$ on $T$ is given by $(fs)(x)=s(f^{-1}x)$ ($f\in \cF$, $x\in E$).

\end{para}

\begin{lem}\label{group2} Let $f:X\to Y$ be as in the hypothesis of \ref{e=B}. Let $\cG, \cH$  be sheaves of groups on $X_{et}$ and let $p,q: \cG\to \cH$ be homomorphisms such that $\nu: \cG\to \cH\;;\;g\mapsto p(g)q(g)^{-1}$ is surjective.
 Let  $\cF\subset \cG$ the equalizer of $(p,q)$.  (That is, we are in the situation of \ref{group1} with $\cS=X_{et}$). Then if  $(R^mf_*\cG)_{cl}\overset{\cong}\to  R^mf_{cl,*}\cG_{cl}$ and $(R^mf_*\cH)_{cl}\overset{\cong}\to R^mf_{cl,*}\cH_{cl}$ for $m=0,1$, we have $(R^1f_*\cF)_{cl}.\overset{\cong}\to R^1f_{cl,*}\cF_{cl}$. 

\end{lem}

\begin{pf} We apply \ref{group1} to the sites $\cS= (X\times_Y Y')_{et}$ for $Y'\in Y_{et}$ and for the pullbacks of $\cF$, $\cG$, $\cH$ on $\cS$. 
We prove  that the map

\medskip

(1) $(R^1f_*\cF)_y \to  (R^1f_{cl,*}\cF_{cl})_y$ 

\medskip
\noindent
is bijective for $y\in Y_{cl}$. We use the commutative diagram of exact sequences of pointed sets.
 $$\begin{matrix}  (f_*\cG)_y&\overset{\nu}\to &  (f_*\cH)_y & \to &(R^1f_*\cF)_y &\to & (R^1f_*\cG)_y \\
  \downarrow && \downarrow &&  \downarrow&&\downarrow \\
 (f_{cl,*}\cG_{cl})_y&\overset{\nu}\to & (f_{cl,*}\cH_{cl})_y& \to & (R^1f_{cl,*}\cF_{cl})_y&\to &(R^1f_{cl,*}\cG_{cl})_y\end{matrix}$$
The vertical arrows are assumed to be isomorphisms except the third arrow. Note that for an element $E$ of $(R^1f_*\cF)_y$, we have the twists of $\cF$, $\cG$, $\cH$ by $E$ and the images of $E$ in $(R^1f_*\cG)_y$ and $(R^1f_*\cH)_y$ (\ref{group1} (2)), defined on some \'etale neighborhood of $y$, and 
 we have the corresponding diagram as above for these twists. 

We prove the injectivity of (1). Let $a,b\in (R^1f_*\cF)_y$ and assume that they have the same image in $(R^1F_{cl,*}\cF_{cl})_y$. By  twisting $\cF$ by $b$, we may assume $b$ is the trivial element of $(R^1f_*\cF)_y$. Then since the image of $a$ under $(R^1f_*\cF)_y \to  (R^1f_*\cG)_y
 \overset{\cong}\to (R^1f_{cl,*}\cG_{cl})_y $ is the trivial element, the image of $a$ in $(R^1f_*\cG)_y$ is the trivial element. The rest is the  diagram chasing using the above  diagram.

We prove the subjectivity of (1).  Let $a\in (R^1f_{cl,*}\cF_{cl})_y$. The image $a_2$ of $a$ in $(R^1f_{cl,*}\cG_{cl})_y$ comes from an element $a_3$ of $ (R^1f_*\cG)_y$.
 Since $p(a_2)=q(a_2)$ in $(R^1f_{cl,*}\cH_{cl})_y$, we have $p(a_3)=q(a_3)$ in $(R^1f_*\cH)_y$.
 By (1.2) in \ref{group1}, $a_3$ comes from an element $a_4$ of $(R^1f_*\cF)_y$. By twisting $\cF$ by $a_4$,  we may assume that $a_4$ is the trivial element of $(R^1f_*\cF)_y$. Then $a_3$ and hence $a_2$ are trivial elements. Thus we are reduced to the case $a_2$ is the trivial element. The rest is the diagram chasing using the above diagram. 
\end{pf}

\begin{para}\label{pfe=B} We prove \ref{e=B}. 
The proof is similar to the proof of \ref{fin2} given in \ref{pffin2}.

We first prove the part of \ref{e=B} for sheaves of abelian groups. 

By using the exact sequence (1) in \ref{pffin2} and by the arguments as in \ref{pffin2} using this exact sequence, we are reduced to the following case: $X$ and $Y$ are prime non-commutative schemes, the image of $X$ in $Y$ is dense, and \ref{e=B} for sheaves $\cF$ of abelian groups on $X_{et}$ is known to be true if $\cF$ has supports in  some closed subset $C\neq X$ of $X$. 
By using the exact sequence (2) in \ref{pffin2} and the similar exact sequence 
$\dots \to R^mf_{cl,*}\cF_{cl}\to R^mf_{U_{cl}, *}(\cF_U)_{cl})\oplus \dots$, we are reduced to the situation \ref{spcase3}. Assume we are in the situation \ref{spcase3}, and consider the distinguished triangle (3) in \ref{pffin2}. We have an isomorphism 
$(Rf_* Rj_*(\cF_U))_{cl}\overset{\cong}\to Rf_{cl,*}(Rj_*(\cF_U))_{cl}$ because 
$(Rf_*Rj_*(\cF_U))_{cl}\cong R(f\circ j)_{cl,*}(\cF_{U, cl})\cong Rf_{cl,*}Rj_{*,cl}(\cF_{U, cl}) \cong Rf_{*,cl}(Rj_*(\cF_U))_{cl}$, where the first isomorphism is by \ref{spcase4} and the third isomorphism is by the case $X=Y$ of \ref{spcase4}. Furthermore, $Rj_*(\cF_U)$ is constructible (\ref{fin2}) and hence $(Rj_*(\cF_U))_C$ is constructible. Hence by the distinguished triangle (3) in \ref{pffin2} and by Noetherian induction, we obtain the comparison theorem for $\cF$.

The proof of the part of \ref{e=B} for sheaves of sets  is similar to that for sheaves of abelian groups by using (4), (5), (6) in \ref{pffin2}.  

The proof of  the part of \ref{e=B} for sheaves of groups is also similar to that for sheaves of abelian groups, but we replace (1), (2), (3) in \ref{pffin2} by the following (1), (2), (3), respectively.

(1) Let $C$ and $C'$ be closed subsets of $X$ and let $\cF$ be a sheaf of groups on $X_{et}$. Let $\cG:= \cF_{C} \times \cF_{C'}$, let $\cH= \cF_{C\cap C'}$, and let $p, q: \cG\to \cH$ be the two projections. Then the morphism $\cG\to \cH\;;\; g\mapsto p(g)q(g)^{-1}$ is surjective and $\cF_{C\cup C'}$ coincides with the equalizer of $(p, q)$.

(2) Let $U_i$  ($i=1,2$)  be open sets of $X$  such that $X=U_1\cup U_2$ and let $U_3=U_1\cap U_2$. Let $\cF$ be a sheaf of groups on $X_{et}$, let $\cG= j_{1,*}(\cF_{U_1})\times j_{2,*}(\cF_{U_2})$, $\cH=j_{3,*}(\cF_{U_3})$, where $j_i$ is the inclusion morphism $U_i\to X$, and let $p, q:\cG\to \cH$ be the two canonical homoorphisms. 
Then the morphism $\cG\to \cH\;;\; g\mapsto p(g)q(g)^{-1}$ is surjective and $\cF$ coincides with the equalizer of $(p, q)$.

(3) Let the situation be as in \ref{spcase3}. Let $\cF$ be a sheaf of groups on $X_{et}$, let $\cG=j_*(\cF_U)\times \cF_C$, $\cH=(j_*(\cF_U))_C$, and let $p, q: \cF\to \cG$ be the two canonical homomorphisms. Then the morphism  $\cG\to \cH\;;\; g\mapsto p(g)q(g)^{-1}$ is surjective and $\cF$ coincides with the equalizer of $(p, q)$.

For a quasi-compact non-commutative scheme $X$ over $\C$  satisfying ({\bf F}$_\C$), let  $S_X$ (resp. $S_{X,C}$ for a closed subset $C$ of $X$) be the statement that \ref{e=B} for sheaves of groups with $m=1$ is  true for  every $Y$ and every $f:X\to Y$ and every $\cF$ (resp. every $\cF$ with supports in $C$)  as in the hypothesis in \ref{e=B}.

By (1) and \ref{group2}, if
$S_{X, C}$ and $S_{X, C'}$ are true, then $S_{X, C\cup C'}$ is true. By this, the proof of $S_X$ is reduced to the following case: $X$ and $Y$ are prime non-commutative schemes, the image of $X$ in $Y$ is dense, and $S_{X,C}$  is true if $\cF$ has supports in some closed subset $C\neq X$ of $X$. 

We prove that in (2), if $S_{U_i}$ are true for $i=1,2,3$, then $S_X$ is true. By (2) and by \ref{group2}, it is sufficient to prove that  the comparison theorem is true for the sheaf of groups $\cF_i=j_{i,*}(\cF_{U_i})$ on $X_{et}$ for each $i$.  Write $\cF_{U_i}$ by $\cE$ for simplicity. Then 
$(R^1f_*\cF_i)_{cl}$ is the kernel of $(R^1(f\circ j_i)_*\cE)_{cl} \to (f_*R^1j_{i,*}\cE)_{cl}$, 
$R^1f_{cl,*}\cF_{i, cl}= R^1f_{cl,*}j_{i,cl, *}\cE_{cl}$ (we use \ref{e=B} for sheaves of sets in this identification) is the kernel of $R^1(f\circ j_i)_{cl,*}\cE_{cl}\to f_{cl,*}R^1j_{i,cl,*}\cE_{cl}$. We have   $(R^1(f\circ j_i)_*\cE)_{cl}\cong R(f\circ j_i)_{cl,*}\cE_{cl}$ by $S_{U_i}$,  $(f_*R^1j_{i,*}\cE)_{cl}\cong f_{cl,*}(R^1j_{i,*}\cE)_{cl}\cong f_{cl,*}R^1j_{i, cl,*}\cE_{cl}$ by \ref{e=B} for sheaves of sets and by the case $Y=X$ of $S_{U_i}$.

By this, we are reduced to the situation of \ref{spcase3}. Assume we are in the situation \ref{spcase3}, and consider 
 the above (3). Note that $j_*(\cF_U)$ is constructible (\ref{fin2}) and hence $(j_*(\cF_U))_C$ is constructible. By Noetherian induction, the comparison theorem is true for $\cF_C$ and $(j_*(\cF_U))_C$. We show that the comparison theorem is also true for $j_*(\cF_U)$. In fact, $R^1f_*j_*(\cF_U)$ is the kernel of $R^1(f\circ j)_*(\cF_U)\to f_*R^1j_*(\cF_U)$ and the above proof for the comparison theorem for $j_{i,*}(\cF_{U_i})$ works by replacing $S_{U_i}$ with \ref{spcase4}. By (3) and \ref{group2}, we have the comparison theorem for $\cF$. This completes the proof of \ref{e=B}.

\end{para}

From the case $Y=\Spec(\C)$ of \ref{e=B} for sheaves of sets (with $m=0$) and for sheaves of groups (with $m=1$), we obtain

\begin{prop}\label{toreq}  Let $X$ be a non-commutative scheme over $\C$ satisfying ({\bf F$_\C$}). Let $\cF$ a constructible sheaf of groups on $X_{et}$. Then the category of $\cF$-torsors on $X_{et}$ and the category of $\cF_{cl}$-torsors on $X_{cl}$ are equivalent via $T\mapsto T_{cl}$. 

\end{prop}

We give some examples of the comparison of the \'etale topology and the classical topology. 

First, we compare an algebraic theory \ref{1dim3} and an analytic theory \ref{E1}. 

\begin{para}\label{1dim3} Example. Let $R$ be a (commutative) Dedekind domain with field of fractions $F$, let $A$ be an $R$-algebra which is  finitely generated and torsion free as an $R$-module,  and assume that $A\otimes_R F$ is a central simple algebra over $F$. Let $a: \Spec(A)\to \Spec(R)$ be the canonical morphism, and let $n\geq 1$. 
Then for the \'etale topology, we have:

(1) 
$R^0a_*(\Z/n\Z)=\Z/n\Z$ and $R^1a_*(\Z/n\Z)=0$. 

(2) $H^1(R_{et},\Z/n\Z) \overset{\cong}\to H^1(A_{et}, \Z/n\Z)$.

(3) Assume  $n$ is invertible in $R$. Then  $R^2a_*(\Z/n\Z) \cong \oplus_s \;  i_{s,*}(\Z/n\Z(-1))^{\sharp(\Sig(s))-1}$ where $s$ ranges over all points of $\Spec(R)$ of codimension one, $i_s: s \to \Spec(R)$ is the inclusion morphism, and $\Sig(s)$ is the inverse image of $s$ in $\Spec(A)$. We have $R^ma_*(\Z/n\Z)=0$ for $m\geq 3$.

\begin{pf} (2) follows from (1) by the exact sequence $0\to H^1(R_{et}, R^0a_*(\Z/n\Z))\to H^1(A_{et}, \Z/n\Z)\to H^0(R_{et}, R^1a_*(\Z/n\Z))$.

To prove (1) and (3), we may assume that $R$ is a strict local discrete valuation ring. Let $s$ be the closed point of $\Spec(R)$. Then $\Sig(s)$ is the set of all maximal two-sided ideals of $A$. For $x\in \Sig(s)$, let $\cU(x)\subset X$ be the smallest  neighborhood (\ref{cU(x)})  of $x$. We have an exact sequence $0\to H_x^0(\cU(x))\to H^0(\cU(x))\to H^0(F) \to H_x^1(\cU(x))\to  H^1(\cU(x))\to H^1(F)\to \dots $ of \'etale cohomology groups, where the coefficients of the cohomology are $\Z/n\Z$. From this and from $H^m(\cU(x))=0$ for $m\geq 1$ and $H^0(\cU(x))=H^0(F)=\Z/n\Z$, we have $H^m_x(\cU(x))=0$ for $m=0,1$,   and $H_x^m(\cU(x))\cong H^{m-1}(F)$ for $m\geq 2$. If $n$ is invertible in $R$, then  $H^1(F, \Z/n\Z)=\Z/n(-1)$, $H^m(F, \Z/n\Z)=0$ for $m\geq 2$, and hence $H^m_x(\cU(x))$ is $\Z/n\Z(-1)$ if $m\geq 2$  and is $0$ if $m\geq 3$. We have $H^m_x(X)\cong H^m_x(\cU(x))$, and an exact sequence 
$0\to \oplus_{x\in \Sig(s)}  H_x^0(X) \to H^0(X) \to H^0(F) \to \oplus_{x\in \Sig(s)}  H^1_x(X)\to H^1(X) \to H^1(F) \to\dots$. (1) and (3) follow from this. 
 \end{pf}

\end{para}

\begin{para}\label{E1} Let $G$ be the group defined by generators $\alpha$, $\beta$ and the relations $\alpha^2=1$,  $\alpha\beta \alpha^{-1}=\beta^{-1}$. Let $R$ be a commutative ring in which $2$ is invertible, and let $A$ be the group ring $R[G]$. Then 
$Z(A)$ is the polynomial ring $R[T]$ with $T=\beta+\beta^{-1}$.  We have an injective ring homomorphism $$A\to M_2(R[T])\;;\; \alpha\mapsto \begin{pmatrix} 1&0\\ 0& -1\end{pmatrix}, \;\beta\mapsto \frac{1}{2}\begin{pmatrix}  T& 1 \\ T^2-4 & T\end{pmatrix}$$ over $R$ which is injective and whose image consists of all matrices ($(a_{ij})_{1\leq i, j\leq 2})$ such that $a_{21}\equiv 0 \bmod T^2-4$. Hence $A[1/(T^2-4)] \cong M_2(R[T][1/(T^2-4)])$. 

These facts about general $R$ will be used in \ref{9Ex}, but now let $R=\C$, and let $X=\Spec(A)=\Spec(\C[G])$, $S=\Spec(Z(A))=\Spec(\C[T])$. 

For $n\in \Z$, we  have 
$$H^0(X_{cl}, \Z/n\Z)=\Z/n\Z, \quad H^2(X_{cl}, \Z/n\Z)\cong (\Z/n\Z)^2, \quad H^m(X_{cl}, \Z/n\Z)=0 \;\text{for}\; m \neq 0, 2.$$

The case $n\neq 0$ of this is the consequence of \ref{1dim3}  and the comparison theorem  \ref{e=B}. In fact, for $a:X\to S$ and  for $R^ma_*$ for the \'etale topology, by \ref{1dim3},  $R^ma_*(\Z/n\Z)=0$ for $m\neq 0, 2$,  $R^0a_*(\Z/n\Z)=\Z/n\Z$, 
 and $R^2a_*(\Z/n\Z)$ is the sheaf which is zero outside $\{(T-2), (T+2)\}$ and whose  stalks at $(T+2)$ and at $(T-2)$ are $\Z/n\Z$. 
    
 We give an analytic proof of the above results on $H^m(X_{cl}, \Z/nZ)$ including the case $n=0$.    
    Note that $X_{cl}\to \C=S(\C)$ is a homeomorphism outside $\{2, -2\}\subset \C$  and that the inverse image of $2\in \C$ (resp. $-2\in \C$) in $X_{cl}$ is a two point set.  
 Let $E\subset \C\times \R$ be the union of $\{(z,0)\;|\; |z-2|\geq 1, |z+2|\geq 1\}$ and two spheres 
 $\{(z,t)\;|\; |z-2|^2+t^2=1\}$ and  $\{(z,t)\;|\; |z+2|^2+t^2=1\}$. 
We have a continuous map $\pi: E\to X_{cl}$ defined as follows. If $z\neq 2$ and $z\neq -2$, $\pi(z,t)$ is the unique point lying over $z\in \C=S(\C)$ (that is, the maximal two-sided ideal of $A$ generated by $T-z$). If $z=2$ (resp. $z=-2$), $\pi(z,t)$ (here $t=\pm 1$) is the maximal two-sided ideal of $A$ generated by $(\beta-1, \alpha-t)$ (resp. $(\beta+1, \alpha-t)$ of $A$. Then the pullback by $\pi$ induces isomorphisms $H^m(X_{cl}, \Z/n\Z)\overset{\cong}\to H^m(E, \Z/n\Z)$ for all $m$. The isomorphism $H^2(X_{cl}, \Z/n\Z) \overset{\cong}\to (\Z/n\Z)^2$ is given by the pullback to the $H^2$ of the two spheres inside $E$.

   \end{para}
 
 \begin{para}\label{Oet2} 
    We give an example announced in \ref{Oet}. 
    
    Let $G$ be the group defined by generators $\alpha, \beta$ and the relation $\alpha\beta\alpha^{-1}=\beta^{-1}$ and let $A=\C[G]/(\beta-1)^2$, $X=\Spec(A)$.
     Then $\Gamma(X, \cO_X)=A$,  and   $\Gamma(X, Z(\cO_X))=Z(A)=\C[t^{\pm 1}]$ with $t=\alpha^2$. Hence $\Gamma(X, Z(\cO_X))$    does not contain a square root of $t$. We show that $\Gamma(X_{et}, Z(\cO_{X_{et}}))$ contains a square root ot $t$. After that, we will understand in an analytic way that a square root of $t$ exists in $\Gamma(X_{cl}, Z(\cO_X^{an}))$ where $\cO_X^{an}$ is as in \ref{antop}, and then explain that the existence of a square root of $t$ in $\Gamma(X_{et}, Z(\cO_{X_{et}}))$ can be deduced also from it.

Let  $A'=\C[t^{\pm 1/2}] \otimes_{\C[t^{\pm 1}]} A$ and let $Y=\Spec(A')$. Let $H=\{1,\sig\}$ be the Galois group of $\C[t^{\pm 1/2}]$ over $\C[t^{\pm 1}]$. Then the canonical morphism $Y\to X$ is a covering in the \'etale site of $X$, $H$ acts on $Y$ over $X$, and $H\times Y \overset{\cong}\to Y \times_X Y\;;\; (h, a) \mapsto (a, ha)$. To prove that a square root of $t$ exists in 
$\Gamma(X_{et}, Z(\cO_{X_{et}}))$, it is sufficient to prove that there is a square root of $t$ in $\Gamma(Y, Z(\cO_Y))$ which is invariant under the action of $H$. 
We have $\Gamma(Y, \cO_Y)= A' \times A'$ on which $\sig$ acts as $(a,b)\mapsto (\sig(b), \sig(a))$ ($a,b\in A'$). Then $(t^{1/2}, -t^{1/2})\in A'\times A'$ is $H$-invariant and is a square root of $t$.

We next show  that $\Gamma(X_{\cl}, Z(\cO_X^{\an}))$ contains a square root of $t$. We have a homeomorphism $X_{\cl} \cong \C^\times$ in which $u\in \C^\times$ corresponds to the maximal two-sided ideal $(\alpha-u, \beta-1)$ of $A$. Let  $S=\Spec(Z(A))$. Then we have a homeomorphism $S(\C)\cong \C^\times$ in which $v\in \C^\times$ corresponds to the maximal ideal $(t-v)$, and the canonical map  $f: X_{\cl}\to S(\C)$ is identified with $\C^\times\to \C^\times\; ;\; z\mapsto z^2$. We have $\cO^{\an}_X = f^{-1}(\cO_S^{\an} \otimes_{Z(A)} A)$ and $Z(\cO_X^{\an})=f^{-1}(\cO_S^{\an})$. There is an element of $\Gamma(X_{\cl}, f^{-1}(Z(\cO^{\an}_S)))$ which has value $u$ at the point $u\in \C^\times = X_{\cl}$. This element is a square root of $t$. 

Now in general, let $X$ be a non-commutative scheme over $\C$ satisfying ({\bf F}$_\C$). Consider the exact sequence of sheaves 
$0 \to \Z/2\Z \to Z(\cO_{X_{et}})^\times \overset{2}\to Z(\cO_{X_{et}})^\times  \to 0$  on $X_{et}$ and 
$0 \to \Z/2\Z \to Z(\cO_X^{\an})^\times \overset{2}\to Z(\cO_X^{\an})^\times  \to 0$ on $X_{\cl}$. These exact sequences induce a commutative diagram
$$\begin{matrix} \Gamma(X_{et}, Z(\cO_{X_{et}}))^\times/(\Gamma(X_{et}, Z(\cO_{X_{et}}))^\times)^2&\overset{\subset}\to & H^1(X_{et}, \Z/2\Z)\\
\downarrow && \downarrow\\
\Gamma(X_{\cl}, Z(\cO^{\an}_X))^\times/(\Gamma(X_{\cl}, Z(\cO^{\an}_X))^\times)^2& \overset{\subset}\to & H^1(X_{\cl}, \Z/2\Z).
\end{matrix}$$
Since $H^1(X_{et}, \Z/2\Z)\overset{\cong}\to H^1(X_{\cl}, \Z/2\Z)$ (\ref{e=B}), this diagram shows that for $h\in \Gamma(X, Z(\cO_X))^\times$, 
a square root of $h$ exists in $ \Gamma(X_{et}, Z(\cO_{X_{et}}))$ if and only if a square root of $f$ exists in 
$\Gamma(X_{\cl}, Z(\cO^{\an}_X))$.

\end{para}

\section{Fundamental groups}\label{sec7}

\begin{para}\label{etcon} Let $X$ be a non-commutative scheme. We say $X_{et}$ is connected if the topos $\tilde X_{et}$ of sheaves on $X_{et}$ is connected. This means that the following equivalent conditions (i)--(iii) are satisfied. (i) The canonical map $\Sig\to \Gamma(X_{et}, \Sig)$ is bijective for every set $\Sig$. (ii) The canonical map $\{0,1\}\to \Gamma(X_{et}, \{0,1\})$ is bijective. (iii) $X$ is not empty, and the constant sheaf on $X_{et}$  associated to a one-point set is not a disjoint union of two non-empty subsheaves. 

These conditions (i)--(iii) are equivalent to the following condition (iv). (iv) $X$ is not empty, and there is no idempotent in $\Gamma(X_{et}, Z(\cO_{X_{et}}))$ other than $0$, $1$. In fact, since $Z(\cO_{Y, y})$ are local rings for all non-commutative schemes $Y$ and for all $y\in Y$, the sheaf of idempotents in $Z(\cO_{X_{et}})$ on $X_{et}$ is identified with the constant sheaf associated to the set $\{0,1\}$. Hence (iv) is equivalent to (ii).

If we replace $\Gamma(X_{et}, Z(\cO_{X_{et}}))$ in the above condition (iv) by $ \Gamma(X, Z(\cO_X))$, then the modified condition  becomes equivalent to the condition $X$ is connected. Hence if $X_{et}$ is connected, then $X$ is connected. Note that $\Gamma(X, Z(\cO_X))$ and $\Gamma(X_{et}, Z(\cO_{X_{et}}))$ may be different (\ref{Oet2}). However, 
the author does not know any example of a non-commutative scheme $X$ such that $X$  is connected but $X_{et}$ is not connected.

\end{para}

\begin{prop}\label{con} Assume $X$ satisfies ({\bf F}). Then  $X$ is connected if and only if  $X_{et}$ is connected. 

\end{prop}

\begin{pf} Let $\frak e$ be the constant sheaf on $X_{et}$ associated to a one-point set, and assume $\frak e=\frak e_1 \coprod \frak e_2$ for subsheaves $\frak e_i$ of $\frak e$. Then there are $U_i\in X_{et}$ ($i=1,2$) such that $U_1\coprod U_2 \to X$ is a covering and the canonical element of $\Gamma(U_{i, et}, \frak e)$ comes from   $\Gamma(U_{i, et}, \frak e_i)$ for $i=1,2$. 
Let $V_i$ be the image of $U_i$ in $X$. Then $V_i$  are open by \ref{a3}, and $V_1\cap V_2$ is empty. If $X$ is connected, we have either $V_1$ is empty or $V_2$ is empty and hence have either $U_1$ is empty or $U_2$ is empty. Hence either $\frak e_1=\frak e$ or $\frak e_2=\frak e$. 
\end{pf}

\begin{para}\label{pi1}
Assume  $X_{et}$ is connected. 

Let $k$ be a separably closed field, let $n\geq 1$, and let $a:\Spec(M_n(k))\to X$ be a morphism (we call this morphism a base point). (We assume here that such $a$ exists.) By \ref{Azet}, the topos of sheaves on $\Spec(M_n(k))_{et}$ is equivalent to the category of sets, and hence the pullback by $a$ is a point of the topos $\tilde X_{et}$ of sheaves on $X_{et}$ which we call the point $a$.

We define the fundamental group $\pi_1(X,a)$ as the pro-finite fundamental group of the connected topos $\tilde X_{et}$ with the point $a$ (\cite{Jo} 8.49). More precisely:

Let $\cC$ be the category of sheaves on $X_{et}$ which are locally constant and finite. Then 
with the pullback functor by $a$, $\cC$ is a Galois category in the sense of \cite{SGA1} Section V 5. The pro-finite group $\pi_1(X,a)$ is  the Galois group of this Galois category. If $X=\Spec(A)$, we denote $\pi_1(X,a)$ also as $\pi_1(A,a)$. 

If $(X', a')$ is another pair as $(X, a)$, a commutative diagram
$$\begin{matrix} a'& \to & X' \\
\downarrow & & \downarrow \\
a&\to & X\end{matrix}$$
induces a homomorphism $\pi_1(X', a') \to \pi_1(X,  a)$. 

If we only assume that a base point $a\to X$ exists but do not fix it, we have $\pi_1(X)$ as an object of the category of pro-finite groups in which morphisms are continuous homomorphisms considered modulo conjugacy. If we have a morphism $X'\to X$ (we assume $(X')_{et}$ is connected and has a base point which is not fixed), we have a morphism $\pi_1(X') \to \pi_1(X)$ in this category.

\end{para}

\begin{thm}\label{e=B3} Let $X$ be a non-commutative scheme over $\C$ satisfying (F$_\C$).

(1) The map from the set of all connected components of $X_{cl}$ to the set of all connected components of $X$ is bijective.

(2) The category of locally constant finite sheaves on $X_{et}$ and that on $X_{cl}$ are equivalent via the functor $\cF\mapsto \cF_{cl}$.

(3) Assume $X$ is connected and let $a\in X_{cl}$. Then $\pi_1(X,a)$ is canoniacally isomorphic to the pro-finite completion of $\pi_1(X_{cl}, a)$.

\end{thm}

\begin{pf}  We prove (1). Let $\Sig=\{0,1\}$. By \ref{con}, $\Gamma(X, \Sig) \to \Gamma(X_{et}, \Sig)$ is bijective. By  \ref{e=B} (2), $\Gamma(X_{et}, \Sig) )\to \Gamma(X_{cl}, \Sig)$ is bijective. (1) follows from these. 

We prove (3). Assume $X$ is connected. For a finite group $G$, the set of all continuous homomorphisms  (resp. all homomorphisms) from $\pi_1(X, a)$ (resp. $\pi_1(X_{cl},a)$) to $G$ is identified with the set of all isomorphisms of pairs $(T, \iota)$ where $T$ is a $G$-torsor on $X_{et}$ (resp. $X_{cl}$) and $\iota$ is a trivialization of the $G$-torsor $T_a$. Hence (2) follows from the equivalence \ref{toreq} of the categories of $G$-torsors. 

(2) follows from (1) and (3). 
\end{pf}

\begin{para}\label{Oet4}  For a scheme $S$, we have an equivalence between the category of finite \'etale schemes over $S$ and the category of locally constant finite sheaves on $S_{et}$,  which sends an object $T$ of the former category to the sheaf $\text{Mor}_S(\;,T')$ of the latter category.

We show an example of  a non-commutative scheme $X$ satisfying ({\bf F}) such that:  

(1) Let $S$ and $T$ be as above, and let $X\to S$ be a morphism and let $Y=X\times_S T\in X_{et}$. Then it can happen that the presheaf $\text{Mor}_X(\;, Y)$ on $X_{et}$ is not a sheaf.

(2) Let $\cF$ be a locally constant finite sheaf on $X_{et}$. It can happen that $\cF$ is not isomorphic to the sheaf associated to the presheaf $\text{Mor}_X(\;, Y) $ for any object $Y$ of $X_{et}$.

Let $A$, $X=\Spec(A)$ and  $S=\Spec(\C[t^{\pm 1}])$  be as in \ref{Oet2}. 

For (1), let  $T=\Spec(\C[t^{\pm 1/2}])$ and let $Y= X\times_S T$. Let $\cF$ be the presheaf $\Mor_X(\;,Y)$ and let $\tilde \cF$ be the associated sheaf. Then $\cF(X)=\text{Mor}_X(X, Y)$ is empty because $\Gamma(X, Z(\cO_X))$ does not contain a square root of $t$, but $\tilde \cF(X)$ is not empty because $\Gamma(X_{et}, Z(\cO_{X_{et}}))$ has a square root of $t$. 

For (2), fix a square root $t^{1/2}$ of $t$ in $\Gamma(X_{et}, Z(\cO_{X_{et}}))$.  Let $\cG$ be the sheaf of square roots of this $t^{1/2}$ in $Z(\cO_{X_{et}})$ on $X_{et}$.   Then $\cG$ is  isomorphic locally on $X_{et}$ to the constant sheaf associated to a  set of order $2$. However, there is no $Y\in X_{et}$ such that  $\cG$ is the sheaf associated to the presheaf $\text{Mor}_X(\; ,Y)$. Assume $Y$ exists. Since $A/(\beta-1)\cong \C[\alpha^{\pm 1}]$ and $\alpha$ is a square root of $t$ in this ring, $Y \otimes_A A/(\beta-1)$ must be isomorphic to $\Spec(\C[\alpha^{\pm 1}][T]/(T^2-\alpha))$. Take any point $y$ of $Y$ and let $x$ be the image of $y$ in $X$.  Let $A'=\C(t) \otimes_{\C[t^{\pm 1}]} \cO_{X,x}$ and let $B'= \C(t) \otimes_{\C[t^{\pm 1}]} \cO_{Y,y}$.  
 Since $Y\to X$ is flat, $A\to B$ and $A'\to B'$ are flat, and $B'$ is a free $A'$-module of rank $2$. Let $Z$ be the center of $B'$. Since $B'$ is r.c. over $A'$, $Z\otimes_{\C(t)} A'\to B'$ is surjective. Note $B'/(\beta-1)=\C(\alpha^{1/2})$.

{\bf Claim.} The ring $Z$ is a local ring. Let $Z_1$ be the residue field of $Z$. Then  $\dim_{\C(t)}(Z_1)\leq 2$.

We prove Claim. If $Z$ is not a local ring, it is a product of two non-zero rings, and hence $B'$ is also a product of two non-zero rings, but this contradicts the fact
  $B'$ divided by some two-sided nilpotent ideal is $\C(\alpha^{1/2})$. 
There is a subring $Z_0$ of $Z$ over $\C(t)$ such that $Z_0\to Z_1$ is an isomorphism. We show $1$, $\alpha$, $\beta-\beta^{-1}$ are linearly independent over $Z_0$. Assume $a+b\alpha + c(\beta-\beta^{-1})=0$ with $a,b,c\in Z_0$. Then by applying $\alpha(-)\alpha^{-1}$ to this, we have $a+b\alpha-c(\beta-\beta^{-1})=0$. Hence $c=0$ and $a+b\alpha=0$. If $b\neq 0$, $\alpha=-a/b\in Z_0$ is central, a contradiction, and hence $b=0$ and hence $a=0$. This linear independence shows  $\dim_{\C(t)}(Z_0) \cdot 3 \leq \dim_{\C(t)} B'=8$, and hence $\dim_{\C(t)}(Z_0)\leq 8/3<3$. Hence $\dim_{\C(t)}(Z_0)\leq 2$. This completes the proof of Claim. 
We have a surjection $Z\otimes_{\C(t)} \C(\alpha) \to B'/(\beta-1)= \C(\alpha^{1/2})$. This induces a surjection $Z_1 \otimes_{\C(t)} \C(\alpha) \to \C(\alpha^{1/2})$. This is impossible because $\alpha$ is a square root of $t$ and $\dim_{\C(t)}(Z_1)\leq 2$.

\end{para}

\section{Chow groups, relation to  class field theory}\label{sec8}

We fix a quasi-compact excellent scheme $T$.

We define the Chow groups  $CH_i(X)$ for non-commutative schemes $X$ over $T$ satisfying ({\bf F}$_T$) (\ref{bfF}), and show that $CH_i$ is a covariant functor for $X$ which are proper (\ref{proper1}) over $T$.  In the case $T=\Spec(\Z)$, we relate $CH_0$ to the unramified class field theory like in \cite{KS}. 

\begin{para}\label{S8.1} For $i\geq 0$, let $T_i$ be the set of points $t$ of $T$ such that the closure of $t$ in $T$ is $i$-dimensional. 

\end{para}

\begin{para}\label{S8.12}

Let $X$ be a non-commutative scheme over $T$ satisfying ({\bf F}$_T$). For $x\in X$, let $\kappa(x)=\cO_{X,x}/\frak p(x)$. Then $\kappa(x)$ is a finite-dimensional simple algebra over a commutative field (\ref{Fres}). If $t$ denotes the image of $x$ in $T$, the center of $\kappa(x)$ is a finitely generated field over the residue field $\kappa(t)$ of $t$. This is because locally on $X$, there  are an affine scheme $S=\Spec(R)$ of finite type over $T$ and an $R$-algebra $A$ such that  $R$ is the center of $A$ and such that $A$ is finitely generated as an $R$-module, and an open immersion $X\to \Spec(A)$ over $T$. 
\end{para}

\begin{para}\label{S8.13} Let $X$ be as in \ref{S8.12}. 
Let $X_i$ be the set of all points $x$ of $X$ such that  the image $t$ of $x$ in $T$ belongs to  $T_j$ and  the center of $\kappa(x)$ is of transcendence degree $k$ over $\kappa(t)$ for some $j,k\geq 0$ such that $i=j+k$. 

If $S$, $R$, $A$ and $X\to \Spec(A)$ are as in \ref{S8.12}, then for $x\in X$ with image $s$ in $S$, $x\in X_i$ if and only if $s\in S_i$. Here $S_i$ is defined by taking $X=S$, and not defined as the set of points of $S$ whose closures in $S$ are of dimension $i$. (For example, if $T=\Spec(\Z_p)$ and $S=\Spec(\Q_p)$, the unique point of $S$ belongs to $S_1$ though its closure on $S$ is of dimension $0$. But see  \ref{dim}  concerning this point. ) 
\end{para}

\begin{lem}\label{S8.2}  Let $X$ be as in \ref{S8.12}. Let $x'\in X_i$, $x\in X_j$ and assume that $x$ belongs to the closure $X'$ of $x'$ in $X$. Endow $X'$ with the structure of a prime non-commutative scheme (\ref{spr3}). Then the center $Z(\cO_{X',x})$ is an excellent local ring of dimension $i-j$. In particular, we have $i\geq j$. We have $i>j$ if $x'\neq x$.

\end{lem}

\begin{pf} By replacing $X$ by  $X'$, we may assume that $X$ is the closure of $x'$. Working locally on $X$, we may assume that there are $S$, $R$, $A$, and $X\to \Spec(A)$ as in \ref{S8.12}. By replacing $R$ by $Z(A)$, we may assume $R=Z(A)$.  Let $s'$ be the image of $x'$ in $S$ and let $s$ be the image of $x$ in $S$. Then $\cO_{S,s}=Z(\cO_{X, x})$. Since  $s'\in S_i$ and $s\in S_j$, $\cO_{S,s}$ is of dimension $i-j$. If $i=j$, we have $s'=s$ and $x$ and $x'$ correspond to maximal two-sided ideals of $\cO_{X,x}$ by \ref{fin} (1). Hence $x=x'$. 
  \end{pf}

  \begin{para} Let $X$ be as in \ref{S8.12}. For $i\geq 0$, we define 
$$CH_i(X):= \text{Coker}(\partial: \oplus_{u\in X_{i+1}} \;\kappa(u)^\times \to \oplus_{x\in X_i} \;\Z)$$
$$=  \text{Coker}(\partial: \oplus_{u\in X_{i+1}} \;K_1(\kappa(u)) \to \oplus_{x\in X_i} \;K_0(\kappa(x))),$$
where $K_1$ and $K_0$ are  the algebraic $K$-groups of finitely generated projective left modules and $\partial$ is defined as below.

Note that for a ring $A$, $K_0(A)$ is generated by the classes $[P]$ and $K_1(A)$ is generated by the classes $[P, \alpha]$, where $P$ is a finitely generated projective left $A$-modue and $\alpha$ is automorphism of $P$. We have a homomorphism $A^\times \to K_1(A)$ which sends $a\in A^\times$ to $(A, r_a)$  where $r_a(b)=ba$ ($b\in A$). For $x\in X$, we have 
$\Z\overset{\cong}\to K_0(\kappa(x))$ by sending $1\in \Z$ to the class $[E]$ of a simple $\kappa(x)$-module (any left $\kappa(x)$-module is projective, and the class $[E]$ is unique) and we have an isomorphism $\kappa(x)^\times/[\kappa(x)^\times,\kappa(x)^\times]\overset{\cong}\to K_1(\kappa(x))$.

For $x\in X_i$ and $u\in X_{i+1}$, the component $\partial_{u,x}: K_1(\kappa(u)) \to K_0(\kappa(x))$ of $\partial$ is as follows. 

If $x$ is not in the closure of $u$ in $X$, $\partial_{u,x}=0$. Assume that $x$ is contained in the closure of $u$. Let  $C$ the closure of $u$ in $X$ endowed with the prime non-commutative scheme structure (\ref{spr3}). Let $\cM$ be the category of finitely generated left $\cO_{C,x}$-modules. Since $\cO_{C,x}$ is a left (and also right) Noetherian ring, $\cM$ is an abelian category. Let $\cM_{\fin}\subset \cM$ be the subcategory of $\cM$ consisting of $\cO_{C,x}$-modules of finite length. By \ref{S8.2}, $Z(\cO_{C,x})$ is a Noetherian local ring of dimension $1$  and $\cO_{C,x}$ is a finitely generated module over $Z(\cO_{C,x})$. Hence the quotient category $\cM/\cM_{\fin}$ is equivalent to the category of finitely generated $\kappa(u)$-modules. By the localization theory in $K$-theory (Quillen \cite{Qu}), we have an exact sequence
$$\dots \to K_1(\cM/\cM_{\fin}) \overset{\partial}\to K_0(\cM_\fin) \to K_0(\cM)\to \dots$$ and isomorphisms 
$$K_1(\kappa(u))\cong K_1(\cM/\cM_{\fin}), \quad \oplus_{\m}\; K_0(\cO_{C,x}/\m)\overset{\cong}\to K_0(\cM_\fin),$$ where $\m$ ranges over all maximal two-sided ideals of $\cO_{C,x}$. 
 The map $\partial_{u,x}$ is defined to be the composition $K_1(\kappa(u))\overset{\partial}\to K_0(\cM_\fin)\cong \oplus_{\m}\; \Z \to \Z$, where the last arrow is to take the $\frak p(x)$-component. 

The above map $K_1(\kappa(u)) \overset{\partial}\to K_0(\cM_\fin)$ is described as $[P,\alpha] \mapsto [M_2/\alpha(M_1)] - [M_2/M_1]$ where $M_i$ are finitely generated left $\cO_{C,x}$-submodules of $P$ such that $\kappa(u)\otimes_{\cO_{C,x}} M_i=P$ and such that $M_2\supset M_1$ and $M_2\supset \alpha(M_1)$. 
 \end{para}

\begin{lem}\label{dim2} Let $X$ and $Y$ be non-commutative schemes over $T$ satisfying ({\bf F}$_T$),  let $f:X\to Y$ be a morphism over $T$, and let $x\in X_i$. Then $y:=f(x)$ belongs to $Y_j$ for some $j\leq i$. 
 If $j=i$, $\kappa(x)$ is a finitely generated projective left (resp. right)  $\kappa(y)$-module. 
\end{lem}

\begin{pf} Working locally on $X$ and $Y$, we may assume that there are affine schemes $S=\Spec(R)$ and $S'=\Spec(R')$ over $T$ of finite type, a morphism $S'\to S$ over $T$, an $R$-algebra $A$ which is finitely generated as an $R$-module, an $R'$-algebra $A'$ which is finitely generated as an $R'$-module, and open immersions $Y\to \Spec(A)$ and $X\to \Spec(A')$ via which the morphisms $S'\to S$ and $X\to Y$ are compatible. Let $s$ be the image of $y$ in $S$ and let $s'$ be the image of $x$ in $S'$. Then $s\in S_j$ and $s'\in S'_i$,  and the image of $s'$ under $S'\to S$ is $s$. Hence we have that $i\geq j$ and that if $i=j$, $\kappa(s')$ is a finite extension of $\kappa(s)$. Since $\kappa(y)$ is a finitely generated $\kappa(s)$-module and $\kappa(x)$ is a finitely generated $\kappa(s')$-module, the last fact shows that if $i=j$,  $\kappa(x)$  a finitely generated as a $\kappa(y)$-module. Every finitely generated projective $\kappa(y)$-module is projective. 
\end{pf}

\begin{para}\label{push}

For non-commutative schemes $X$ and $Y$ over $T$ satisfying ({\bf F}$_T$) and for a morphism $f:X\to Y$ over $T$,  and for $i, m\geq 0$, we have a homomorphism $$\oplus_{x\in X_i} K_m(\kappa(x)) \to \oplus_{y\in Y_i} K_m(\kappa(y))$$ defined as follows. For $x\in X_i$ and $y\in Y_i$, the $(x,y)$-component $K_m(\kappa(x))\to K_m(\kappa(y))$ of this map is $0$ unless $y=f(x)$. In the case $y=\kappa(x)$, $\kappa(x)$ is a finitely generated projective left $\kappa(y)$-module by \ref{dim2}, and hence any finitely generated projective left $\kappa(x)$-module  is regarded as a finitely generated projective left $\kappa(y)$-module. This defines an exact functor from the exact category of finitely generated projective left $\kappa(x)$-modules to that for $\kappa(y)$, and defines $K_m(\kappa(x)) \to K_m(\kappa(y))$. \end{para}

\begin{para}\label{proper1} Let $X$ be a non-commutative scheme over $T$ satisfying ({\bf F}$_T$). 

We say $X$ is proper over $T$ if the following conditions (i) and (ii) are satisfied. 

(i)  $X$ is l.f.p. over $T$ and quasi-compact.

(ii) Let $C$ be a one-dimensional Noetherian integral  scheme over $T$ with function field $K$, let $L$ be a finite-dimensional simple $K$-algebra, and assume that a morphism $\Spec(L)\to X$ over $T$ is given.  Then the closure of the image of  $\Spec(L)\to C \times_T X$ which is endowed with the canonical structure of a prime non-commutative scheme (\ref{spr3})  is isomorphic to $\Spec(\cB)$ over $C$ for some coherent $\cO_C$-algebra $\cB$. 
 
 If $X$ is a scheme, this definition of properness coincides with the usual one. In fact, if $X$ is  proper over $T$  in the usual sense, the above closure in $C\times_T X$ is proper and quasi-finite over $C$, and hence finite over $C$. If the above condition of properness  is satisfied, $X$ is proper over $T$ in the usual sense by the dvr criterion of properness for Noetherian schemes. 
 
We denote by $\cP_T$ the category of non-commutative schemes over $T$ which satisfy ({\bf F}$_T$) and which are proper over $T$.

 For example, if $S$ is a proper scheme over $T$ and $X=\Spec(\cB)$ for a coherent $\cO_S$-algebra on $S$, then $X\in \cP_T$. If $X$ is proper over $T$ and $Y\overset{\subset}\to X$ is an l.f.p. closed immersion, then $Y$ is proper over $T$.

\end{para}

\begin{lem}\label{dim} Let $X$ be a non-commutative scheme over $T$ satisfying ({\bf F$_T$}). Assume that either one of the following (i) and (ii) is satisfied. 

(i) $X$ is proper over $T$.

(ii) $T$ is covered by open subschemes which are Spec of Jacobson rings. (Recall that a commutative ring is called a  Jacobson ring if every prime ideal $\frak p$  is the intersection of all maximal ideals which contain $\frak p$.)  

Let $x\in X_i$. Then $i$ is  the largest integer such that  there are points $x_j$ of $X$ for $0\leq j\leq i$ with the following properties: $x_d=x$, and for $0\leq j<d$, $x_j$ belongs to the closure of $x_{j+1}$ in $X$ and $x_j\neq x_{j+1}$.

\end{lem}

\begin{pf} We may replace $X$ by the closure of $x$ with the prime non-commutative scheme structure. It is sufficient to prove that $X_{i-1}$ is not empty if $i>0$. 
Since  $X$ is a Noetherian space,  there is a closed point $y$ of $X$. We have $y\in X_j$ for some $j\leq i$. We may assume that either $j\leq i-2$ or $j=i$. 
There is are open neighborhood $U$ of $y$ in $X$, an affine scheme $S=\Spec(R)$ over $T$ of finite type, an $R$-algebra $A$ such that $R$ is the center of $A$, $A$ is finitely generated as an $R$-module, and $A$ is a prime ring, and an open immersion $U\overset{\subset}\to \Spec(A)$ over $T$.

First assume $j\leq i-2$. 
There is a non-zero element $a\in R$ such that $A[1/a]$ is an Azumaya algebra over $R[1/a]$ and such that $U\otimes_S \Spec(R[1/a]) = \Spec(A[1/a])$. Since $\dim(R)\geq 2$, there is a prime ideal $\frak p$ of height one which does not contain $a$. Take a maximal two-sided ideal $\m$ of the non-zero ring $A\otimes_R R_{\frak p}$. Then by \ref{fin} (1), the inverse image of $\frak m$ in $R_{\frak p}$ is the unique maximal ideal of $R_{\frak p}$. Let $\frak q$ be the inverse image of $\m$ in $A$. Since $A\otimes_R R_{\frak p}=A[1/a] \otimes_R R_{\frak p}$, $\frak q$ gives a point $z$ of $X$ which lies over $\frak p$. We have $z\in X_{i-1}$ by \ref{S8.2}.  

Now assume $i=j$. Then $x$ is a closed point of $X$. We prove $x\in X_0$. 

We first assume that the condition (ii) is satisfied. Since the image $s$ of $x$ in $S$  is a closed point of $S$ by \ref{fin} (1), the condition (ii) shows $s\in S_0$ and hence $x\in X_0$ by \ref{S8.2}. 
We assume the condition (i) is satisfied. 
 Let $t$ be the image of $x$ in $T$. If $t\in T_k$ with $k>0$, let $t'$ be an element of $T_{k-1}$ which belongs to the closure $T'$ of $t$ in $T$, endow  $T'$ with a reduced scheme structure, and let $C=\Spec(\cO_{T',t'})$. If $t\in T_0$ and the center $Z(\kappa(x))$ of $\kappa(x)$ is of transcendence degree $>0$ over  $\kappa(t)$, take a curve $C$ over $\kappa(t)$ whose function field is contained in $Z(\kappa(x))$. In both cases, let $D$ be the closure of the image of $x$ in $C\times_T X$. Since $X$ is proper over $T$, $D$ is isomorphic to $\Spec(\cB)$ for some coherent $\cO_C$-algebra $\cB$ on $C$. Take any closed point $w$ of $D$ and let $z$ be the image of $w$ in $X$. Then $z$ belongs to the closure of $x$ in $X$ and $z\neq x$, a contradiction. 
\end{pf}

 We show that $X\mapsto CH_i(X)$ is a covariant functor from $\cP_T$ to the category of abelian groups.

\begin{lem}\label{CHpush}  For a morphism $X\to Y$ of $\cP_T$, the following diagram is commutative, where the vertical arrows are given by \ref{push}.

$$\begin{matrix}  \oplus_{u\in X_{i+1}} \ K_1(\kappa(u))& \overset{\partial}\to & \oplus_{x\in X_i}\; K_0(\kappa(x))\\
\downarrow &&\downarrow\\
\oplus_{v\in Y_{i+1}} \;K_1(\kappa(v)) &\overset{\partial}\to& \oplus_{y\in Y_i}\; K_0(\kappa(y)).\end{matrix}$$

Hence we have the induced homomorphism $CH_i(X) \to CH_i(Y)$. 
\end{lem}

\begin{pf}
Let $y\in Y_i$,  $u\in X_{i+1}$, and consider the composite map 

\medskip

(1) $K_1(\kappa(u)) \overset{\partial}\to \oplus_x \; K_0(\kappa(x)) \to K_0(\kappa(y))$,

\medskip
\noindent
 where $x$ ranges over all elements of $X$ which lie over $y$ and belong to the closure of $u$. Let $v$ be the image of $u$ in $Y$.

 It is sufficient to prove the following Claim 1 and Claim 2.

{\bf Claim 1.}  Assume  $v\in Y_{i+1}$ and $y$ belongs to the closure of $v$ in $Y$. Then the composite map (1) coincides with the composite map $K_1(\kappa(u))\to K_1(\kappa(v)) \overset{\partial}\to K_0(\kappa(y))$.  

{\bf Claim 2.}  Assume $v=y$. Then the  composite map (1) is the zero map.

We prove Claim 1.

We may replace $X$ (reso, $Y$) by the closure of $u$ (resp. $v$) in $X$ (resp. $Y$) with the prime non-commutative structure.
We may assume that there are an affine scheme $S=\Spec(R)$ over $T$ of finite type, an $R$-algebra $A$ such that $A$  is finitely generated as an $R$-module and such that $R$ is the center of $A$, an open neighborhood $V$ of $y$ in $Y$, and an open immersion $V\to \Spec(A)$ over $T$. Let $s$ be the image of $y$ in $S$. Then   $\cO_{S,s}$ is the center of $\cO_{Y,y}$. Let $C= \Spec(\cO_{S,s})$, $Y'= Y\times_S C=\Spec(\cO_{Y,y})$, $X'=X\times_S C= X\times_Y \Spec(\cO_{Y,y})$.   
 Let $D\subset C\times_T Y$ and  $E\subset C\times_T X$ be the closures of the images of $v$ and of $u$, respectively,  endowed with the prime non-commutative scheme structures (\ref{spr3}). By the properness of $Y$ and of $X$ over $T$, we have $D=\Spec(P)$ and $E=\Spec(Q)$  for some $\cO_{S,s}$-algebras $P$ and $Q$ which are finitely generated as $\cO_{S,s}$-modules. We have $P\subset \kappa(v)$, $Q\subset \kappa(u)$, and  $P\subset Q$ in $ \kappa(u)$. 
 Since the image of $v$ in $Y'$ is dense in $Y'$ and the image of $u$ is dense  in $X'$,  the  morphisms $Y'\to C \times_T Y$ and $X'\to C\times_T X$ induce morphisms $Y'\to D$ and $X'\to E$, respectively.
 
 {\bf Claim 3.}  We have  $P=\cO_{Y,y}$ in $\kappa(v)$. 
 
 We prove Claim 3. Let $z$ be the image of $y\in Y'$ in $P$. Then $\cO_{D,z}\subset \kappa(v)$. The morphisms $Y'\to D\to Y$ induce homomorphisms $\cO_{Y, y}\to \cO_{D, z}\subset \cO_{Y, y}$ whose composition is the identity morphism.  Hence $\cO_{D,z}=\cO_{Y,y}$ in $\kappa(v)$. We prove that $P\to \cO_{D,z}$ is an isomorphism. The image of the center of $P$ in $\cO_{Y,y}$ is contained in the center of $\cO_{Y,y}$ (\ref{ZZ}) which is $\cO_{S,s}$,  and since $\cO_{S,s}$ is contained in $P$, we have that the center of $P$ is $\cO_{S,s}$. This proves $\cO_{D,z}=P$ and hence Claim 3. 
 
 {\bf Claim 4.} The morphism $X'\to E\times_D Y'$ is an isomorphism. 
 
 We prove Claim 4. The morphisms $D\to Y$ and $E\to X$ induce a morphism $E\times_D Y'\to X\times_Y Y'=X'$. The composition $X'\to E\times_D Y'\to X'$ is the identity morphism. The morphism $\Mor(\;, E\times_D Y')\to \Mor(\;, X')$ of functors on the category of non-commutative schemes is injective because
   $\Mor(\;, E)\subset \Mor(\;, C)\times \Mor(\;, X)$ and hence  $\Mor(\;,E)\subset \Mor(\;, Y') \times \Mor(\;, X)$. This proves Claim 4.
  
  We complete the proof of Claim 1. By Claim 4, the set of points of $X$ lying over $y$ is identified with the set of two-sided maximal ideals of $Q$. By the localization theory in $K$-theory, the commutative diagram of categories $$\begin{matrix}   \cM_{Q,\fin} &\to& \cM_Q\\ \downarrow && \downarrow \\ \cM_{P,\fin} &\to& \cM_P
  \end{matrix}$$
  induces a commutative diagram 
 $$\begin{matrix}  K_1(\kappa(u))=K_1(\cM_Q/\cM_{Q,\fin}) &\to& K_0(\cM_{Q, \fin}) =\oplus_{\m\in \max(Q)} K_0(Q/\m)\\
\downarrow &&\downarrow \\
 K_1(\kappa(v))=K_1(\cM_P/\cM_{P,\fin}) &\to& K_0(\cM_{P,\fin})=\oplus_{\m\in \max(P)}\; K_0(P/\m).\end{matrix}$$
  Let $\nn$ be  the maximal two-sided ideal of $P$ corresponding to $y$, and let $\Sig$ be the set of all maximal two-sided ideals of $Q$ lying over $\nn$.
 Then the two arrows in the composition (1) are identical with the two arrows in the part $K_1(\kappa(u)) \to \oplus_{\m\in \Sig} K_0(Q/\m)\to K_0(P/{\nn})$ of this diagram, respectively,  and the map $K_1(\kappa(v))\overset{\partial}\to K_0(\kappa(y))$ is identical with the part $K_1(\kappa(v))\overset{\partial}\to K_0(P/\nn)$ of this diagram. Hence the commutativity of this diagram proves Claim 1.

We prove Claim 2.  Let $k$ be the center of $\kappa(y)$. 
We replace $T$ by $t=\Spec(k)$, $Y$  by $y$, and $X$  by $X\times_Y y$. 
Let $C$ be the proper regular curve over $T=\Spec(k)$ whose  function field is the center of $\kappa(v)$. Let $X'\subset C\times_T X$  be the closure of the image of $u$ in $C\times_T X$ endowed with the prime non-commutative scheme structure (\ref{spr3}). Since $X$ is proper over $T$, $X'=\Spec(\cB)$ for some coherent $\cO_C$-algebra $\cB$ on $C$, and $X'$ also belongs to $\cP_T$. 

Claim 2 for $(X, Y, u, y)$ is reduced to Claim 2 for $(X', Y, u, y)$  and Claim 1 for $(X', X, u, x)$ for $x\in X_0$.  Claim 1 is already proved. Claim 2 for $(X', Y, u, y)$ is reduced to Claim 2 for $(X', T, u, t)$, and this is reduced to Claim 1 for $(X', C, u, w)$ where $w$ is the generic point of $C$ and Claim 2 for $(C, T, w, t)$. This Claim 2 for $(C, T, w, t)$ is just the usual theorem for $C$ that the degree of a principal divisor is $0$. 
\end{pf}

\begin{para}\label{tpi}  Now we take $T=\Spec(\Z)$. Let $X\in \cP_\Z$. 

We  will relate $CH_0(X)$ with abelian coverings of $X$ following the  higher dimensional unramified class field theory in \cite{KS}.

Let $\pi_1^{ab}(X):=\Hom(H^1(X_{et}, \Q/\Z), \Q/\Z)$. In the case $X$ is connected (then $X_{et}$ is connected by \ref{con}), $\pi_1^{ab}(X)$ is identified with the abelianization of the fundamental group of $X$ in Section \ref{sec7}. 
Define a quotient group $\tilde \pi_1^{ab}(X)$ of $\pi_1^{ab}(X)$ as follows. 
Let $\tilde H^1(X_{et}, \Q/\Z)$ be the subgroup of  $H^1(X_{et}, \Q/\Z)$ consisting of all elements which are killed by $H^1(X_{et}, \Q/\Z) \to H^1(\Spec(M_r(\R))_{et}, \Q/\Z)= H^1(\Spec(\R)_{et}, \Q/\Z)=\Z/2\Z$  for every integer $r$ and for every morphism $\Spec(M_r(\R))\to X$. 
 Let $\tilde \pi_1^{ab}(X)=\Hom(\tilde H^1(X_{et}, \Q/\Z), \Q/\Z)$. Thus $\pi_1^{ab}(X) \to \tilde \pi_1^{ab}(X)$ is surjective and its kernel is killed by $2$.

\end{para}

\begin{para}\label{E(x)}

Let $X\in \cP_\Z$. For $x\in X_0$, we have $\kappa(x)\cong  M_{r(x)}(\F_q)$ for some $r(x)\geq 1$ and for a finite field $\F_q$ of $q$ elements. 
The inclusion morphism $x=\Spec(\kappa(x))\to X$ induces $\pi_1^{ab}(x)\to \pi_1^{ab}(X)$. Let $\varphi_x\in \pi_1^{ab}(X)$ (called the Frobenius at $x$) be the image of the canonical generator of $\pi_1(x) \cong \pi_1(\F_q)= \Gal(\bar \F_q/\F_q)$, the $q$-th power map $\bar \F_q\to \bar \F_q\;;\;x\mapsto x^q$.

The image of $\varphi_x$ in $\tilde \pi_1^{ab}(X)$ is also denoted by $\varphi_x$.
\end{para}

\begin{thm}\label{CFT}  Let $X\in \cP_\Z$. 

(1) There is a unique homomorphism $$CH_0(X) \to \tilde \pi_1^{ab}(X)  \quad (X\in \cP_\Z)$$
which sends the class of $1\in \Z$ at $x\in X_0$ to  $\varphi_x^{r(x)}$. Here $r(x)$ and  $\varphi_x$ are as in \ref{E(x)}.  

(2) For a morphism $X\to Y$ in $\cP_\Z$, we have a commutative diagram
$$\begin{matrix}  CH_0(X)&\to& \tilde \pi_1^{ab}(X)\\
\downarrow&& \downarrow\\
CH_0(Y) &\to & \tilde \pi_1^{ab}(Y).\end{matrix}$$

\end{thm}

The proof is given below.

\begin{lem}\label{clpt}  For a morphism $X\to Y$ in $\cP_\Z$, we have a commutative diagram $$\begin{matrix} \oplus_{x\in X_0} \;K_0(\kappa(x))  &\to & \pi_1^{ab}(X)\\
\downarrow && \downarrow \\
\oplus_{y\in Y_0}\; K_0(\kappa(y))&\to & \pi_1^{ab}(Y) \end{matrix}$$
\end{lem}

\begin{pf} Let $x\in X_0$, $y=f(x)\in Y_0$. Then $\kappa(y)=M_r(\F_{p^a})$ and $\kappa(x)=M_r(M_s(\F_{p^{ab}}))= M_{rs}(\F_{p^{ab}})$ for some $r,s,a,b\geq 1$. The map $K_0(\kappa(x))\to K_0(\kappa(y))$ sends the generator $[\F_{p^{ab}}^{rs}]\in K_0(\kappa(x))$ to $bs$ times the generator $[\F_{p^a}^r]\in K_0(\kappa(y))$. On the other hand, the homomorphism $\pi_1^{ab}(X) \to \pi_1^{ab}(Y)$ sends $\varphi_x$ to $\varphi_y^b$ and hence sends $\varphi_x^{rs}$ to $\varphi_y^{rsb}= (\varphi_y^r)^{bs}$. 
\end{pf}

\begin{para} We prove \ref{CFT}. 
By \ref{clpt}, for the proof of \ref{CFT},  it is sufficient to prove that the composite map $K_1(\kappa(u))\to\oplus_x \;\Z \to  \tilde \pi_1^{ab}(X)$ is the zero map for $X\in \cP_\Z$ and for $u\in  X_1$, where $x$ ranges over all points in $X_0$ which belong to the closure of $u$ in $X$. This is reduced to the classical reciprocity law of class field theory as follows.

 Let $F$ be the center of $\kappa(u)$. If $F$ is of positive characteristic, let 
 $C$ be the proper smooth curve over a finite field with function field $F$. If $F$ is of characteristic $0$, let $C= \Spec(O_F)$, where $O_F$ is the integer ring of the number field $F$. Let $X'$ be closure of the image of 
  $u \to C \times_{\Spec(\Z)} X$ endowed with the structure of a prime non-commutative scheme. Then $X'\in \cP_\Z$, $u\in X'_1$,  and by  \ref{CHpush} and \ref{clpt}, the above map $K_1(\kappa(u))\to \tilde \pi_1^{ab}(X)$ is the composition $K_1(\kappa(u))\to \tilde \pi_1^{ab}(X')\to \tilde \pi_1^{ab}(X)$. Hence we may assume $X=X'$. 
  Then $\pi_1^{ab}(X) \overset{\cong}\to \pi_1^{ab}(C)$ by \ref{1dim3} (2). We have a commutative diagram 
$$ \begin{matrix} K_1(\kappa(u)) & \to & \pi_1^{ab}(X)&\to & \tilde \pi_1^{ab}(X) \\
\downarrow && \downarrow &&\\
K_1(F) &\to & \pi_1^{ab}(C). &&\end{matrix}$$
Hence it is sufficient to prove that the composition $K_1(\kappa(u)) \to K_1(F) \to \pi_1^{ab}(C)\cong \pi_1^{ab}(X) \to \tilde \pi_1^{ab}(X)$ is the zero map. 
Let $\Sig_1$ be the set of real places of $F$ such that $F_v \otimes_F \kappa(u)\cong M_r(\R)$ for some  $r$.  Then $\tilde \pi_1^{ab}(X)$ is regarded as the quotient of $\pi_1^{ab}(C)$ by the images of $\Gal(\bar F_v/F_v)$ for $v\in \Sig_1$. Let $\Sig_2$ be the complement of $\Sig_1$ in the set of real places of $F$. By the reciprocity law of the class field theory of $F$, elements of $K_1(F)=F^\times$ which are totally positive at $\Sig_2$ are killed in $\tilde \pi_1^{ab}(X)$. All elements of the  image of  $K_1(\kappa(u))\to F^\times$ are totally positive at $\Sig_2$ because  $F_v \otimes_F \kappa(u)\cong M_r({\mathbb H})$ for some $r$ if $v\in \Sig_2$. 
\end{para}

\section{Appendix: Zeta and $L$-functions of non-commutative schemes. By Takako Fukaya and Kazuya Kato}\label{App}

 In the paper \cite{Fu}, the first author of this Appendix defined the Hasse zeta function of a finitely generated non-commutative ring over $\Z$. In this Appendix, we generalize 
it to the zeta function of a con-commutative scheme. $L$-functions for non-commutative rings were not considered in \cite{Fu}. We discuss here $L$-functions 
for non-commutative schemes using Sections \ref{sec5}, \ref{sec7}, \ref{sec8} of this paper. 
We also give a 
partial result on the relation of zeta functions of non-commutative schemes over a finite field and \'etale cohomology (Thm. \ref{Hczeta}). 
We obtain it by extending  the definition of $Rf_!$ slightly to $f$ which need not be r.c.

\begin{para}\label{Z1} Let $X$ be a quasi-compact non-commutative scheme locally  of finite presentation over $\Z$.

Let $X_0$ be the set of all points of $X$ such that $\cO_{X,x}/\frak p(x)$ is finite as a set. For $x\in X_0$, $\cO_{X,x}/\frak p(x)\cong M_r(\F_q)$ for some $r\geq 1$ and for some finite field $\F_q$ of $q$-elements. Let $N(x)=q$ for these $x$ and $q$.

We define the zeta function $\zeta(X, s)$ of $X$ as 
$$\zeta(X, s) = \prod_{x\in X_0} (1-N(x)^{-s})^{-1}.$$

If $X$ is a scheme, this coincides with the Hasse zeta function of $X$. If $X=\Spec(A)$ for a (not necessarily commutative) ring $A$, this coincides with the zeta function $\zeta_A(s)$ defined in \cite{Fu}.

\end{para}

\begin{para}\label{Z2} It can happen that $\zeta(X, s)$ diverges for all $s\in \C$. For example, this happens in the case $X=\Spec(\F_p\langle T_1, T_2\rangle)$ 
with  $\F_p\langle T_1, T_2\rangle$ the non-commutative polynomial ring in two variables over $\F_p$. If $R$ is a finitely generated commutative ring over $\Z$ and $A$ is an $R$-algebra, $\zeta_A(s)$ absolutely converges when $\text{Re}(s)\gg 0$ for example, if $A$ is finitely generated as an $R$-module (see \ref{Z10}),  or if $A$ is the group ring $R[G]$ for a group $G$ which has a finitely generated nilpotent subgroup of finite index (\cite{Fu2}).

\end{para}

\begin{rem}\label{Z3}  After the author of \cite{Fu} studied  zeta functions of non-commutative rings, she wondered whether we can have  $L$-functions for non-commutative rings. For $A=\Z$ and $B=\Z[i]$, we have the well-known decomposition $\zeta_B(s)= \zeta_A(s)L(s, \chi)$, where $L(s, \chi)$ is the unique Diriclet character of conductor $4$. 
For a homomorphism $A\to B$ of rings which need not be commutative such that $B$ is, say,  free of rank $2$ as a left $A$-module,  the hope was to decompose $\zeta_B(s)$ similarly into the product of $\zeta_A(s)$ and some $L$-function of $A$. The issue as in the following examples (1) and (2) soon appeared. 

\medskip

(1) $A=\F_p\times \F_p$ embedded in $B=M_2(\F_p)$ as the diagonal.

\medskip

(2) $A=\F_{p^2}$ embedded in $B=M_2(\F_p)$ as a subring.

\medskip

In these (1) and (2),  $\zeta_B(s)= (1-p^{-s})^{-1}$. In the case (1), $\zeta_A(s)= (1-p^{-s})^{-2}$ , and in the case (2), $\zeta_A(s)=(1-p^{-2s})^{-1}=(1-p^{-s})^{-1}(1+p^{-s})^{-1}$. Thus in both cases,  it is impossible to decompose  $\zeta_B(s)$ as the product of $\zeta_A(s)$ and  an $L$-function of $A$. 

Now we see that $B$ is not an $A$-algebra in these examples, and we do not have a morphism $\Spec(B)\to \Spec(A)$ of non-commutative schemes. 
The following Claim shows that a bad thing like in the above  (2) does not happen for a morphism of non-commutative schemes.

{\bf Claim}  Let $f: X\to Y$ be a morphism between quasi-compact non-commutative schemes locally of finite presentation over $\Z$. Let $x\in X_0$ and $y=f(x)\in Y$. Then $y\in Y_0$ and  $N(x)=N(y)^n$ for some integer $n\geq 1$.

Proof of Claim. The fact $y\in Y_0$ is clear. Assume $\cO_{Y, y}/\frak p(y)\cong M_r(\F_q)$. Since  $\cO_{X,x}/\frak p(x)$ is a $\cO_{Y,y}/\frak p(y)$-algebra, we have by \ref{Az} that  $\cO_{X,x}/\frak p(x)\cong M_r(B)$ for an $\F_q$-algebra $B$. Since $\cO_{X,x}/\frak p(x)$  is simple, $B$ must be simple and hence $B\cong M_s(\F_{q^n})$ for some $s, n\geq 1$. Hence $\cO_{X,x}/\frak p(x)\cong M_{rs}(\F_{q^n})$.

Our hope is that for a morphism $X\to Y$ as above, 
$\zeta(X, s)$ would be expressed by using  $L$-functions of $Y$, and more generally, $L$-functions of $X$ would be expressed by using  $L$-functions of $Y$. We  study this question. 

\end{rem}

\begin{para}\label{Z4} Assume $X_{et}$ is connected (\ref{etcon}). (If $X$ satisfies ({\bf F}), this is the same as $X$ is connected, by \ref{con}). Let $x\in X_0$. Then we have the Frobenius $$\varphi_x\in \pi_1(X)$$ defined modulo conjugacy as follows. Assume $\cO_{X,x}/\frak p(x)\cong M_r(\F_q)$. The morphism $\Spec(\cO_{X,x}/\frak p(x))\to X$ induces a homomorphism $\pi_1(\F_q)=\pi_1(M_r(\F_q)) \to \pi_1(X)$  of pro-finite groups modulo conjugacy (\ref{pi1}). We define $\varphi_x$ as the image of the canonical generator $\bar \F_q \to \bar \F_q\; ;\; x\mapsto x^q$ of $\pi_1(\F_q)=\Gal(\bar \F_q/\F_q)$.

 If $X\to Y$ is a morphism between quasi-compact non-commutative schemes locally of finite presentation over $\Z$ which sends $x\in X_0$ to $y\in Y_0$, $\pi_1(X) \to \pi_1(Y)$ sends 
 $\varphi_x$ to $\varphi_y^f$ where $f$ is the integer such that $N(x)=N(y)^f$. 
\end{para}

\begin{para}\label{Z6}  The following questions arise, but the authors do not know the answers.

{\bf Question.} Is there a Chebotarev density theorem for non-commutative schemes about the distribution of $\varphi_x$?

Fix an isomorphism of commutative fields $\bar \Q_{\ell}\cong \C$.

{\bf Question.} Assume $X$ is over a finite field $k$, $X$ is irreducible, $\ell$ is not the characteristic of $k$, and let $\rho: \pi_1(X) \to GL_n(\bar \Q_{\ell})$ be a continuous irreducible representation. Then is $\rho$ pure? That is, is there an integer $w$ such that the complex eigenvalues of $\rho(\varphi_x)$ are $N(x)^{w/2}$ for all $x\in X_0$?

For a finite-dimensional continuous representation $\rho$ of $\pi_1(X)$ over $\bar \Q_{\ell}$, we define the  $L$-function $L(X, \rho, s)$ of $\rho$ as
$$L(X, \rho, s):= \prod_{x\in X_0} \; \text{det}(1-\rho(\varphi_x)^{-1} N(x)^{-s})^{-1}.$$

If $\rho$ is the trivial one-dimensional representation, we have $L(X, \rho, s)=\zeta(X, s)$. 

If $\zeta(X, s)$ absolutely converges when $Re(s)\gg 0$ and if $\rho$ has bounded weights (this means that  there are $a,b\in \R$ such that for every $x\in X_0$,  all eigen values $\alpha\in \C$ of $\rho(\varphi_x)$ satisfy $N(x)^a\leq |\alpha|\leq N(x)^b$), then $L(X, \rho, s)$ absolutely converges in $\C$ when $Re(s)\gg 0$.

\end{para}

\begin{para}\label{Art0} In the following \ref{ThArt}, let $X\to Y$ be a morphism of quasi-compact non-commutative schemes which are l.p.f over $\Z$. Assume that there is a covering $Y'\to Y$ in $Y_{et}$ such that $X\times_Y Y'$ is isomorphic over $Y'$ to a disjoint union of $m$ copies of $Y'$ for an integer $m\geq 1$. Assume that $X_{et}$ and $Y_{et}$ are connected (\ref{etcon}) and that $X_0$ is non-empty, so $\pi_1(X)$ and $\pi_1(Y)$ are defined. 

Then the homomorphism $\pi_1(X) \to \pi_1(Y)$ is injective, the image of it is an open subgroup of $\pi_1(Y)$ of index $m$, and the sheaf on $Y_{et}$ associated the presheaf $\text{Mor}_Y(\;, X)$ corresponds to the $\pi_1(Y)$-set $\pi_1(Y)/\pi_1(X)$. 

Let $\ell$ be a prime number which is invertible on $Y$, and let $\rho: \pi_1(X) \to GL_n(\bar \Q_{\ell})$ be a continuous representation. 
Fixing an isomorphism $\bar \Q_{\ell}\cong \C$ of commutative fields, we assume that $\rho$ has bounded weights and $\zeta(X, s)$ absolutely converges when $Re(s)\gg 0$, so,  $L(X, \rho, s)$ absolutely converges in $\C$ when $Re(s)\gg 0$. 
Let $\tilde \rho: \pi_1(Y) \to GL_{mn}(\bar \Q_{\ell})$ be the induced representation of $\rho$. 

\end{para}

\begin{prop}\label{ThArt} Let the notation be as in \ref{Art0}. Then  $L(Y, \tilde \rho, s)$ absolutely converges when $Re(s)\gg 0$, and we have 
$$L(X, \rho,s)= L(Y, \tilde \rho, s).$$

\end{prop} 

\begin{pf}  For $y\in Y_0$, let $X_0(y)$ be the inverse image of $y$ in $X_0$. For $y\in Y_0$ and for $x\in X_0(y)$, let $\rho_x$ be the composition $\pi_1(x) \to \pi_1(X) \to GL_n(\bar \Q_{\ell})$ and let 
$\tilde \rho_y$ be the composition $ \pi_1(y) \to \pi_1(X) \to GL_{mn}(\bar \Q_{\ell})$.  We have $L(X, \rho, s)= \prod_{y\in Y_0, x\in X_0(y)} L(x, \rho_x, s)$ and $L(Y, \tilde \rho, s)=\prod_{y\in Y_0} L(y, \tilde \rho_y, s)$, and hence it is sufficient to prove $L(y, \tilde \rho_y, s)= \prod_{x\in X_0(y)} L(x, \rho_x, s)$ for each $y\in Y_0$. 

We have $y\cong \Spec(M_r(\F_y))$  for a finite field $\F_x$ of order $N(y)$ and for some $r\geq 1$. For each $x\in X_0(y)$, $x\cong \Spec(M_s(\F_x))$ for a finite field $\F_x$ of order $N(x)$ and for some $s\geq 1$, the homomorphism $\pi_1(x)= \Gal(\bar \F_x/\F_x) \to \pi_1(y)= \Gal(\bar \F_y/\F_y)$ is injective, and the $\pi_1(y)$-set $\pi_1(Y)/\pi_1(X)$ is identified with the disjoint union of the $\pi_1(y)$-sets $\pi_1(y)/\pi_1(x)$ for $x\in X_0(y)$. From this, we see that $\tilde \rho_y$ is the direct sum of the representations $\tilde \rho_x$ of $\pi_1(y)$ induced from $\rho_x$ for $x\in X_0(y)$. Hence $L(y, \tilde \rho_y)= \prod_{x\in X_0(y)} L(y, \tilde \rho_x, s)= \prod_{x\in X_0(y)} L(x, \rho_x, s)$. 
\end{pf}

\begin{para}\label{LF}  Fix an isomorphism of commutative fields $\bar \Q_{\ell}\cong \C$. Let $X$ be a quasi-compact non-commutative scheme satisfying  ({\bf F}$_\Z$), let $\ell$ be a prime number which is invertible on $X$, and let $\cF$ be a mixed $\bar \Q_{\ell}$-sheaf on $X$. 

We define the  $L$-function of $\cF$ by 
$$L(X, \cF, s)=\prod_{x\in X_0}  \text{det}(1-\varphi_x^{-1}N(x)^{-s}\;;\; \cF_{\bar x})^{-1}.$$
This absolutely converges in $\C$ when $Re(s)\gg 0$. If $X$ is connected and $\cF$ is smooth, it is the $L$-function of the corresponding representation of $\pi_1(X)$ considered in 
\ref{Z6}.

For an exact sequence $0\to \cF'\to \cF \to \cF''\to 0$ of mixed $\bar \Q_{\ell}$-sheaves, we have $L(X, \cF,s)=L(X,\cF', s)L(X, \cF'', s)$.

\end{para}

\begin{prop}\label{Z10}  Let $X$ be a quasi-compact non-commutative scheme satisfying ({\bf F}$_\Z$), let $\ell$ be a prime number which is invertible on $X$,  and let $\cF$ be a mixed $\bar \Q_{\ell}$-sheaf on $X$. Then $L(X, \cF, s)=\prod_{i=1}^n L(S_i, \cF_i, s) ^{m(i)}$ for some $n\geq 0$, for some schemes $S_i$  of finite type over $\Z$ on which $\ell$ is invertible, for some mixed $\bar \Q_{\ell}$-sheaves on $S_i$, and for some $m(i)\in\Z$. If $X$ is over a finite field $k$, we can take as $S_i$ schemes of finite type over $k$ and hence $L(X, \cF, s)$ is a rational function in $q^{-s}$. 

\end{prop}

\begin{pf} For closed subsets $C, C'$ of $X$, we have 

$L(X, \cF_{C\cup C'}, s)= L(X, \cF_C, s)L(X, \cF_{C'}, s)L(X, \cF_{C\cap C'}, s)^{-1}$. 

\noindent
By this and by Noetherian induction, we may assume that $X$ is a prime non-commutative scheme and the statement of \ref{Z10} is true if $\cF$ has supports in a closed subset $C\neq X$ of $X$. 
For open subsets $U, U'$ of $X$ such that $X=U\cup U'$, we have 

$L(X, \cF, s)= L(U, \cF, s)L(U', \cF, s)L(U\cap U', s)^{-1}$.

\noindent 
 By this, we may assume that there are a prime ring $A$ such that $Z(A)$ is a finitely generated $\Z$-algebra and $A$ is finitely generated as a $Z(A)$-module and an open immersion $X\overset{\subset}\to \Spec(A)$. There is a non-zero element $a$ of $Z(A)$ such that $A[1/a]$ is an Azumaya algebra over $Z(A)[1/a]$ and $U:=X_{\Spec(Z(A)} \Spec(Z(A)[1/a])$ coincides with $\Spec(A[1/a])$. Let $V=\Spec(Z(A)[1/a])$ and let $\cG$ be the mixed $\bar \Q_{\ell}$-sheaf on $V$ corresponding to $\cF$ via \ref{Azet}. Let $Y=X \times_{\Spec(Z(A)} \Spec(Z(A)/aZ(A))$. Then 
 
 $L(X, \cF, s)=L(U, \cF, s)L(X, \cF_Y, s)=L(V, \cG, s)L(X, \cF_Y, s)$, 
 
 \noindent
 and \ref{Z10} is true for $L(X, \cF_Y, s)$ by Noetherian induction. 
  \end{pf}

\begin{prop}\label{Z30}  Let $X$ and $Y$ be quasi-compact non-commutative schemes satisfying ({\bf F}$_\Z$) and let $X\to Y$ be a morphism 
satisfying the condition \ref{indep}.1. Let $\ell$ be a prime number which is invertible on $Y$ and let $\cF$ be a mixed $\bar \Q_{\ell}$-sheaf on $X$. Then we have 
$$L(X, \cF, s)= L(Y, Rf_!\cF, s).$$
\end{prop}
\begin{pf} We are reduced to the case $X=Y\times_S T$ with $T\to S$ a proper morphism of schemes. By \ref{pbc2}, we have $L(Y, Rf_!\cF, s)= \prod_{y\in Y_0} \; L(y, Rf^y_!(\cF|_{X(y)}), s)$ where $X(y)=X\times_Y y$, $f^y: X(y)\to y$, and $\cF|_{X(y)}$ is the inverse image of $\cF$ on $X(y)$. Hence we are reduced to the case $Y= y=\Spec(M_r(\F_q))$ for $r\geq 1$ and for a finite field $\F_q$. The morphism $Y=\Spec(M_r(\F_q))\to S$ factors as $Y \to \Spec(\F_q)\to S$. Hence by \ref{Azet},  we are reduced to the case $Y=\Spec(\F_q)=S$ and $X$ is a proper scheme over $\F_q$, and hence we are reduced to the classical formula of Grothendieck on $L$-functions of $\bar \Q_{\ell}$-sheaves.  
\end{pf}

 The $L$-function of an abelian character of the Galois group of a number field  is understood by class field theory.
  The following \ref{L=L} is a version of it.

 \begin{lem}\label{lcm} Let $X$ be a quasi-compact non-commutative scheme satisfying ({\bf F}$_\Z$). Then $r(x)$ for $x\in X_0$ are bounded.
 
 \end{lem} 
 
 \begin{pf} We may assume that $X=\Spec(A)$ for a finitely generated commutative ring $R$ over $\Z$ and for an $R$-algebra $A$ which is finitely generated as an $R$-module. 
 Then for some $n\geq 1$, $A$ is generated by $n$ elements as an $R$-module. Let $\m\in \max(A)$ and let $\m'$ be the image of $\m$ under $\max(A)\to \max(R)$. Then 
 $A/\m\cong M_{r(x)}(k)$ is generated by $n$ elements as a module over $k'=R/\m'$, where $k$ and $k'$ are finite fields such that $k\supset k'$. Hence $r(x)^2\leq n$. 
  \end{pf}

\begin{prop}\label{L=L} Let $X$ be an object of $\cP_\Z$ (\ref{proper1}) and let  $\tilde \pi_1^{ab}(X)$ as in \ref{tpi}. 

 Let $\rho$ be a  continuous homomorphism $\tilde \pi_1^{ab}(X) \to \bar \Q_{\ell}^\times$ of bounded weights.   Let  $m=l.c.m$ of $r(x)$ for $x\in X_0$ (\ref{lcm}). 
Assume that $\rho=(\rho')^m$ for some 
continous homomorphism $\rho': \tilde \pi_1^{ab}(X) \to \bar \Q_{\ell}^\times$.
  Then there is a homomorphism 
$\chi: CH_0(X) \to \bar \Q_{\ell}^\times$ 
 such that 
 $$L(X, \rho, s)= \prod_{x \in X_0} (1-\chi(x)^{-m/r(x)}N(x)^{-s})^{-1}.$$

\end{prop}

\begin{pf} Let $\chi$ be the composite map $CH_0(X) \to \tilde \pi_1^{ab}(X) \overset{\rho'}\to \bar \Q_{\ell}^\times$ (\ref{CFT}). Then $\rho(\varphi_x)^{-1}=\rho'(\varphi^m_x)^{-1}=\chi(x)^{-m/r(x)}$.  
\end{pf}

 Now we develop a theory of $Rf_!$ which is different form that in Section \ref{sec5} and which works for some non-r.c. morphisms $f$.

\begin{prop}\label{homeo}  Let $S$ be an excellent scheme, let $X$ be a non-commutative scheme, let $f: X\to S$ be a morphism, and assume that for some  coherent $\cO_S$-algebra $\cB$ on $S$, $X$ is isomorphic over $S$ to an open subspace of $\Spec(\cB)$. Then the following conditions (i)--(iii) are equivalent.

(i)  For every morphism $S'\to S$ from an excellent scheme $S'$, the map $X\times_S S'\to S'$ is a homeomorphism. 

(ii) For every \'etale morphism $S'\to S$, the map $X\times_S S'\to S'$ is a homeomorphism. 

(iii) For every $s\in S$, if $\cO_{S, \bar s}$  denote the strict henselization of the local ring $\cO_{S,s}$, there are an $\cO_{S,\bar s}$-algebra $A$ which is finitely generated as an $\cO_{S, \bar s}$-module and a closed point $\bar x$ of $\Spec(A)$ such that $X\times_S \Spec(\cO_{S, \bar s})$ is isomorphic to $\cU(\bar x)$ (\ref{cU(x)}) over $\Spec(\cO_{S, \bar s})$. 
\end{prop}

\begin{pf} 

The implication  (i) $\Rightarrow $ (ii) is clear. 

(ii) $\Rightarrow $ (iii). Let $S'=\Spec(\cO_{S, \bar s})$. By the limit argument, we see that the map $X':=X\times_S S'\to S'$  is a homeomorphism. Let $\bar x$ be the unique  point of $X'$ lying over the closed point $\bar s$ of $S'$. Let $A$ be the stalk $\cB_{\bar s}$ of $\cB$ at $\bar s$. Then $X'$ is an open set of $\Spec(A)$ containing $\bar x$. Since $S'$ is the smallest neighborhood of $\bar s$ in $S'$, $X$ must be the smallest neighborhood $\cU(\bar x)$ of $x$ in $\Spec(A)$.

(iii) $\Rightarrow$ (i). For any $s\in S$, there is only one point of $X\times_S \Spec(\cO_{S, \bar s})$ lying over the closed point $\bar s$ of $\Spec(\cO_{S, \bar s})$. This shows that for every morphism $S'\to S$ from an excellent scheme $S'$, the map $X':=X\times_S S'\to S'$ is bijective. We prove that this is a homeomorphism. 
Let $x', y'\in X'$, let $s',t'\in S'$ be the images of $x', y'$ in $S'$, respectively, and assume $t'$ converges to $s'$. Since $X'$ is locally a Noetherian space, 
it is sufficient to prove that $y'$ converges to $x'$. Let $s$ be the image of $s'$ in $S$, let $\cO_{S,\bar s}$ be the strict henselization of $\cO_{S,s}$, and let $\cO_{S', \bar s'}$ be the strict henselization of $\cO_{S',s'}$. We have a homomorphism $\cO_{S,\bar s} \to \cO_{S', \bar s'}$ over $\cO_{S,s}$ which is compatible with $S'\to S$. Since $X\times_S \Spec(\cO_{S,\bar s})$ is $D(e)$ in $\Spec(\cB_{\bar s})$ for an idempotent $e$ for $(\cB_{\bar s}, \m)$ where $\m$ is a maximal ideal of $\cB_{\bar s}$, $X'\times_{S'} \Spec(\cO_{S', s'})$ is $D(e)$ in $\Spec(\cB_{\bar s}\otimes_{\cO_{S, \bar s}} \cO_{S', \bar s'})$, and hence by \ref{idem2} (3),   all points of $X'\times_{S'} \Spec(\cO_{S', s'})$ converges to the unique closed point  $\bar x'$ of  $\Spec(\cB_{\bar s}\otimes_{\cO_{S, \bar s}} \cO_{S', \bar s'})$. Since $\Spec(\cO_{S', \bar s'})\to \Spec(\cO_{S', s'})$ is surjective, there is a point $u$ of $\Spec(\cO_{S', \bar s'})$ whose image in  $\Spec(\cO_{S',s'})$ is $t'$. Let $z$ be the point of $X'\times_{\Spec(\cO_{S', s'})} \Spec(\cO_{S', \bar s'})$ whose image in $\Spec(\cO_{S', \bar s'})$ is $u$. Then the image of $z$ in $X'$ is $y'$ because its image in $S'$ is $t'$ and $X' \to S'$ is bijective. Since $z$ converges to $\bar x'$, the image $y'$ of $z$  converges to the image $x'$ of $\bar x'$. 
\end{pf}

  \begin{para}\label{H11}  
    We will denote the equivalent conditions in \ref{homeo} by  (H$_1$). 
   
   \end{para}
   \begin{prop}\label{H12} Let $f: X\to S$ be a morphism as in \ref{homeo} satisfying the condition (H$_1$).

   (1)  The category of sheaves on $X_{et}$ and the category of sheaves on $S_{et}$ are equivalence by $f_*$ and its inverse $f^{-1}$. 
   
   (2) Let $x\in X$ and $s=f(x)\in S$, and let $\cF$ be a sheaf on $X_{et}$. Then we have $(f_*\cF)_{\bar s}\overset{\cong}\to \cF_{\bar x}$.

   (3) Over $\C$:  $X_{cl} \to S_{\cl}$ is a homeomorphism.    
 \end{prop}
 
 \begin{pf} (2) follows from the condition (iii) in \ref{homeo} by \ref{small4} (2). (1) follows from (2) by \ref{stalk}.
 (3) follows from the condition (iii) in \ref{homeo} by  (5) of \ref{cl1}.
  \end{pf}
 
 \begin{para}\label{H21}  The condition (H$_2$) for a morphism $f$ from a non-commutative scheme $X$ to an excellent  scheme $S$:
 \medskip  

(H$_2$)   \'Etale locally on $S$, there are an open covering $(U_i)_{i\in I}$ of  $X$, an open covering $(V_i)_{i\in I}$ of $S$,  and a morphism $f_i:U_i\to V_i$ which is compatible with $f$ and which satisfies (H$_1$) for each $i\in I$.

\medskip

If $f$ satisfies (H$_2$), $X\times_S S'\to S'$ satisfies (H$_2$) for every excellent scheme $S'$ and a morphism $S'\to S$. 

    \end{para}

 \begin{prop}\label{H22} 

Let $f: X\to S$ be a morphism as in \ref{H21} satisfying (H$_2$).

  (1) The inverse image functor $f^{-1}$ from the category of sheaves of abelian groups on $S_{et}$ to the category of sheaves of abelian groups  on $X_{et}$ has a left adjoint functor $f_!$.

  (2) For  $s\in S$,  we have an isomorphism  $\oplus_{\bar x} \; \cF_{\bar x} \overset{\cong}\to (f_!\cF)_{\bar s}$, where $\bar x$ ranges over all points of $X \times_S \bar s$. 
  
\end{prop}
  
 Here the map in (2) is defined as follows. By (1) of \ref{H22}, there is a canonical map $\cF \to  f^{-1}f_!\cF$ corresponding to the identity map $f_!\cF\to f_!\cF$. Hence for each $\bar x$, an element of $\cF_{\bar x}$  gives an element of $(f^{-1}f_!\cF)_{\bar x} = (f_!\cF)_{\bar s}$. This gives the canonical map in (2). 
  
  \begin{pf} Work \'etale locally on $S$, let $U_i$ and $V_i$ be as in \ref{H21} and let  $V_{ij}$ be the image of $U_{ij}:= U_i\cap U_j$. Then $U_{ij}\to V_{ij}$ also satisfies (H$_1$). Denote the morphisms $U_i \to X$, $U_{ij}\to X$, $V_i\to Y$, $V_{ij}\to Y$, $U_i\to V_i$, and $U_{ij}\to V_{ij}$ by $a_i$, $a_{ij}$, $b_i$, $b_{ij}$, $f_i$, and $f_{ij}$, respectively. For a sheaf $\cF$ of abelian groups on $X_{et}$,  define $f_!\cF$ to  be  the cokernel of $\oplus_{i,j}\;  b_{ij,!}f_{ij, *}a_{ij}^{-1}\cF\to \oplus_k\;  b_{k,!}f_{k,*}a_k^{-1}\cF$.  Here the map from the $(i, j)$-component to the $k$-component is the inclusion map if $k=i\neq j$, the minus of the
 inclusion map if $k=j\neq i$, and is the zero map otherwise. Since $f_i$ and $f_{ij}$ satisfy  (H$_1$),   $$\Hom_{\La}(f_!\cF, \cG)= \text{Ker}(\oplus_k \Hom(f_{k,*}a_k^{-1}(\cF), b_k^{-1}(\cG)) \to \oplus_{i,j} \Hom(f_{ij,*}a_{ij}^{-1}(\cF), b_{ij}^{-1}(\cG)))$$ $$= \text{Ker}(\oplus_k \Hom(a_k^{-1}(\cF), f_k^{-1}b_k^{-1}(\cG)) \to \oplus_{i,j} \Hom(a_{ij}^{-1}(\cF), f_{ij}^{-1}b_{ij}^{-1}(\cG)))= \Hom(\cF, f^{-1}\cG),$$ where the second $=$ follows from   \ref{H12} (1). This proves (1).

      (2) follows from the definition of $f_!\cF$ by \ref{H12} (2). 
     \end{pf}
  
  This $f_!$ is compatible with  base changes by a morphism $S'\to  S$ of excellent schemes. 
      
  \begin{prop}\label{Hdim1}   
    
  Let $S$ be an excellent scheme of dimension $1$ such that the strict henselization $\cO_{S,\bar s}$ of the local ring $\cO_{S,s}$ is an integral domain for every $s\in S$. Let  $X=\Spec(\cB)$ for a coherent  $\cO_S$-algebra $\cB$ on $S$ such that $\cO_S\overset{\cong}\to Z(\cB)$. Then the morphism $X\to S$ satisfies (H$_2$).
    
  \end{prop}
  
\begin{pf}  Working \'etale locally on $S$, we may assume that for every $x\in X$ with the image $s$ in $S$, the center of $\cO_{X,x}/\frak p(x)$ is a purely inseparable extension of the residue field $\kappa(s)$ of $s$. Then for some dense open set $V$ of $S$, $\cB|_V$ 
 is an Azumaya algebra over $\cO_V$. Let $U$ be the inverse image of $V$ in $X$. For each $x\in X$, let $U_x:=U\cup \{x\}$ and let $V_{x}=V\cup\{s\}$ where $s$ is the image of $x$ in $S$.  Then $(U_x)_{x\in X}$ is an open covering of $X$, $(V_x)_{x\in X}$ is an open covering of $S$, and for any \'etale morphism $S'\to V_x$ from an excellent scheme $S'$, $U_x \times_{V_x} S'\to S'$ is a homeomorphism, that is, $U_x\to V_x$ satisfies (H$_1$).  
 \end{pf}

\begin{para}\label{H3} The condition (H$_3$)  on a morphism $f:X\to S$ from a non-commutative scheme $X$ to an excellent scheme $S$.

(H$_3$)  There is a factorization $X\overset{f_1}\to Y \overset{f_2}\to T \overset{f_3}\to S$ of $f$ such that

$T$ is a scheme and $f_3$ is  compactifiable, that is, via $f_3$, $T$ is isomorphic to an open subscheme of a proper scheme over $S$,

$f_2$ satisfies (H$_2$), and

$f_1$ is an l.f.p. immersion.

Such  decomposition $f=f_3\circ f_2\circ f_1$ will be called a good factorization of $f$. 

For $f$ satisfying (H$_3$), for a sheaf of torsion abelian groups  $\cF$ on $X_{et}$. we define $Rf_!\cF$ on $S_{et}$ relative to a good factorization of $f$ as above by $$Rf_!\cF= Rf_{3,!}\circ f_{2,!}\circ f_{1,!}\cF.$$ 

The authors can not prove that this $Rf_!\cF$ is independent of the choice of a good factorization of $f$. The issue is that for another good factorization $X\to Y'\to T'\to S$ of $f$, we do not have the diagonal embedding $X\to Y\times_S Y'$ because $Y\times_S Y'$ is not defined (though $Y\otimes_S Y'$ is defined) and hence it is hard to compare two factorizations. The comparison with the classical topology in the case over $\C$ in \ref{Hcl2} below suggests that $Rf_!\cF$ may be independent of the choice of a good factorization. 

We give three remarks.

(1) Assume $f:X\to S$ satisfies (H$_3$)  and let $g:X'\to X$ be a morphism satisfying \ref{indep}.1 which is expressed as $X'\to X\times_S T'\to X$ where $T'$ is a proper scheme over $S$ and the first arrow is an l.f.p. immersion. Then $f':=f\circ g: X'\to S$ satisfies (H$_3$). In fact, we have a good factorization $f'=f'_3\circ f'_2\circ f'_1$ of $f'$ where $f'_1$ is the composition $X'\to X\times_S T'\to Y\times_S T'$, $f_2'$ is the morphism $Y\times_S T'\to T\times_S T'$, and $f'_3$ is the morphism $T\times_S T'\to S$, which are induced by $f_1, f_2, f_3$, respectively. We have  $Rf'_!=Rf_!Rg_!$ where $Rf_!$ and $Rf'_!$ are defined with respect to the good factorizations $f=f_3\circ f_2\circ f_1$ and $f'=f'_3\circ f'_2\circ f'_1$. respectively.

(2) If $X\to S$ satisfies  (H$_3$), then $f'=X\times_S S'\to S'$ satisfies (H$_3$) for every excellent scheme $S'$ with a morphism  $S'\to S$. If we define $Rf_!$ and $Rf'_!$ by using a good factorization of $f$ and the induced good factorization of $f'$, then they commute with the  pullbacks of torsion sheaves of abelian groups. 

(3) In the case $S=\Spec(k)$ for a separably closed commutative field $k$, if $f$ satisfies (H$_3$), $R^mf_!\cF$ (defined by fixing a good factorization of $f$) is denoted by $H^m_c(X_{et}, \cF)$. This will be used in \ref{Hczeta} to have  cohomological expressions of  $L$-functions.

\end{para}

\begin{prop}\label{Lf!2} Let $f:X\to S$ be a morphism from a non-commutative scheme to an excellent scheme satisfying  (H$_3$). Fix a good factorization of $f$ to define $Rf_!$. Let $\ell$ be a prime umber which is invertible on $S$. 

(1) $R^mf_!$ sends constructible sheaves of abelian groups on $X_{et}$ to constructible sheaves of abelian groups on $S_{et}$. If $m\gg 0$, $R^mf_!\cF=0$ for all sheaves of torsion abelian groups on $X_{et}$. It sends constructible $\bar \Q_{\ell}$-sheaves on $X$ to constructible $\bar \Q_{\ell}$-sheaves on $S$.

(2) If $S$ is of finite type over $\Z$, $R^mf_!$ sends mixed $\bar \Q_{\ell}$-sheaves on $X$ to mixed $\bar \Q_{\ell}$-sheaves. For a mixed $\bar \Q_{\ell}$-sheaf on $X$, we have $L(X, \cF, s)= L(S, Rf_!\cF, s)$.

\end{prop}

\begin{pf} This is reduced to the cases $f=f_1$, $f=f_2$, and $f=f_3$. The case $f=f_1$ is easy and the case $f=f_3$ is the usual theory for schemes. In the case $f=f_2$, (1) is  clear and (2) follows from  \ref{H22} (2). (The statement in (1) about $\bar \Q_{\ell}$-sheaves  follows also from \ref{lfin}.)
\end{pf}

\begin{para}\label{Hcl1} In the following, when we say locally compact, the Hausdorff condition is already included. 

To compare the above $Rf_!$ with the classical topology over $\C$, we  slightly generalize $Rf_!$ in Verdier \cite{Ve} for locally compact spaces to non-Hausdorff spaces. 

If $X$ is a non-commutative scheme over $\C$ satisfying (H$_3$), then by \ref{H12} (3), each point of $X_{cl}$ has a locally compact open neighborhood. However, $X_{cl}$ itself need not be Hausdorff.

Let $X, Y$ be topological spaces and let $f: X \to Y$ be a continuous map. Assume $Y$ is locally compact and  $X$ has a finite covering by locally compact open sets (but $X$ itself need not be Hausdorff). We define $Rf_!\cF$ for a sheaf $\cF$ of abelian groups on $X$. Write 
$X=\cup_{i\in I} U_i$ with $I$ a totally ordered finite set and with $U_i$ a locally compact open set of $X$.  Let 
$\cF\to \cI^{\bullet}$ be an injective resolution of $\cF$. 
Then $Rf_!\cF$ is defined as the simple complex associated to the double complex $(\cG^{s,t})_{s,t}$ defined as follows.  $\cG^{s,t}=0$  if $s>0$ or if $t<0$. Assume $s\leq 0$ and $t\geq 0$. For a subset $J$ of $I$,  let $U_J=\cap_{i\in J} U_i$, $f_J:U_J\to Y$ the induced map. Then $\cG^{s,t}= \oplus_J (f_J)_!(\cI^t|_{U_j})$ where $J$ ranges over all subsets of $I$ of order $1-s$ and $(f_J)_!$ is the part of compact supports of $(f_J)_*$. The differential $\cG^{s,t}\to \cG^{s, t+1}$ is induced by $d: \cI^t\to \cI^{t+1}$, and the differential $\cG^{s,t}\to \cG^{s+1,t}$ is defined in the following way. For a subset $J$ of $I$ of order $1-s$ and a subset $K$ of order $-s$, the part $(f_J)_!(\cI^t) \to (f_K)_!(\cI^t)$ is $0$ if $K$ is not contained in $J$ and $(-1)^r$ times the canonical map if $J=\{a_1, \dots, a_{1-s}\}$ ($a_1<\dots<a_{1-s}$) and $K=\{a_1, \dots, a_{r-1}, a_{r+1}, \dots, a_{1-s}\}$. 

If $X$ is locally compact, this $Rf_!\cF$ coincides with that considered  in \cite{Ve}.

In the case $Y$ is a one point set, the unique stalk of $R^mf_!\cF$ will be denoted by $H^m_c(X, \cF)$.

\end{para}

\begin{prop}\label{Hcl2}  Let $f:X\to S$ be as in \ref{H3} satisfying (H$_3$), and assume that $S$ is of finite type over $\C$ and separated. Let $\cF$ be either a constructible sheaf of abelian groups on $X_{et}$  or a constructible  $\bar \Q_{\ell}$-sheaf on $X$,  and let $\cF_{cl}$ be the pullback of $\cF$ on $X_{cl}$. Fix a good factorization of $f$. Then we have 
$$(Rf_!\cF)_{\cl}\overset{\cong}\to  Rf_{cl,!}(\cF_{cl}).$$ Here $Rf_{cl, !}$ is defined  as in  \ref{Hcl1}.

\end{prop}

\begin{pf} By the classical comparison theorem of Artin (\cite{SGA4} vol. 3, Exp. XVI, 4), we are reduced to the case $f$ satisfies (H$_2$), and this case follows from \ref{H12} (3). \end{pf}

\begin{thm}\label{Hczeta} Let $X$ be a non-commutative scheme over a finite field $k=\F_q$ such that the morphism $X\to \Spec(k)$ 
satisfying (H$_3$). Fix a good factorization of this morphism. 
 Let $\ell$ be a prime number which is different from the characteristic of $k$ and let $\cF$ be a mixed $\bar \Q_{\ell}$-sheaf on $X$. Then 
 $$L(X, \cF, s)=\prod_m \; \text{det}(1-\varphi^{-1}_qu\;|\; H_c^m((X\otimes_k \bar k)_{et}, \cF))^{(-1)^{m-1}}\quad \text{with} \;\;u=q^{-s}.$$

Here $\varphi_q$ denotes the action of the generator $\bar k \to \bar k\;;\;x\mapsto x^q$ of $\Gal(\bar k/k)$. 
\end{thm}

\begin{pf} This follows from Prop. \ref{Lf!2} (4). 
\end{pf}

\begin{para}\label{9Ex} Example. Let the group $G=\langle \alpha,\beta\;|\; \alpha^2=1, \alpha\beta\alpha^{-1}=\beta^{-1} \rangle$ be as
 in \ref{E1}.
 
 Then we have 
 $$\zeta(\Spec(\F_q[G]), s)=(1-q^{-s})^{-2}(1-q^{1-s})^{-1}$$
 for every finite field $\F_q$ of characteristic $\neq 2$. On the other hand, 
 $$H^m_c(\Spec(\C[G])_{cl}, \Z)\cong \Z^2 \;(\text{resp}. \;  \Z, \; \text{resp}. \; 0) \quad  \text{if}\; m=0 \;(\text{resp}. \; m=2, \;\text{resp.}\;m\neq 0, 2).$$
 Thus like in Weil conjectures, the zeta function of a space $\Spec(\F_q[G])$ over a finite field and the cohomology of a space $\Spec(\C[G])_{cl}$  over $\C$ 
  have similar shapes, and this is explained as follows.

 Let $R$ be a commutative ring in which $2$ is invertible. The center of $R[G]$ is $R[T]$ where $T=\beta+\beta^{-1}$. 
The morphism  $\pi: \Spec(R[G])\to \Spec(R[T])$ satisfies (H$_2$). In fact, for $1\leq i\leq 4$, let 
  $U_i=\Spec(R[G])-V(I_i)$ where 
  $$I_1= (\alpha-1, T-2), \; I_2=(\alpha+1, T-2), \; I_3= (\alpha-1, T+2), \; I_4= (\alpha+1, T+2),$$
  and let $V_i$ be the open set $D(T-2)$ (resp. $D(T+2)$) of $\Spec(R[T])$ for $i=1,2$ (resp. $i=3,4$). 
Then $\Spec(R[G])= \cup_{i=1}^4 U_i$, and  $U_i\to V_i$ satisfies (H$_1$). 
 
 By using this, for $\La=\Z/n\Z$ with an integer  $n$ which is  invertible in $R$, we have an exact sequence $0 \to \cF \to \pi_!(\La) \to \La\to 0$ on $\Spec(R[T])$, where $\cF$  restricted to  $D(T-2)\cap D(T+2)$ in $\Spec(R[T])$ is $0$ and the pullback of $\cF$ to $V(T-2)$ and that to $V(T+2)$ in $\Spec(R[T])$ are  isomorphic to $\La$. (The fact $2$ is invertible in $R$ is used for the fact $V(T-2)$ and $V(T+2)$  in $\Spec(R[T])$ are disjoint.) From this we see that $f: \Spec(R[G]) \to \Spec(R)$ satisfies (for the \'etale topology)
 $R^0f_!\La=\La^2$, $R^2f_!\La= \La(-1)$, $R^mf_!\La=0$ for $m\neq 0, 2$.  Via \ref{Hczeta} and \ref{Hcl2}, this relates the zeta function of $\Spec(\F_q[G])$ and the compact support cohomology of $\Spec(\C[G])$.

\end{para}

\begin{para} 

For a group $G$ which has a finitely generated nilpotent subgroup of finite index, it is proved in \cite{Fu2} that the zeta function $\zeta_A(s)$ of the group ring $A=\F_q[G]$
converges when $\text{Re}(s)\gg 0$. We expect that this zeta function is expressed by some compact support \'etale cohomology theory, but the method of this paper is limited to the case $A$ is a finitely generated module over its center. For example, if $G$ is the Heizenberg group (the group of unipotent upper triangular $(3,3)$-matrices with entries in $\Z$), $A$ is not a finitely generated module over its center, and we ask whether its zeta function $\zeta_A(s)= (1-q^{-s})(1-q^{s-1})^{-3}(1-q^{s-2})^3(1-q^{s-3})^{-1}$ can be explained by some compact support etale cohomology theory.
 
\end{para}


\begin{thebibliography}{99}

\bibitem{SGA4} 
    {\sc M. Artin,  A. Grothendieck, J.-L. Verdier},
    {\em S\'eminaire de G\'eom\'etrie Alg\'ebrique du Bois Marie -- 1963--64 -- Th\'eorie des topos et cohomologie \'etale des sch\'emas  (SGA 4)},  vol. 1. Lecture Notes in Math. {\bf 269} (1972),  vol. 2. Lecture Notes in Math. {\bf 270} (1972),  vol. 3. Lecture Notes in Math. {\bf  305} (1972), Springer. 
    
        
    
 \bibitem{Bo}
{\sc N. Bourbaki}, 
{\em Elements of mathematics,  Commutative algebra},
Translated from the French. Hermann, Paris; Addison-Wesley Publishing Co., Reading, Mass. (1972).

\bibitem{De}
{\sc P. Deligne},
{\em La conjecture de Weil, II},
IHES.  Publ. Math. {\bf 43}  (1980), 157--252. 

\bibitem{De2}
{\sc P. Deligne},  
{\em S\'eminaire de G\'eom\'etrie Alg\'ebrique du Bois Marie - Cohomologie \'etale  (SGA 4$\frac{1}{2}$)},  
Lecture notes in Math. {\bf 569}, Springer  (1977). 


\bibitem{Fu}
{\sc T. Fukaya},
{\em Hasse zeta functions of non-commutative rings},
J. of Algebra, {\bf 208} (1998), 304--342.

\bibitem{Fu2}
{\sc T. Fukaya},
{\em On Hasse zeta functions of group algebras of almost nilpotent groups},
Kodai Math. J. {\bf 22} (1999), 140--151.



\bibitem{EGA4.3}
{\sc A. Grothendieck, J. Dieudonn\'e}, 
{\em \'El\'ements de g\'eom\'etrie alg\'ebrique: IV. \'Etude locale des sch\'emas et des morphismes de sch\'emas, Troisi\`eme partie}, Publ. Math. IHES. {\bf  28} (1966) 5--255. 


\bibitem{SGA1} 
 {\sc  A. Grothendieck},
 {\em S\'eminaire de G\'eom\'etrie Alg\'ebrique du Bois Marie -- 1960--61 -- Rev$\hat{e}$tements \'etales et groupe fondamental  (SGA 1)}, Lecture Notes in Math. {\bf 224}, Springer  (1971). 


\bibitem{SGA5}
{\sc L. Illusie},    
{\em S\'eminaire de G\'eom\'etrie Alg\'ebrique du Bois Marie -- 1965--66 -- Cohomologie l-adique et Fonctions L  (SGA 5)}, Lecture Notes in Math. {\bf 589}, Springer  (1977). 

\bibitem{ILO}
{\sc L. Illusie, Y. Laszlo, F. Orgogozo},
{\em Travaux de Gabber},
Ast\'erisque 363--364 (2014). 


\bibitem{Jo}
{\sc P. T. Johnstone}, 
{\em Topos theory}, 
London Mathematical Society Monographs, Vol. 10. Academic Press (1977). 

\bibitem{KS}
{\sc K. Kato, S. Saito},
{\em Global class field theory of arithmetic schemes}, Contemporary Math. 55 (1986), 255--331.

\bibitem{La0}
{\sc T. Y. Lam}
A first course in noncommutative rings. Graduate Texts in Mathematics, {\bf 131}, Springer, 1991.


\bibitem{La}
{\sc T. Y. Lam},
{\em Lectures on modules and rings}, Graduate Texts in Math. {\bf 189}, Springer (1999). 


\bibitem{Po}
{\sc D. Popescu},
{\em General N\'eron desingularization and approximation}, 
Naguya Math. J. 104 (1986), 85--115. 

\bibitem{Pr}
{\sc C. Procesi},
{\em Non commutative Jacobson rings}, Ann. Scola Norm. Sup. Pisa (3) 21 (1967), 281--290. 
 
 \bibitem{Qu}
 {\sc D.  Quillen}, 
 {\em  Higher algebraic K -theory. I},  In Algebraic K -theory, Lecture Notes in Math., {\bf  341}, Springer 
 (1973), 85--147. 



\bibitem{Ve} 
{\sc J.-L. Verdier}, 
{\em Dualit\'e dans la cohomologie des espaces localement compacts}, S\'eminaire Bourbaki, Vol. 9, Paris: Soc. Math\'ematique de France,  Exp. 300,  (1966) 337--349.

\end{thebibliography}
\end{document}